\documentclass[twoside, 10pt]{report}

\usepackage[english]{babel}
\usepackage[utf8]{inputenc}
\usepackage{geometry}
\usepackage{amssymb,amsmath, amsthm, amsfonts, mathtools, mathrsfs}
\usepackage{stmaryrd}
\usepackage{enumerate}
\usepackage{graphicx}
\usepackage{color}
\usepackage{pgfplots}
\pgfplotsset{compat=1.16}
\usepackage[onehalfspacing]{setspace}
\usepackage[framemethod=tikz]{mdframed}
\usepackage{caption, subcaption}
\usepackage{pdfpages}
\allowdisplaybreaks
\usepackage[labelfont={bf},font={footnotesize},labelsep=space]{caption}

\usepackage{fancyFormat}

\mdfdefinestyle{theoremstyle}{%
linecolor=black!,
linewidth=0.5pt,%
frametitlerule=true,%
frametitlebackgroundcolor=gray!10,
innertopmargin=\topskip,%
theoremseparator={.},%
font = \itshape}

\mdtheorem[style=theoremstyle]{theorem}{Theorem}%[chapter]

\theoremstyle{definition}

\newtheorem*{remark*}	{Remark}

\DeclareMathOperator*{\esssup}{ess\,sup}

\newcommand {\D}       {\mathcal D}
\newcommand {\E}       {\mathcal E}
\newcommand {\Q}       {\mathcal Q}
\newcommand {\R}       {\mathcal R}
\newcommand {\T}       {\mathcal T}
\newcommand {\sfE}     {\mathsf  E}

\newcommand {\CC}      {\mathbb C}
\newcommand {\HH}      {\mathbb H}
\newcommand {\KK}      {\mathbb K}
\newcommand {\NN}      {\mathbb N}
\newcommand {\PP}      {\mathbb P}
\newcommand {\RR}      {\mathbb R}
\renewcommand{\SS}     {\mathbb S}
\newcommand {\VV}      {\mathbb V}
\newcommand {\WW}      {\mathbb W}
\newcommand {\XX}      {\mathbb X}
\newcommand {\YY}      {\mathbb Y}

\newcommand {\uloc}    {u^{\mathrm{loc}}_{\WW}}
\newcommand {\ut}      {\widetilde u_{\WW}}
\newcommand {\res}     [2][\WW] {\rho_{#1}(#2)}

\newcommand {\dd}      {\mathrm{\; d}}
\newcommand {\Div}     {\operatorname{div}}
\newcommand {\tr}      {\operatorname{tr}}
\newcommand {\osc}     {\operatorname{osc}}

\newcommand{\abs}      [1] {\left\vert #1 \right\vert}
\newcommand{\norm}     [1] {\left\Vert #1 \right\Vert}
\newcommand{\wt}       [1] {\widetilde{#1}}
\newcommand{\wh}       [1] {\widehat{#1}}
\newcommand{\ul}       [1] {\underline{#1}}
\newcommand{\myind}    [1] {\mathsf{#1}}         % multi-indices
\newcommand{\myvec}    [1] {\mathfrak{#1}}       % vectors
\newcommand{\mymatrix} [1] {\boldsymbol{#1}}     % matrices

\begin{document}

\setlength\abovedisplayskip{6pt}
\setlength\belowdisplayskip{5pt}

\renewcommand{\figurename}{Fig.}
\renewcommand{\tablename}{Tab.}

\pagenumbering{Roman}

% -----------------------------------------------------------------------------------------
%    Titelseite:
% -----------------------------------------------------------------------------------------

\begin{titlepage}
\begin{center}
{\Large Paris Lodron Universit\"at Salzburg}\vspace{0.1cm} \\
Fakult\"{a}t f\"{u}r Digitale und Analytische Wissenschaften \\
Fachbereich Mathematik
\vspace{4.5cm}

{\large Dissertation zur Erlangung des akademischen Grades Dr.~rer.~nat.}\vspace{0.5cm} \\
{\LARGE \bf $\boldsymbol{hp}$-FEM for Elastoplasticity \vspace{0.2cm} \\
\& $\boldsymbol{hp}$-Adaptivity Based on Local Error Reductions} \\
\vspace*{2cm}

{\bf Eingereicht von} \vspace{0.1cm}\\
{\it Patrick Bammer \,}\\ 
%{\it Patrick Bammer, MSc}\\ 
\vspace*{5cm}

{\bf Hauptbetreuer}\vspace{0.1cm} \\
{\it Univ.-Prof.~Dr.~Andreas Schr\"{o}der}\vspace{0.25cm} \\
{\bf Nebenbetreuer}\vspace{0.1cm} \\
{\it Assoz.-Prof.~Dr.~Lothar Banz} \\

\vspace*{3cm}

{\large Salzburg, J\"{a}nner 2024}

\end{center}
\end{titlepage}

% -----------------------------------------------------------------------------------------
%    Adresskopf:
% -----------------------------------------------------------------------------------------

\newpage 
\thispagestyle{empty}
\vspace*{17cm}

\quad \\
{\small
\textbf{Patrick Bammer, 01221123} \\
\textit{$hp$-FEM for Elastoplasticity \& $hp$-Adaptivity Based on Local Error Reductions} \\
Dissertation zur Erlangung des akademischen Grades Dr.~rer.~nat., J\"{a}nner 2024 \\
Hauptbetreuer: \textit{Univ.-Prof.~Dr.~Andreas Schr\"{o}der} \\
Nebenbetreuer: \textit{Assoz.-Prof.~Dr.~Lothar Banz} \vspace{0.5cm} \\
\textbf{Paris Lodron Universit\"{a}t Salzburg} \\
Fakult\"{a}t f\"{u}r Digitale und Analytische Wissenschaften\\
Fachbereich Mathematik\\
Hellbrunner Straße 34\\
5020 Salzburg
}

% -----------------------------------------------------------------------------------------
%    Abstract:
% -----------------------------------------------------------------------------------------

\chapter*{Abstract}

The \emph{finite element method} represents a certain discretization of weak formulations related to boundary value problems, which are, for instance, frequently arising in problems of solid or fluid mechanics. Thereby, the finite element solution is sought in a finite-dimensional approximation space, the definition of which is based on a decomposition of the domain associated with the boundary value problem. By enriching the underlying approximation space one can improve the numerical solution. In the case of $hp$-adaptive strategies, i.e.~applying (usually isotropic) refinements ($h$-refinements) and varying the local polynomial degree ($p$-refinements) on selected mesh elements, respectively, one can obtain notable efficient methods leading to high algebraic or even exponential convergence. In this context, usually an \emph{a posteriori error estimator} is used to steer the automatic adaptive mesh refinement. 

The first part of the present thesis consists of the papers [P1, P2, P3] and contains the numerical analysis of different $hp$-finite element discretizations related to two different weak formulations of a model problem in \emph{elastoplasticity with linearly kinematic hardening}. Thereby, the weak formulation either takes the form of a variational inequality of the second kind, including a non-differentiable plasticity functional, or represents a mixed formulation, in which the non-smooth plasticity functional is resolved by a Lagrange multiplier. As the non-differentiability of the plasticity functional causes many difficulties in the numerical analysis and the computation of a discrete solution it seems advantageous to consider discretizations of the mixed formulation. In [P2], an \emph{a priori error analysis} of an higher-order finite element discretization of the mixed formulation (explicitly including the discretization of the Lagrange multiplier) is presented. The relations between the three different $hp$-discretizations are studied in [P3] and a reliable \emph{a posteriori error estimator} that also satisfies some (local) efficiency estimates is derived. In [P1], an efficient \emph{semi-smooth Newton solver} is proposed, which is obtained by reformulating a discretization of the mixed formulation as a system of decoupled nonlinear equations. The paper [P4] represents the second part of the thesis and introduces a new \emph{$hp$-adaptive algorithm} for solving variational equations, in which the automatic mesh refinement \emph{does not} rely on the use of an a posteriori error estimator or smoothness indicators but is based on comparing locally predicted error reductions. More precisely:

\begin{itemize}
\item In [P1], the weak formulation of the model problem in form of the mixed formulation is discretized by $hp$-finite elements. Thereby, the use of biorthogonal basis functions for the discretization of the plastic strain and the Lagrange multiplier, respectively, allows to decouple the inequality constraints associated with the discrete Lagrange multiplier. Therefore, the discrete formulation can be reformulated as a system of decoupled nonlinear equations, which enables the application of various solution schemes. The numerical examples demonstrate the applicability of the proposed semi-smooth Newton solver, which is also used for the numerical experiments in the papers [P2, P3]. In particular, the examples in [P1] show the robustness to mesh size, polynomial degree and projection parameter of the solver.

\item In [P2], again a higher-order finite element method for the mixed formulation is considered. In contrast to [P1], the Frobenius norm of the discrete Lagrange multiplier is only constrained in a certain set of Gauss quadrature points instead of enforcing it in a weak sense over the entire domain. The discretization is conforming in the displacement field and the plastic strain but non-conforming in the Lagrange multiplier (except for the lowest order case). After proving a uniform discrete inf-sup constant of One and the well posedness of the discrete mixed problem the convergence and guaranteed convergence rates of the method with respect to the mesh size and polynomial degree are proved. Though the non-conformity in the Lagrange multiplier leads to an implementable discretization scheme, it causes a reduction of the guaranteed convergence rates, which, however, is common for higher-order mixed methods for variational inequalities. Indeed, optimal convergence rates are achieved for the lowest order case. Finally, numerical experiments underline the theoretical results.

\item In [P3], a reliable a posteriori error estimator is proposed, which is applicable to any discretization of the model problem that is conforming with respect to the displacement field and the plastic strain. The residual-based estimator is derived from upper and lower error estimates relying on a suitable variational equation as auxiliary problem and satisfies some (local) efficiency estimates. Beside the two $hp$-finite element discretizations of the mixed formulation, which are already introduced in the papers [P1, P2], an $hp$-finite element discretization for the variational inequality is presented. Thereby, the non-differentiable plasticity functional is approximated by an appropriate quadrature rule. Under a slight limitation of the elements' shapes all three discretizations turn out to be equivalent. Numerical experiments underline the theoretical findings and demonstrate the potential of $h$- and $hp$-adaptive finite element discretizations for problems in elastoplasticity.

\item In [P4], a new $hp$-adaptive strategy for variational equations associated with elliptic boundary value problems is introduced, which does not use classical a posteriori error estimators or smoothness indicators to steer the adaptivity. Instead, the proposed algorithm compares the predicted reduction of the energy error that can be expressed in terms of local modifications of the degrees of freedom in the underlying discrete approximation space. Thereby, the predicted error reduction can be computed by solving computationally inexpensive, low-dimensional linear problems. The concept is first presented in an abstract Hilbert space framework, before it is applied to $hp$-finite element discretizations. For the latter, an explicit construction of $p$- and $hp$-enrichment functions in any dimension associated with one element is given and a constraint coefficient technique allows an highly efficient computation of the predicted error reductions. The applicability and effectiveness of the resulting $hp$-adaptive strategy is finally illustrated with some one- and two-dimensional numerical examples.
\end{itemize}

% -----------------------------------------------------------------------------------------
%    Acknowledgements:
% -----------------------------------------------------------------------------------------

\chapter*{Acknowledgements}

First and foremost, I would like to thank my supervisor \emph{Univ.-Prof.~Dr.~Andreas Schr\"{o}der}, who sparked my interest in the field of numerical mathematics. He not only became the superviser of my master's thesis but also offered me the great opportunity to do a dissertation under his mentoring at the \emph{Paris Lodron Universit\"{a}t Salzburg}. I am very grateful for his guidance, motivation and support during the last years. I would also like to thank my co-supervisor \emph{Assoz.-Prof.~Dr.~Lothar Banz} for his support, especially, in relation to the implementation of the finite element method.
\vspace{0.2cm}

It has been a great pleasure to work with my co-author \emph{Univ.-Prof.~Dr.~Thomas P.~Wihler} from the \emph{Universit\"{a}t Bern} and I would like to express my gratitude for many interesting discussions and a fruitful collaboration. 
\vspace{-0.25cm}

I am very grateful for the support of my family during my studies and would also like to thank \emph{Let\'{i}cia} for her patience and motivation. Finally, I would like to give thanks to my friends \emph{Florian}, \emph{Miriam}, \emph{Rudolf}, \emph{Thimo}, \emph{Tobias} and, in particular, to \emph{Paolo} for his help, many pleasant conversations and for being a great colleague.

%\newpage
%\thispagestyle{plain}
%\quad
%\newpage

% -----------------------------------------------------------------------------------------
%    Inhaltsverzeichnis:
% -----------------------------------------------------------------------------------------

\newpage

\tableofcontents

%\newpage
%\thispagestyle{plain}
%\quad
%\newpage

\pagenumbering{arabic}

% -----------------------------------------------------------------------------------------
%    List of Publications:
% -----------------------------------------------------------------------------------------

\chapter{List of Publications}

This cumulative thesis consists of the papers [P1, P2, P3, P4]. The first one is already published in the book of selected papers from the ICOSAHOM conference, which took place in Vienna from 12--16 July, 2021 while the remaining papers are submitted and available as arXiv-preprints. \\

{\small
\begin{itemize}
 \item[{[P1]}] P.~Bammer, L.~Banz and A.~Schr\"{o}der, $hp$-Finite Elements with Decoupled Constraints for Elastoplasticity, published in: \emph{Spectral and High Order Methods for Partial Differential Equations ICOSAHOM 2020+1, Springer} (2023) 141–153.
 
 \item[{[P2]}] P.~Bammer, L.~Banz and A.~Schr\"{o}der, Mixed Finite Elements of Higher-Order in Elastoplasticity, submitted to: \emph{Applied Numerical Mathematics} (under review), 2024.
 
 \item[{[P3]}] P.~Bammer, L.~Banz and A.~Schr\"{o}der, A Posteriori Error Estimates for $hp$-FE Discretizations in Elastoplasticity, submitted to: \emph{Computers \& Mathematics with Applications}, 2024.
 
 \item[{[P4]}] P.~Bammer, A.~Schr\"{o}der and T.P.~Wihler, An $hp$-adaptive strategy based on locally predicted error reductions, submitted to: \emph{Computational Methods in Applied Mathematics} (in revision), 2023.
\end{itemize}
}

\newpage
\thispagestyle{plain}
\quad
\newpage

% -----------------------------------------------------------------------------------------
%    CHAPTER :  Introduction
% -----------------------------------------------------------------------------------------

\chapter{Introduction}\label{chapter:introduction}

The modeling of problems arising in physics and engineering often leads to partial differential equations on some given domain $\Omega$ with certain conditions on its boundary $\partial\Omega$. In particular, many problems of solid and fluid mechanics can be formulated as such \emph{boundary value problems}, see e.g.~\cite{ref:Bazilevs_2007, ref:Bazilevs_2011, ref:Belytschko_2000, ref:Chen_1988, ref:Ciarlet_1988, ref:Girault_2012, ref:Gwinner_2013, ref:Han_2013, ref:Kikuchi_1988, ref:Schroeder_2011}. Thereby, a \emph{classical formulation} of such problems -- in which the partial differential equation as well as the set of boundary conditions is understood to be satisfied point-wise -- asks for high smoothness requirements on the solution and the involved data for the classical formulation to make sense. These strong smoothness assumptions, however, do not guarantee the existence of a classical solution and are often unrealistic from a physical point of view as well.

% -----------------------------------------------------------------------------------------
%   Weak Formulation of Boundary Value Problems and its Approximation
% -----------------------------------------------------------------------------------------

\section{Weak Formulation of Boundary Value Problems and its Approximation}

Removing the (possibly unrealistic) high regularity requirements on the classical solution of a boundary value problem leads to a so-called \emph{weak} or \emph{variational formulation} of it, for which it is easier to show existence results for a corresponding \emph{weak solution}. Thereby, \emph{weak} refers to the lower regularity of the solution, which is usually sought in an appropriate Sobolev space. Nevertheless, the classical and the weak formulation of a boundary value problem turn out to be equivalent in the sense that any classical solution of the problem solves the weak formulation, and, conversely, any weak solution being sufficiently smooth solves the classical formulation. The weak formulation frequently takes the form of a variational equation or a variational inequality; for instance, the weak formulation of a boundary value problem in \emph{elastoplasticity with linearly kinematic hardening} can be formulated as a variational inequality of the second kind, see e.g.~\cite{ref:Carstensen_1999, ref:Chen_1988, ref:Han_1991, ref:Han_2013}. In general, the weak formulation of a boundary value problem, however, still represents an infinite-dimensional problem and, thus, it may not be possible to explicitly determine an analytic solution. For this reason, the development of numerical methods giving an approximation of the weak solution of a boundary value problem, which is frequently called a \emph{discrete solution}, is essential.
\vspace{0.1cm}

The basic idea for the numerical approximation of a variational formulation is to determine a discrete solution not in the full dimensional Hilbert space, in which the weak solution is sought, but in an appropriate subspace. If the weak formulation is projected to a \emph{finite-dimensional} subspace, the resulting finite-dimensional problem is called a \emph{Riesz-Galerkin discretization} of the original variational formulation, and the corresponding discrete solution is called a \emph{Riesz-Galerkin approximation}. In the case of a variational equation, its Riesz-Galerkin discretization turns out to be equivalent to a linear system of equations, exploiting the fact that any element of the finite-dimensional subspace can uniquely be represented in terms of a linear combination of finitely many basis functions. Thus, the applicability of efficient solution schemes for calculating the Riesz-Galerkin approximation asks for a suitable choice of basis functions for the finite-dimensional subspace. The \emph{finite element method}, see e.g.~\cite{ref:Babuska_2001, ref:Braess_2013, ref:Brenner_2008, ref:Brezzi_2012, ref:Ciarlet_2002, ref:Hackbusch_2017, ref:Karniadakis_2013, ref:Knabner_2000, ref:Schwab_1998, ref:Szabo_1991}, or, more precisely, the \emph{conforming} finite element method represents a certain Riesz-Galerkin discretization for weak formulations of boundary value problems. Thereby, the construction of the finite-dimensional subspace relies on a decomposition of the domain $\Omega$, on which the original boundary value problem is given. If the finite-dimensional approximation space is not a subspace of the full dimensional Hilbert space the method is called \emph{non-conforming}. Other prominent methods for computing an approximative solution are the \emph{finite difference method}, see e.g.~\cite{ref:Hackbusch_2017, ref:Knabner_2000, ref:Smith_1985, ref:Strikwerda_2004}, the \emph{finite volume method}, see e.g.~\cite{ref:Knabner_2000, ref:LeVeque_2002} and the \emph{boundary element method}, see e.g.~\cite{ref:Hackbusch_2012, ref:Sauter_2004, ref:Steinbach_2013, ref:Wrobel_2002}. In the field of engineering and engineering applications, however, the finite element method turned out to be the method of choice.

% -----------------------------------------------------------------------------------------
%    The Finite Element Method
% -----------------------------------------------------------------------------------------

\section{The Finite Element Method}\label{sec:FEM}

The development of a finite element algorithm for computing a discrete solution of a boundary value problem in its weak formulation follows several general steps: In a first step, the bounded domain $\Omega\subseteq \RR^d$ ($d\in\NN$) is decomposed into finitely many non-overlapping closed subdomains $K\in\D$, often called \emph{physical elements}, such as segments in the one-dimensional case, triangles and quadrilaterals in two dimensions, and tetrahedrons and hexahedrons in 3D. Then, the corresponding finite-dimensional \emph{finite element space} is defined to be the collection of functions that are piecewisely smooth (on the physical elements) and may satisfy an additional global smoothness requirement on the entire domain $\Omega$. Thereby, the global smoothness requirement depends on the regularity of the Hilbert function space, in which the weak solution is sought. If $\XX = \XX(\Omega)$ denotes the full dimensional Hilbert function space over $\Omega$, a conforming finite element space, in general, takes the form
\begin{align}\label{eq:abstract_FE_space}
 \WW 
  := \big\lbrace w\in\XX \; ; \; w_{\, | \, K}\in\Psi_K  \text{ for all } K\in\D \big\rbrace 
 \subseteq \XX,
\end{align}
where $\Psi_K$ denotes a local finite-dimensional space of smooth functions on the physical element $K\in\D$. Therefore, each $\Psi_K$ is spanned by finitely many local basis functions, which are frequently called \emph{shape functions}. In this context, a prominent choice is the finite element space consisting of piecewise polynomials (or more precisely, piecewise images of polynomials) on the physical elements which are globally continuous. The finitely many functions of a basis of the finite element space $\WW$ are usually called \emph{degrees of freedom}. In a next step, a \emph{finite element discretization} of the variational problem is obtained by projecting the weak formulation to the finite element space. In this way, a \emph{finite element system} is formed, which is finally solved, for instance, by some iterative numerical method.
\vspace{0.1cm}

The difference $u_\XX - u_\WW$ between the weak solution $u_\XX$ and its finite element approximation $u_\WW$ is called \emph{discretization error} and is usually measured in some norm $\norm{\cdot}_\XX$ on the Hilbert space $\XX$ to quantify the approximation quality of the discrete solution. The a priori error analysis for the finite element method suggests that the accuracy of the discrete solution can be improved by increasing the dimension of the underlying finite element space. Essentially, there are two strategies for such an increase:
\begin{enumerate}[(i)]
 \item Refining physical elements of the decomposition $\D$.
 \item Increasing the local approximation property on physical elements.
\end{enumerate}

It is customary, to call $\D$ a \emph{mesh} for $\Omega$ and use the lowercase $h$ as a parameter for the mesh size, for instance, $h := \max_{K\in\D} h_K$, where $h_K$ denotes the diameter of the physical element $K\in\D$. Frequently, $h$ is added as an index to $\D$, i.e.~$\D = \D_h$. Therefore, (i) is typically referred to as an \emph{$h$-(adaptive) refinement}. The local space $\Psi_K$ usually represents the space of polynomials up to degree $p_K$ on the physical element $K\in\D$ (or images of polynomials up to degree $p_K$). In this case, (ii) can be achieved by increasing local polynomial degrees $p_K$ on elements $K\in\D_h$ and therefore is called a \emph{$p$-(adaptive) refinement}. If \emph{all} elements of the mesh $\D_h$ are refined or the polynomial degree of \emph{all} elements in $\D_h$ is increased, one calls the refinement a \emph{uniform} $h$- or a \emph{uniform $p$-refinement}, respectively.

Convergence of the finite element method may be achieved by progressively refining all elements $K\in\D_h$ or increasing the polynomial degrees $p_K$ of all elements $K\in\D_h$, leading to the so-called \emph{$h$-version} or \emph{$p$-version} of the finite element method, respectively. Typically, the convergence of the method is specified as the discretization error per degrees of freedom. Thereby, so-called \emph{a priori error estimates}, which are based on a priori knowledge, i.e.~the characteristics of the weak solution (e.g.~its regularity) and discretization parameters such as $h$ and $p$ but do not include quantities computed by the method, give the asymptotic rate of convergence that can be expected for a method at most. Usually, such estimates take the form
\begin{align*}
 \norm{u_\XX - u_\WW}_\XX
  \leq c \, h^r
\end{align*}
for some constant $c$, where $r$ is the rate of convergence. While the $h$-version leads to an algebraic convergence of the order $h^p$ at best, where $p:=\min_{K\in\D_h} p_K$, the $p$-version of the method results in an exponential decrease of the error if the weak solution is smooth, see e.g.~\cite{ref:Schwab_1998, ref:Szabo_1991}. However, applying uniform refinements leads to a aggravation of the convergence of the discretization error with respect to the degrees of freedom if the weak solution has a low regularity. Such a low regularity may result from the geometry of the domain $\Omega$, for instance, in the presence of a re-entrant corner, changes of material properties inside the domain, such as a transition from pure plastic to elastoplastic behavior, or changes of the boundary conditions. While in regions where the solution is smooth large elements with a high local polynomial degree turn out to be appropriate, in regions where the solution has singularities the mesh should be refined towards these singularities with elements of a low local polynomial degree to obtain highest possible convergence rates. Hence, to recover the optimal algebraic convergence rates it may be necessary to apply $h$-refinements only locally, resulting in a so-called \emph{$h$-adaptive} method. Simultaneously applying $h$- and $p$-refinements leads to the celebrated \emph{$hp$-adaptive finite element method}, which allows to achieve exponential convergence rates even for weak solutions of a low regularity, see e.g.~\cite{ref:Gui_1985, ref:Schwab_1998, ref:Szabo_1991}. Thereby, the a priori results of Babu\v{s}ka and co-authors, cf.~\cite{ref:Babuska_1988, ref:Babuska_1987, ref:Babuska_1989, ref:Babuska_1994, ref:Guo_1986_1, ref:Guo_1986_2}, form the basis of $hp$-adaptive strategies.

% -----------------------------------------------------------------------------------------

\subsubsection*{Adaptivity and A Posteriori Error Estimation}

In order to derive adaptive methods one either needs to have a priori knowledge about the regularity of a weak solution, in particular, the location of possible singularities, or to apply \emph{a posteriori error estimates} to identify regions where the discretization error is large. While a priori error estimates exclusively resort to a priori knowledge and give the convergence rate of a specific method but cannot be used to quantify the local or global error of a discrete solution, \emph{a posteriori error estimates} bound the current error in terms of known computable quantities, which may be determined while or after computing a discrete solution, such as the discrete solution $u_\WW$ itself or given data like a volume force or surface traction in the context of elastoplasticity. As in general the weak solution of a problem is unknown a posteriori error estimators have to be designed to steer the mesh refinements within an adaptive procedure. Thereby, the error estimator has to be expressible as local contributions given on the physical elements $K\in\D_h$ in order to determine those elements on which the discretization error is large and which therefore have to be refined. In the last decades, significant contributions have been made in the field of a posteriori error estimation and automatic mesh adaption, see e.g.~\cite{ref:Ainsworth_1996, ref:becker_2001, ref:eriksson_1995, ref:houston_2002, ref:verfuerth_1996} as well as in the study of optimal convergence rates for adaptive methods, see e.g.~\cite{ref:Carstensen_2014, ref:Carstensen_2006, ref:Cascon_2008, ref:chen_2009, ref:Stevenson_2007}. While for many years such methods usually exclusively applied (mostly isotropic) adaptive $h$-refinements in recent years, however, also $hp$-adaptive methods are derived, see e.g.~the overview article \cite{ref:mitchell_2014} and the references therein or the text books \cite{ref:demkovicz_2007, ref:Schwab_1998, ref:Solin_2004}.
\vspace{0.1cm}

As the weak solution of problems in elastoplasticity typically does not enjoy high regularity properties in this area there is a strong need for adaptive methods to obtain adequate convergence properties. Therefore, the aim of the first part of this thesis is to introduce and analyze $hp$-finite element discretizations of a model problem in elastoplasticity and demonstrate the potential of the resulting $h$- and $hp$-adaptive methods. Thereby, a reliable a posteriori error estimator is derived to steer the adaptive refinements and an efficient semi-smooth Newton solver with superlinear convergence properties is proposed to compute the finite element approximation. Lower order $h$- and $p$-versions of the finite element method in the context of elastoplastic problems were already considered e.g.~in \cite{ref:Alberty_2000, ref:Carstensen_1999, ref:Han_2013, ref:Kienesberger_2006}. For adaptivity methods, see e.g.~\cite{ref:Carstensen_2016}. In the paper [P4], we introduce an alternative strategy to steer the $hp$-adaptivity of a finite element method for discretizing self-adjoint elliptic boundary value problems, which neither relies on a posteriori error estimates nor on smoothness indicators.

% -----------------------------------------------------------------------------------------
%   Some Aspects of Implementation
% -----------------------------------------------------------------------------------------

\section{Some Aspects of Implementation}\label{sec:implementation}

Frequently, the quantities appearing in the finite element system, which results from discretizing the weak formulation of a boundary value problem and is solved by some numerical method, are computed by means of local quantities given on the physical elements $K\in\D_h$ of the mesh $\D_h$. This element-wise computation is generally known as \emph{assembling}. For the ease of presentation, let the weak formulation be given by: Find a $u_\XX\in\XX$ such that the variational equation
\begin{align}\label{eq:abstract_VEQ}
 a(u_\XX, v) = b(v) \qquad
 \forall \, v\in\XX
\end{align}
holds true with a bounded, $\XX$-elliptic bilinear form $a: \XX\times \XX\longrightarrow \RR$ and a bounded linear form $b:\XX\longrightarrow\RR$ on the Hilbert function space $\XX = \XX(\Omega)$ over a domain $\Omega$. Then, the unique existence of a weak solution $u_\XX$ is guaranteed by the \emph{Lax-Milgram-Lemma}, see e.g.~\cite{ref:Braess_2013, ref:Brenner_2008, ref:Han_2013}, and a finite element discretization of \eqref{eq:abstract_VEQ} is given by: Find a $u_\WW\in\WW$ such that the discrete variational equation
\begin{align}\label{eq:abstractDiscrete_VEQ}
 a(u_\WW, w) = b(w) \qquad
 \forall \, w\in\WW
\end{align}
is valid, where the finite element space $\WW$ takes the form \eqref{eq:abstract_FE_space}. As $\WW$ is spanned by finitely many degrees of freedom $\phi_1,\ldots,\phi_N$ the discrete solution $u_\WW$ can be represented in terms of a linear combination $u_\WW = \sum_{i\in\ul{N}} u_i \, \phi_i$ with some coefficients $u_1,\ldots,u_N\in\RR$, where the notation $\ul{n}:=\lbrace 1,\ldots,n\rbrace$ for any positive integer $n\in\NN$ is used. By \eqref{eq:abstractDiscrete_VEQ} and exploiting the fact that also each $w\in\WW$ is a linear combination of the degrees of freedom, the coefficient vector $\myvec{u} = (u_1,\ldots,u_N)^{\top}\in\RR^N$ is uniquely determined by the linear \emph{finite element system}
\begin{align}\label{eq:abstact_linSys}
 \mymatrix{A}_\WW \, \myvec{u} = \myvec{b}_\WW,
\end{align}
where the so-called \emph{(global) stiffness matrix} $\mymatrix{A}_\WW = (a_{ij})\in\RR^{N\times N}$ and the \emph{(global) load vector} $\myvec{b}_\WW = (b_i)\in\RR^N$ are given component-wise by
\begin{align*}
 a_{ij} 
  := a(\phi_j, \phi_i), \qquad
 b_i 
  := b(\phi_i) \qquad
 \forall \, i,j\in\ul{N}.
\end{align*}
As so far only the finite dimension of $\WW$ was used any Riesz-Galerkin discretization of the infinite-dimensional variational equation \eqref{eq:abstract_VEQ} can be formulated as a linear system of equations such as \eqref{eq:abstact_linSys}.

% -----------------------------------------------------------------------------------------

\subsubsection*{Assembling of the Global Quantities}

As for $K\in\D_h$ the local space $\Psi_K$ is spanned by finitely many shape functions $\psi_1^K,\ldots,\psi_{N_K}^K$ the restriction $v_{\, | \, K}$ of any $v\in\XX$ to the physical element $K$ can be represented in terms of a linear combination of these shape functions. In particular,
\begin{align}\label{eq:repres_restriction}
 \phi_{i \, | \, K} = \sum_{j\in\ul{N_K}} c_{ij}^K \, \psi_j^K \qquad
 \forall \, i\in\ul{N}
\end{align}
for some uniquely determined coefficients $c_{ij}^K\in\RR$, which are the entries of the so-called \emph{connectivity matrix} $\mymatrix{C}_K = (c_{ij}^K)\in\RR^{N\times N_K}$ for $K\in\D_h$. Hence, the degrees of freedom of a finite element space might be constructed element-wise with the help of suitable shape functions on the physical elements. For the applicability of efficient solution schemes to solve the finite element system \eqref{eq:abstact_linSys} one obviously will try to achieve a sparse matrix as (global) stiffness matrix. Therefore, the choice of the degrees of freedom, spanning the finite element space, or rather the choice of shape functions on the physical elements to construct the global degrees of freedom in general depends on the bilinear form $a(\cdot,\cdot)$. In either case, it is advantageous if the degrees of freedom have a possible small support. If the bilinear form $a(\cdot, \cdot)$ and the linear form $b(\cdot)$ are decomposable in the sense that
\begin{align}\label{eq:decomp_ab}
 a(v,w) 
  = \sum_{K\in\D_h} a_K \big( v_{\, | \, K}, w_{\, | \, K} \big), \qquad
 b(v) 
  = \sum_{K\in\D_h} b_K \big( v_{\, | \, K} \big) \qquad
 \forall \, v,w \in \XX
\end{align}
for some local bilinear forms $a_K : \XX_K\times \XX_K\longrightarrow\RR$ and local linear forms $b_K:\XX_K\longrightarrow\RR$ on the restriction spaces $\XX_K := \lbrace v_{\, | \, K} \; ; \; v\in\XX \rbrace$ for $K\in\D_h$, the global quantities $\mymatrix{A}_\WW$ and $\myvec{b}_\WW$ can be assembled element-wise by the formulas
\begin{align}\label{eq:abstract_assembling}
 \mymatrix{A}_\WW
  = \sum_{K\in\D_h} \mymatrix{C}_K \, \mymatrix{A}_K \, \mymatrix{C}_K^{\top}, \qquad
 \myvec{b}_\WW
  = \sum_{K\in\D_h} \mymatrix{C}_K \, \myvec{b}_K,
\end{align}
where, for $K\in\D_h$, the \emph{local stiffness matrix} $\mymatrix{A}_K = (a_{ij}^K)\in\RR^{N_K\times N_K}$ and the \emph{local load vector} $\myvec{b}_K = (b_i^K)\in\RR^{N_K}$ are defined component-wise as
\begin{align*}
 a_{ij}^K 
  := a_K \big( \psi_j^K, \psi_i^K \big), \qquad
 b_i^K 
  := b_K \big(\psi_i^K \big) \qquad
 \forall \, i,j\in\ul{N_K}.
\end{align*}

If the local bilinear form $a_K(\cdot,\cdot)$ represents an inner product on $\XX_K$ then for the implementation it is very beneficial to use shape functions which are orthogonal with respect to $a_K(\cdot,\cdot)$. It should also be mentioned that the computation of the local stiffness matrix and the local load vector may be done in parallel, which allows an efficient assembling of the global quantities $\mymatrix{A}_\WW$ and $\myvec{b}_\WW$. In this context, the use of a so-called \emph{reference element} frequently turns out to be particularly beneficial. If for any $K\in\D_h$ the shape functions are the image of (the same) linearly independent functions on some fixed, so-called \emph{reference element} $\wh{K}\subseteq\RR^d$ under a suitable, bijective transformation $\myvec{F}_K : \wh{K}\longrightarrow K$ one can avoid to explicitly construct shape functions on each physical element but only has to determine functions $\wh{\psi}_j$ on the reference element $\wh{K}$ for which $\psi_j^K = \psi_j \circ \myvec{F}_K$ for $j\in\ul{N_K}$ and $K\in\D_h$. In this case, the global degrees of freedom can be constructed via images of the same linearly independent functions on the reference element $\wh{K}$.

% -----------------------------------------------------------------------------------------

\subsubsection*{Specific Physical Elements and Degrees Of Freedom}

\begin{figure}[ht]
 \centering
 \begin{subfigure}[b]{0.4\textwidth}
  \centering
  \includegraphics[trim = 8cm 1cm 8cm 1cm, clip, width = 5cm]{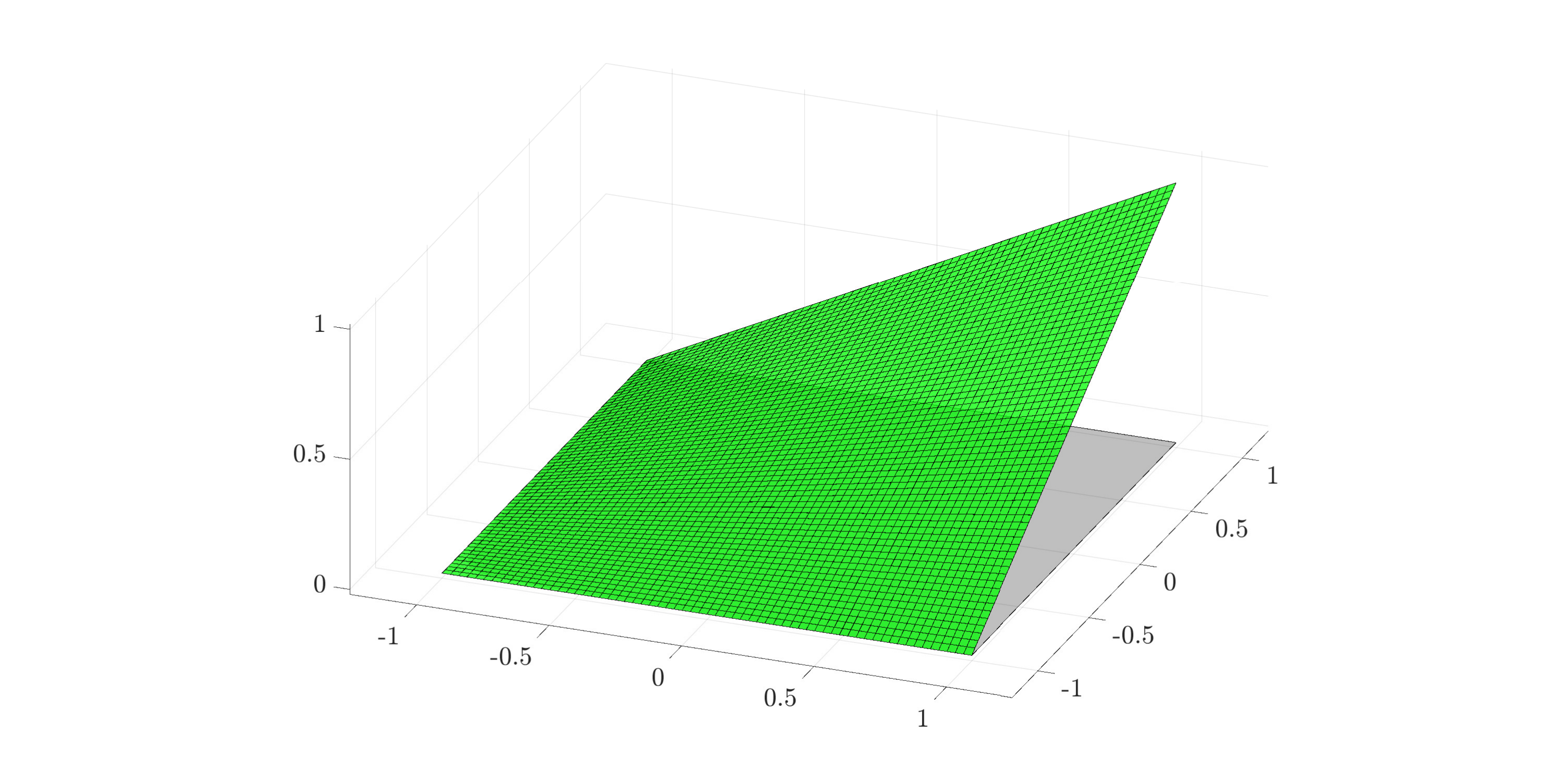}
  \caption{Shape function related to vertex.}
  %\label{fig:vertex_shapeF}
 \end{subfigure}
 %
 %\quad
 %
 \begin{subfigure}[b]{0.4\textwidth}
  \centering
  \includegraphics[trim = 3.5cm 3cm 4cm 3.5cm, clip, width = 5.25cm]{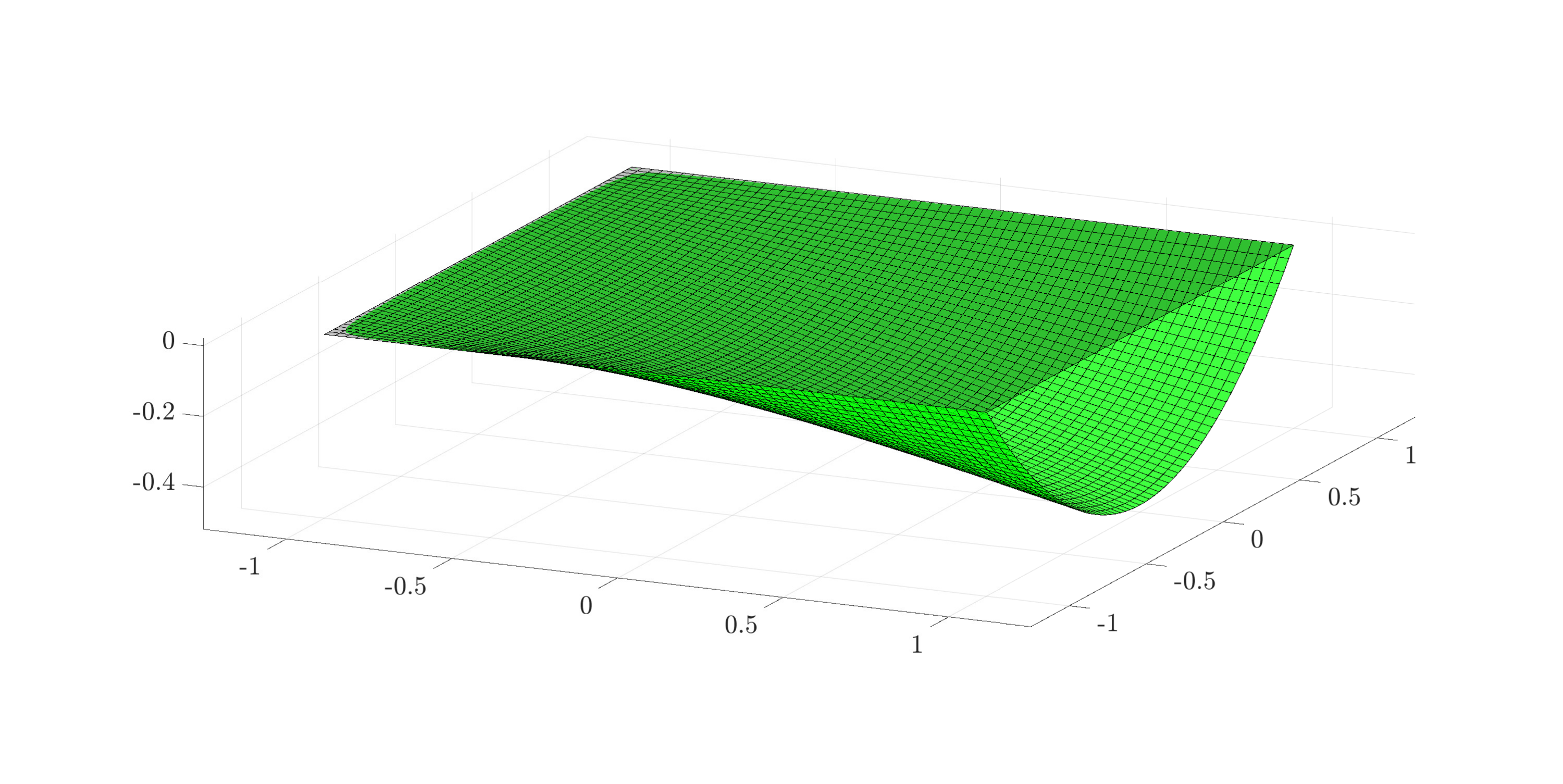}
  \caption{Shape function related to edge.}
  %\label{fig:edge_schapeF}
 \end{subfigure}
 %
 %\quad
 %
 \begin{subfigure}[b]{0.4\textwidth}
  \centering
  \includegraphics[trim = 4cm 4cm 4cm 4cm, clip, width = 5.5cm]{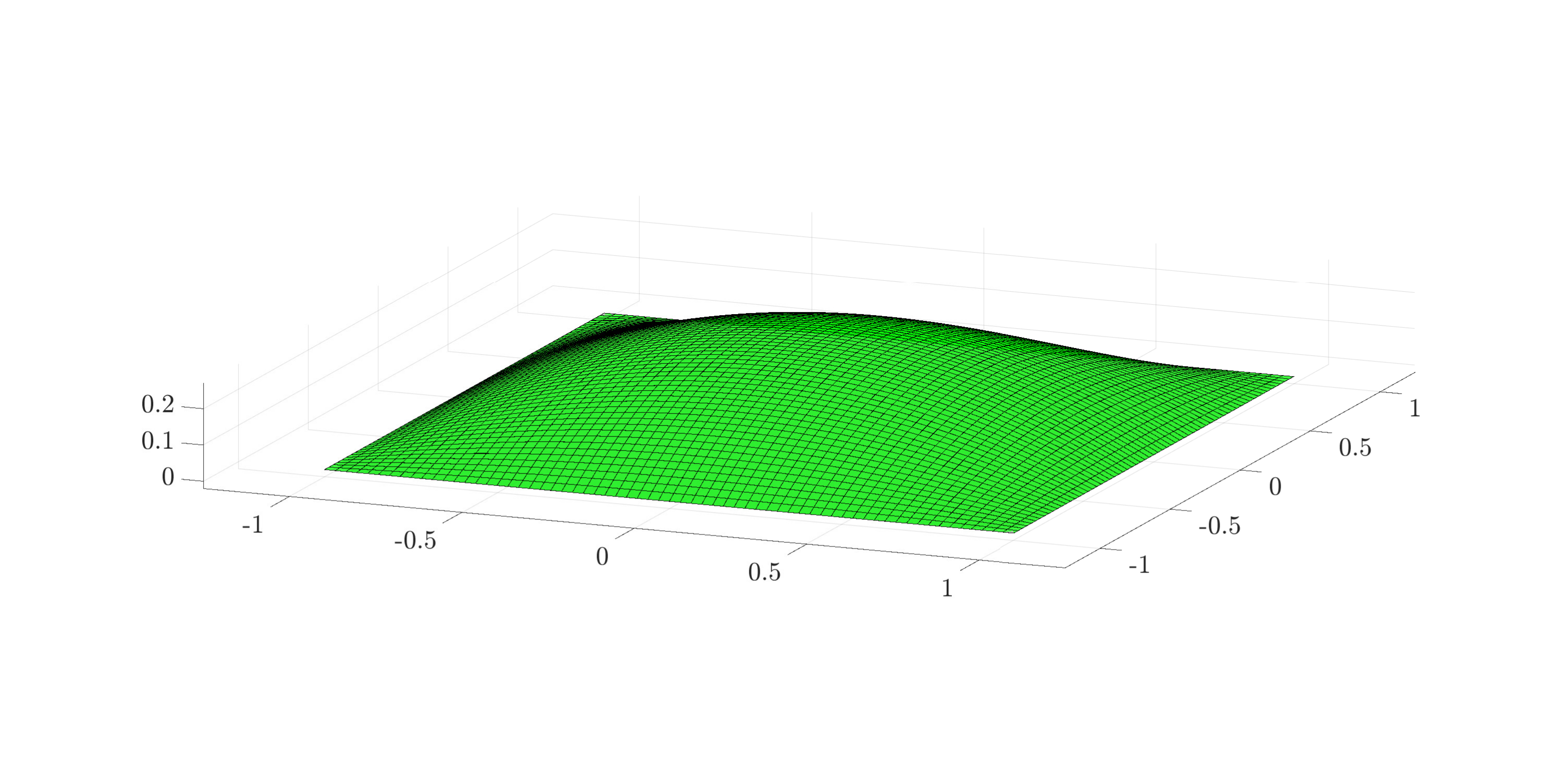}
  \caption{Shape function related to element.}
  %\label{fig:edge_schapeF}
 \end{subfigure}
 \caption{\it Two-dimensional shape functions constructed via images of integrated Legendre polynomials.}
 \label{fig:shape_functions}
\end{figure}

In [P4], the above ideas are used to construct so-called \emph{enrichment functions} on the physical elements and to assemble the global quantities which are needed to compute the predicted local error reduction. Furthermore, the implementation of the assembling of the global quantities included in the finite element system is in all papers based on the above concepts. In [P4], as well as in the papers [P1, P2, P3], \emph{transformed hexahedrons} $Q$ (see [P4, Sec.~3.2] for the definition in arbitrary dimension), are used as physical elements, which are the image of the reference element $\wh{Q} := [-1,1]^d$ under suitable, bijective mappings $\myvec{F}_Q : \wh{Q}\longrightarrow Q$. Moreover, in all papers $hp$-finite element spaces of the form
\begin{align*}
 \WW_{hp} = \Big\lbrace v\in\XX \; ;\; v_{\, | \, Q}\circ\myvec{F}_Q\in\big( \PP_{p_Q}(\wh{Q}) \big) \text{ for all } Q\in\T_h \Big\rbrace
\end{align*}
are considered where $\T_h$ denotes the mesh of $\Omega$ consisting of transformed hexahedrons $Q\in\T_h$, $h = (h_Q)_{Q\in\T_h}$, $p = (p_Q)_{Q\in\T_h}$ and $\PP_{p_Q}(\wh{Q})$ is the space of polynomials up to degree $p_Q$ on the reference element $\wh{Q}$. Using transformed hexahedrons being the image of the reference element $\wh{Q}$ offers the advantage that by the tensor structure of $\wh{Q}$ polynomials on $\wh{Q}$ can easily be constructed via tensor products of polynomials on the segment $[-1,1]$, which span the space of polynomials on $[-1,1]$, cf.~[P4, Sec.~3.1]. This, in particular, is beneficial for applying Gauss quadrature rules in order to evaluate integrals over the physical elements. In all papers, the degrees of freedom that span the finite element spaces consisting of globally continuous piecewise images of polynomials are constructed via shape functions being the images of tensor products of \emph{integrated Legendre polynomials}, see Figure~\ref{fig:shape_functions}. First of all, this is done because while the assembling derivatives of these functions appear in the considered cases, which are given by the \emph{Legendre polynomials} that are orthogonal with respect to the $L^2$-inner product (see e.g.~\cite{ref:Haemmerlin_1991} for details on Legendre polynomials). Secondly, the resulting degrees of freedom can be associated with the so-called \emph{nodes} of the mesh, i.e.~either with a vertex, an edge or an higher-dimensional face (including the $d$-dimensional physical elements) of the mesh $\Omega$. In [P4, Sec.~3.4.3], enrichment functions on an isotropic refined element with this property are constructed.
%\vspace{0.1cm}

% -----------------------------------------------------------------------------------------

% -----------------------------------------------------------------------------------------

If adaptive isotropic $h$-refinements are applied to a mesh of transformed hexahedrons, however, one has to face the inconvenience of \emph{hanging nodes} to preserve the global continuity of the degrees of freedom. Thereby, vertices, edges or higher-dimensional faces (up to dimension $d-1$) that can \emph{not} be associated with degrees of freedom -- so-called \emph{hanging nodes} -- naturally arise in the mesh whenever physical elements stay unrefined while their neighbouring elements are refined. One possible way to handle hanging nodes is the so-called \emph{constrained approximation}, see e.g.~\cite{ref:demkovicz_2007, ref:demkovicz_2008, ref:Karniadakis_2013, ref:Schroeder_2011_constrAppr, ref:Schwab_1998, ref:Solin_2004, ref:Zienkiewicz_2005}, which is used for the implementation in [P1, P2, P3, P4]. To preserve the global continuity of the degrees of freedom the shape functions associated with hanging nodes are constrained in  this approach. In \cite{ref:DiStolfo_2016, ref:DiStolfo_2020} an easy treatment of this challenging to implement technique is presented.

% -----------------------------------------------------------------------------------------
%   Basic Principles of Elasticity
% -----------------------------------------------------------------------------------------

\section{Basic Principles of Elasticity}

In problems of elasticity and elastoplasticity, the behavior of a material body is considered, which is subjected to different kind of forces acting on the entire body and its boundary, respectively. By regarding the body (in a macroscopic level) to be composed of a continuously distributed material one may identify its undeformed and unstressed state -- the so-called \emph{reference configuration} -- with a bounded domain $\Omega\subseteq \RR^d$ (typically of dimension $d\in\lbrace 2,3\rbrace$) and refers to the points $\myvec{x}\in\Omega$ as to \emph{material points}. In order to model the behavior of such a body by a system of partial differential equations its behavior as well as its (material) properties have to be expressed in terms of functions of position $\myvec{x}\in\Omega$ and time $t$. In a first step, the motion and deformation of a body have to be described within the framework of \emph{continuum mechanics}, see e.g.~\cite{ref:Chadwick_1999, ref:Ciarlet_1988, ref:Gurtin_2010}, for further details. For a general theory of elasticity and elastoplasticity, see the monographs \cite{ref:Antman_2005, ref:Chen_1988, ref:Ciarlet_1988, ref:Han_2013}.

% -----------------------------------------------------------------------------------------
\subsubsection*{Kinematics and Stress}

Due to the applied forces the body is moving and deforming with time so that it occupies a domain $\Omega_t\in\RR^d$ at time $t$, which is called \emph{current configuration}. Without loss of generality, $\Omega$ is the configuration at the time $t=0$. Hence, the new position $\myvec{y}\in\Omega_t$ of the material point $\myvec{x}\in\Omega$ may be expressed by the vector-valued \emph{motion} $\myvec{y} : \Omega\times [0,T]\longrightarrow\RR^d$ with $\myvec{y}(\myvec{x},t) = \myvec{y}$. In order to describe the behavior of the body, however, it turns out to be more convenient to consider the vector-valued \emph{displacement} $\myvec{u} : \Omega\times[0,T]\longrightarrow\RR^d$, which is given by
\begin{align}\label{eq:general_strainTensor}
 \myvec{u}(\myvec{x},t) 
  := \myvec{y}(\myvec{x},t) - \myvec{x}.
\end{align}

\begin{figure}[ht]
\centering
\includegraphics[trim = 6cm 1.75cm 5cm 1.5cm, clip, width = 4cm]{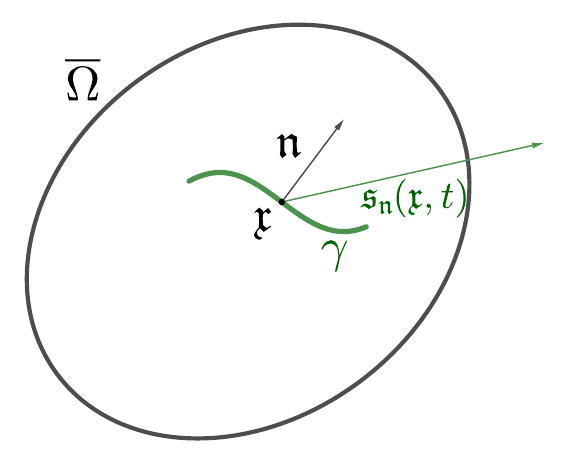}
\caption{\it The stress vector $\myvec{s}_{\myvec{n}}(\myvec{x},t)$.}
\label{fig:stress_vec}
\end{figure}

To distinguish between \emph{rigid body motions}, in which the body is \emph{only} translated and rotated, and real deformations, which change the body's shape, the matrix-valued \emph{strain tensor} $\mymatrix{\eta}$ associated with the displacement field $\myvec{u}=(u_1,\ldots,u_d)^{\top}$ is introduced. Thereby, $\mymatrix{\eta}$ is defined as
\begin{align}\label{eq:def_strain}
 \mymatrix{\eta}(\myvec{u})
  := \frac{1}{2} \Big( \nabla\myvec{u} + (\nabla\myvec{u})^{\top} + (\nabla\myvec{u})^{\top} \nabla\myvec{u} \Big)
\end{align}
and measures the deformation of the body, where the gradient $\nabla\myvec{u}$ is with respect to the variable $\myvec{x} = (x_1,\ldots,x_d)^{\top}$, i.e.~consists of the components $\partial u_i/\partial x_j$ for $i,j\in\ul{d}$. Therewith, the body undergoes a rigid body motion if and only if $\mymatrix{\eta} = \mymatrix{0}$, cf.~\cite{ref:Han_2013}. In many problems of practical interest the deformations can be regarded as \emph{small}, which is also assummed in the model problem considered in [P1, P2, P3]. In the case of so-called \emph{infinitesimal deformations} $\nabla\myvec{u}$ is regarded to be sufficient small to neglect the nonlinear term in \eqref{eq:general_strainTensor}. Thereby, the strain tensor $\mymatrix{\eta}$ is replaced by the \emph{infinitesimal strain tensor} $\mymatrix{\varepsilon}$, given by
\begin{align}\label{eq:def_infinitesimal_strain}
 \mymatrix{\varepsilon}(\myvec{u}) 
  :=  \frac{1}{2} \, \Big( \nabla\myvec{u} + (\nabla\myvec{u})^{\top} \Big).
\end{align}
Then, an infinitesimal rigid body motion is characterized by $\mymatrix{\varepsilon}(\myvec{u}) = \mymatrix{0}$. The state of internal forces acting in the body is described by means of a so-called \emph{stress tensor} for its definition we first classify the applied forces:
\begin{itemize}
 \item A \emph{volume force} $\myvec{f} : \Omega\times[0,T]\longrightarrow\RR^d$ represents the force per unit reference volume. Gravity, for instance, is given by the volume force $\myvec{f}(\myvec{x},t) = \rho \, g \, \myvec{e}$ for $\myvec{x}\in\Omega$ and $t\in[0,T]$, where $\rho$ is the mass density of the body, $g$ is the gravitational acceleration and $\myvec{e}$ denotes the unit vector pointing in the downward vertical direction.
 
 \item A \emph{surface traction} acts on the body's boundary. If $\gamma$ is a regular surface in $\overline{\Omega}$ passing through $\myvec{x}$ with unit normal $\myvec{n}$ at $\myvec{x}$, cf.~Figure~\ref{fig:stress_vec}, then the \emph{stress vector} $\myvec{s}_{\myvec{n}}(\myvec{x},t)$ is defined to be the current force per unit area exerted by the portion of $\Omega$ on the side of $\gamma$ in which $\myvec{n}$ points, on the portion of $\Omega$ which lies on the other side. For an arbitrary subset $\Omega'\subseteq\Omega$ with boundary $\Gamma'$ the \emph{surface traction} at time $t$ is the stress vector $\myvec{s}_{\myvec{n}}(\myvec{x},t)$ ($\myvec{x}\in\Gamma'$) acting on $\Gamma'$.
\end{itemize}

% -----------------------------------------------------------------------------------------

% -----------------------------------------------------------------------------------------

It can be shown that there exists a tensor field, the so-called \emph{first Piola-Kirchoff stress}, $\mymatrix{\tau}:\Omega\times[0,T]\longrightarrow\RR^{d\times d}$ with the property $\mymatrix{\tau} \, \myvec{n} = \myvec{s}_{\myvec{n}}$ for any unit vector $\myvec{n}$, cf.~\cite{ref:Lemaitre_1990}. Therefore, the \emph{divergence theorem} of Gauss yields
\begin{align}\label{eq:prop_stress}
 \int_{\Gamma'} \myvec{s}_{\myvec{n}} \dd \myvec{s}
  = \int_{\Gamma'} \mymatrix{\tau} \, \myvec{n} \dd \myvec{s}
  = \int_{\Omega'} \Div \mymatrix{\tau} \dd \myvec{x},
\end{align}
where the divergence $\Div \mymatrix{\tau}$ is with respect to the variable $\myvec{x}$. Introducing the \emph{acceleration field} $\ddot{\myvec{u}}(\myvec{x},t) := \frac{\partial^2}{\partial t^2} \, \myvec{u}(\myvec{x},t)$ and exploiting \eqref{eq:prop_stress} one can write the \emph{balance of linear momentum} in terms of
\begin{align*}
 \myvec{0} 
  = \int_{\Omega'} \rho \, \ddot{\myvec{u}} \dd \myvec{x} - \bigg( \int_{\Omega'} \myvec{f} \dd \myvec{x} + \int_{\Gamma'} \myvec{s}_{\myvec{n}} \dd \myvec{s} \bigg)
  = \int_{\Omega'} \rho \, \ddot{\myvec{u}} - \myvec{f} - \Div \mymatrix{\tau} \dd\myvec{x}
\end{align*}
for any subset $\Omega'\subseteq\Omega$ with boundary $\Gamma'$, which immediately yields the local form of the \emph{equation of motion}
\begin{align}\label{eq:equation_of_motion}
 \Div \mymatrix{\tau} + \myvec{f} 
  = \rho \, \ddot{\myvec{u}}.
\end{align}
If the given data is independent of time and, thus, the response of the body can be regarded as independent of time as well, the equation \eqref{eq:equation_of_motion} becomes the \emph{equation of equilibrium}
\begin{align*}
 \Div \mymatrix{\tau} + \myvec{f}
  = \myvec{0}.
\end{align*}

The above considerations related to the reference configuration $\Omega$ can also be expressed with respect to the current configuration $\Omega_t$. Then, in analogy to $\mymatrix{\tau}$ there exists the so-called \emph{Cauchy stress} $\mymatrix{\sigma} : \Omega_t\times[0,T]\longrightarrow\RR^{d\times d}$, which is related to the first Piola-Kirchoff stress through the identity
\begin{align*}
 \mymatrix{\sigma} = \frac{1}{\det \nabla\myvec{y}} \, \mymatrix{\tau} \, (\mymatrix{I} + \nabla\myvec{u})^{\top}.
\end{align*}
As a consequence of the \emph{balance of angular momentum}, $\mymatrix{\tau} \, (\mymatrix{I} + \nabla\myvec{u})^{\top}$ is symmetric, see~\cite{ref:Han_2013}, which implies the symmetry of the stress tensor $\mymatrix{\sigma}$. The equation of motion \eqref{eq:equation_of_motion}, with respect to the current configuration, takes the form $\Div \mymatrix{\sigma} + \myvec{f} = \rho_t \, \myvec{a}$, where $\rho_t$ is the mass density per unit current volume, $\myvec{a}$ is the acceleration and $\Div \mymatrix{\sigma}$ is with respect to the variable $\myvec{y}$. For infinitesimal deformations a distinction between the reference and current configuration may be neglected and, thus, also the distinction between the two stress tensors $\mymatrix{\tau}$ and $\mymatrix{\sigma}$. Furthermore, derivatives with respect to $\myvec{y}$ might be replaced by derivatives with respect to $\myvec{x}$, i.e.~for infinitesimal deformations the equations
\begin{subequations}
\begin{align}
 \Div \mymatrix{\sigma} + \myvec{f} &= \rho \, \ddot{\myvec{u}}, \label{eq:infinitesimal_equation_of_motion}\\
 \mymatrix{\sigma} &= \mymatrix{\sigma}^{\top}
\end{align}
\end{subequations}
hold true, where the right-hand side of the first equation equals $\myvec{0}$ if, in addition, the data is assumed to be independent of time, which is the case for the model problem in [P1, P2, P3].
\vspace{0.1cm}

Note that so far only the \emph{kinematics} of a body have been discussed which \emph{do not} contain any particular material behavior. The description is therefore incomplete in a physical point of view. Likewise, the heretofore derived equations show the mathematically incompleteness of the problem. For instance, in three dimensions, i.e.~in the case $d=3$, the equation of motion \eqref{eq:equation_of_motion} as well as its version \eqref{eq:infinitesimal_equation_of_motion} for infinitesimal deformations, when written out component-wise, lead to three equations. The strain-displacement relation \eqref{eq:def_infinitesimal_strain} results in another six equations (taking the symmetry of $\mymatrix{\varepsilon}$ into account), which leads to a total of nine equations. However, the three components of the displacement, the six components of the strain and the six components of the stress (taking the symmetry of $\mymatrix{\varepsilon}$ and $\mymatrix{\sigma}$ into account) give a total number of fifteen unknowns. Hence, six additional equations are missing that the problem becomes at least in principle solvable. These equations are provided by \emph{constitutive equations} of the specific material that is modeled. The constitutive equations of elastoplastic materials can thereby be seen as an extension of the linear elastic case.

% -----------------------------------------------------------------------------------------

\subsubsection*{Linear Elastic Materials}

A material body is said to be \emph{linearly elastic} if the stress $\mymatrix{\sigma}$ depends linearly on the infinitesimal strain $\mymatrix{\varepsilon}$, i.e.~if
\begin{align}\label{eq:constitutie_equation_linElast}
 \mymatrix{\sigma}
  = \CC \, \mymatrix{\varepsilon},
\end{align}
where the so-called \emph{elasticity tensor} $\CC$ is a linear map from the space of symmetric second-order tensors into itself. In general, the elasticity tensor depends on the position $\myvec{x}\in\Omega$ but \emph{does not} depend on the time $t$. If both, the mass density $\rho$ as well as the elasticity tensor $\CC$ are independent of position, the body is called \emph{homogeneous}. In view of \eqref{eq:constitutie_equation_linElast}, the elasticity tensor may be assumed to have the symmetry properties
\begin{align*}
 \CC_{ijkl} 
  = \CC_{jikl} 
  = \CC_{ijlk} \qquad
 \forall \, i,j,k,l\in\underline{d}
\end{align*}
due to the symmetry of $\mymatrix{\varepsilon}$ and $\mymatrix{\sigma}$. In the mathematically modeling, $\CC$ is typically assumed to have the additional symmetry property $\CC_{ijkl} = \CC_{klij}$ for $i,j,k,l\in\ul{d}$ and to be \emph{uniformly elliptic}, i.e.~it exists a positive constant $c_e > 0$ such that
\begin{align*}
 (\CC \, \mymatrix{\tau}) : \mymatrix{\tau} \geq c_e \, \abs{\mymatrix{\tau}}_F^2
\end{align*} 
for all symmetric second-order tensors $\mymatrix{\tau}$, where $:$ and $\abs{\cdot}_F$ denote the Frobenius inner product and its induced norm, respectively. In the special case that a material has no preferred direction and, thus, responds to a force independently of its orientation it is called \emph{isotropic}. In this case, the twenty-one independent components of $\CC$ (taking the symmetry properties of $\CC$ into account) reduce to \emph{two}. The choice of these two material coefficients is, of course, not unique but a common one are the so-called \emph{Lam\'{e} moduli} $\lambda$ and $\mu$, which turn \eqref{eq:constitutie_equation_linElast} into
\begin{align}\label{eq:constitutive_equation_isotropic}
 \mymatrix{\sigma}
  = \lambda \, (\tr \mymatrix{\varepsilon}) \, \mymatrix{I} + 2 \, \mu \, \mymatrix{\varepsilon}.
\end{align}
The numerical examples in [P1, P2, P3] consider isotropic materials and, hence, in the implementation the appearing stress of the elastoplastic model is given by \eqref{eq:constitutive_equation_isotropic}.
\vspace{0.1cm}

A complete mathematical formulation for describing the deformation of and stresses in a linearly elastic body can now be stated where, for simplicity, the data is assumed to be independent of time. Let the body initially occupy the bounded domain $\Omega\subset\RR^d$ with boundary $\Gamma := \partial \Omega$ and outer unit normal $\myvec{n}$, where $\Gamma$ is decomposed into the non-overlapping parts $\Gamma_D$ and $\Gamma_N$ such that $\Gamma = \overline{\Gamma_D} \cup \overline{\Gamma_N}$. Then, for given volume force $\myvec{f} : \Omega\longrightarrow\RR^d$, displacement $\overline{u} : \Gamma_D\longrightarrow\RR^d$ and surface traction $\myvec{g} : \Gamma_N\longrightarrow\RR^d$ the \emph{boundary value problem of linearized elasticity} is: Find a displacement field $\myvec{u} : \Omega\longrightarrow\RR^d$ that satisfies the \emph{equation of equilibrium}
\begin{align*}
 - \Div \big(\mymatrix{\sigma}(\myvec{u})\big) &= \myvec{f} \qquad \text{in } \Omega
\end{align*}
as well as the boundary conditions
\begin{align*}
 \myvec{u} = \overline{\myvec{u}} \quad \text{on } \Gamma_D, \qquad
 \mymatrix{\sigma}(\myvec{u}) \, \myvec{n} = \myvec{g} \quad \text{on } \Gamma_N,
\end{align*}
where $\mymatrix{\sigma}$ and $\mymatrix{\varepsilon}$, for $\myvec{x}\in\Omega$, are given by the \emph{elastic constitutive relation} $\mymatrix{\sigma}(\myvec{u}) = \CC \, \mymatrix{\varepsilon}(\myvec{u})$ and the \emph{strain-displacement relation} $\mymatrix{\varepsilon}(\myvec{u}) = \frac{1}{2} \, \big( \nabla\myvec{u} + (\nabla\myvec{u})^{\top} \big)$, respectively.

% -----------------------------------------------------------------------------------------
%   Elastoplastic Materials
% -----------------------------------------------------------------------------------------

\section{Elastoplastic Materials}

\begin{figure}[ht]
\centering
 \begin{subfigure}[b]{0.4\textwidth}
  \centering
  \includegraphics[trim = 5cm 2cm 9.5cm 2.5cm, clip, width = 4.5cm]{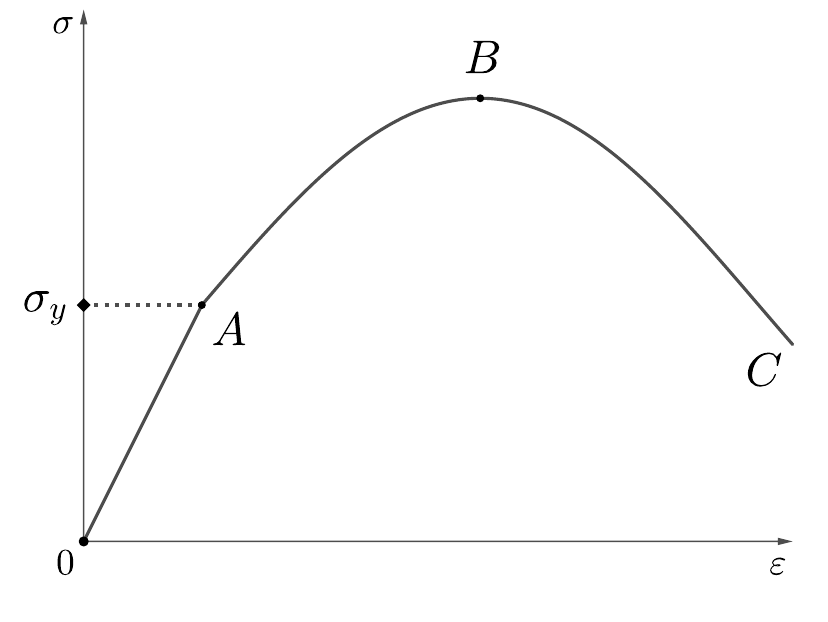}
  \caption{Hardening and softening.}
  \label{fig:stress_strain_relation}
 \end{subfigure}
 %
 %\quad
 %
 \begin{subfigure}[b]{0.4\textwidth}
  \centering
  \includegraphics[trim = 5cm 0.5cm 4cm 1cm, clip, width = 4.5cm]{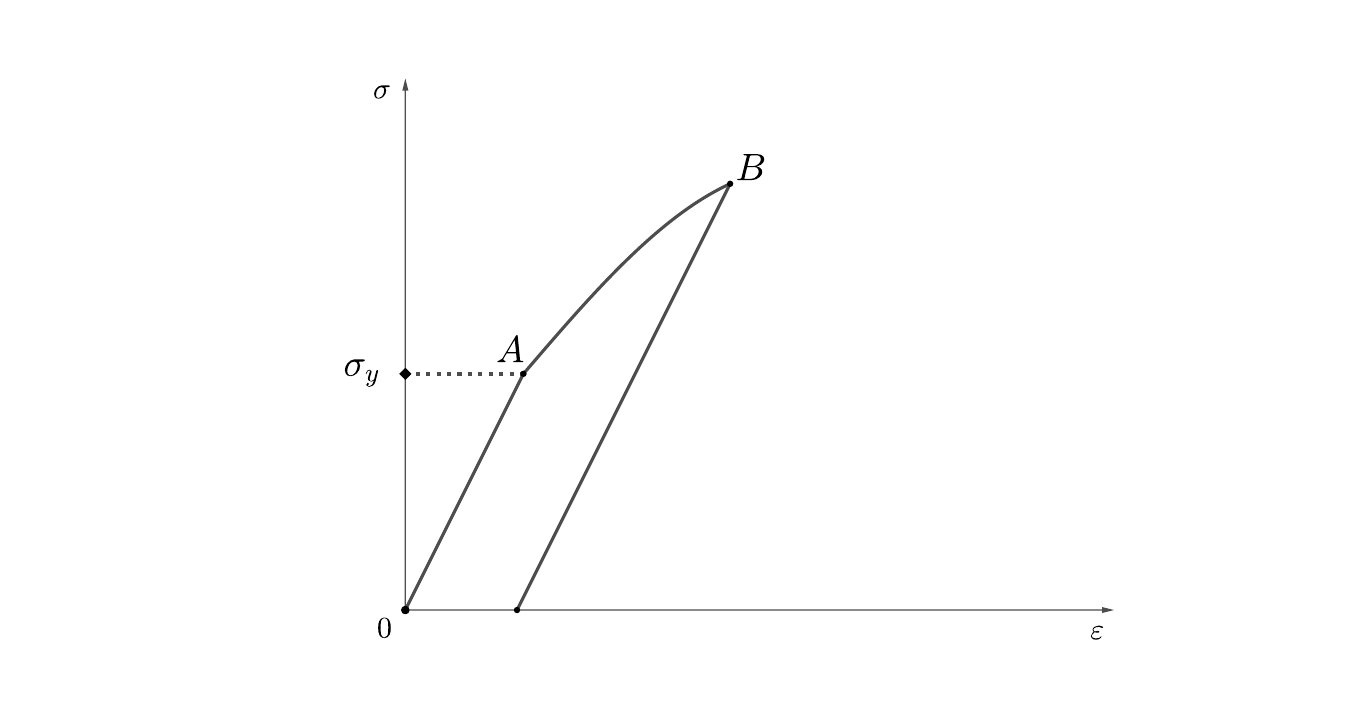}
  \caption{Elastic unloading.}
  \label{fig:elastic_unloading}
 \end{subfigure}
\caption{\it Stress versus strain for an elastoplastic material.}
\end{figure}

While for elastic materials the stress is completely determined by the strain and vice versa (and for linearly elastic ones this dependence is even linear) this one-to-one relation no longer holds true for \emph{elastoplastic} materials. To illustrate the behavior of such materials, consider an elastoplastic body $\Omega$ with uniaxial stress, i.e.~$\sigma_{11}$ is the only nonzero component of $\mymatrix{\sigma}$ and let $\sigma := \sigma_{11}$. In Figure~\ref{fig:stress_strain_relation}, the relation between stress and strain is plotted at a fixed material point $\myvec{x}\in\Omega$. More precisely, the graph shows the history of the stress $\sigma$ versus the strain $\varepsilon := \varepsilon_{11}$ during a process of loading.
%\vspace{0.1cm}

% -----------------------------------------------------------------------------------------

% -----------------------------------------------------------------------------------------

Successively increasing the force acting on the body will change its length and therefore causes a corresponding increase in strain. Up to a certain value $\sigma_y$ of stress -- the so-called \emph{initial yield stress} -- the material reacts in a linearly elastic fashion (see section $0A$ in Figure~\ref{fig:stress_strain_relation}). If then the force and therefore the stress is increased further, elastoplastic materials show a decrease in the slope of the curve. After that, various different phenomena can take place:
\begin{itemize}
 \item If the curve continues to rise with a slope less than that until the initial yield stress $\sigma_y$ the phenomena is known as \emph{hardening} (see section $AB$ in Figure~\ref{fig:stress_strain_relation}).
 \item If the slope of the curve becomes negative the behavior is called \emph{softening} (see section $BC$ in Figure~\ref{fig:stress_strain_relation}).
 \item If the slope starts to increase again the phenomena is called \emph{stiffening}.
\end{itemize}

In the papers [P1, P2, P3], a model problem of elastoplasticity with a certain type of hardening is considered. The transition from an elastic behavior (from a state of zero stress and strain until its limit $\sigma_y$) to a plastic one (at the point $A$ in Figure~\ref{fig:stress_strain_relation}) naturally entails a \emph{nonlinearity} in the mathematically description of elastoplastic problems. It is, however, not the nonlinearity that separates elastoplastic from elastic materials but the property of \emph{irreversibility}. In contrast to elastic materials, if the applied forces are removed the state of stress does not revert to its original state but decreases in an elastic fashion known as \emph{elastic unloading}, cf.~Figure~\ref{fig:elastic_unloading}. Thus, the one-to-one correspondence of stress and strain no longer holds true for elastoplastic materials. Furthermore, plastic behavior, in general, is \emph{rate-dependent}, i.e.~the material's response depends on the rate of a process. However, there are numerous materials that behave essentially \emph{rate-independent} for slow processes. Therefore, in the  elastoplasticity theory one usually requires rate-independence. Its extension to rate-dependent behavior is called \emph{viscoplasticity}, see e.g.~\cite{ref:Gurtin_2010, ref:Lemaitre_1990, ref:Lubliner_1972, ref:Simo_1998}.

% -----------------------------------------------------------------------------------------

\subsubsection*{The Plastic Strain and Internal Variables}

As the elastic and the plastic behavior completely differ at a microstructural level the total strain $\mymatrix{\varepsilon}$ can be decomposed additively into the \emph{elastic strain} $\mymatrix{e}$, which represents the elastic behavior of a material point and only depends on the stress $\mymatrix{\sigma}$, and the \emph{plastic strain} $\mymatrix{p}$, which representing the irreversible part of the deformation, i.e.~$\mymatrix{\varepsilon} = \mymatrix{e}(\mymatrix{\sigma}) + \mymatrix{p}$. If the elastic behavior of the material is linear the relation $\mymatrix{e} = \CC^{-1} \mymatrix{\sigma}$ holds true so that
\begin{align*}
 \mymatrix{\sigma} 
  = \CC \, \mymatrix{e} 
  = \CC \, (\mymatrix{\varepsilon} - \mymatrix{p})
\end{align*}
in that case. Moreover, from the behavior of metals it is observed that plastic deformations essentially do not cause a change of volume of a body and therefore it is generally required that a change of volume is exclusively caused by elastic deformations, which implies that $\tr \mymatrix{p} = 0$, see \cite{ref:Han_2013}.
\vspace{0.1cm}

In order to completely describe the behavior of elastoplastic materials in addition to the primary variables, which are given by the strain $\mymatrix{\varepsilon}$ (characterizing the local deformation) and the stress $\mymatrix{\sigma}$ being a quantity conjugate to the strain, so-called \emph{internal variables} $\mymatrix{\xi}_1,\ldots,\mymatrix{\xi}_m$ and \emph{internal forces} $\mymatrix{\chi}_1,\ldots,\mymatrix{\chi}_m$ have to be introduced, see \cite{ref:Coleman_1967, ref:Gurtin_1974, ref:Han_2013} for details on the notion of internal variables (particularly in the framework of thermoelasticity). Thereby, these internal variables and forces represent kinematic quantities and resulting forces that correspond to the internal restructuring while plastic deformation. Hardening, for instance, is characterized by such internal variables. While in general the plastic strain cannot be assumed to be one of the internal variables in the specific case of \emph{linearly kinematic hardening}, where hardening takes place at a constant rate, there is only one internal variable $\mymatrix{\xi}$ which is generally taken to be the plastic strain, i.e.~$\mymatrix{\xi} = \mymatrix{p}$. Furthermore, the corresponding internal force is given by
\begin{align*}
 \mymatrix{\chi} 
  = \HH \, \mymatrix{\xi}
  = \HH \, \mymatrix{p},
\end{align*}
where $\HH$ denotes the so-called \emph{hardening tensor} (comprising the material specific properties). As a consequence of the \emph{maximum plastic work inequality}, which follows from the \emph{postulate of maximum plastic work}, one obtains so-called \emph{plastic flow laws}, see e.g.~\cite{ref:Han_2013} for details on the maximum plastic work principle and the derivation of plastic flow laws for a general framework. In the specific case of \emph{linearly kinematic} hardening, the resulting flow law takes the form
\begin{align}\label{eq:plastic_work_flow}
 \mymatrix{\sigma} - \HH\, \mymatrix{p} \in \partial j(\dot{\mymatrix{p}}),
\end{align}
where $\dot{\mymatrix{p}} := \frac{\partial}{\partial t} \, \mymatrix{p}$ and $\partial j(\cdot)$ denotes the subdifferential of the \emph{dissipation function} $j(\cdot)$, which, for second-order tensors $\mymatrix{\tau}$, is given by $j(\mymatrix{\tau}) := \sigma_y \, \abs{\mymatrix{\tau}}_F$, cf.~\cite{ref:Wiedemann_2013}. In the following, the model problem of \emph{elastoplasticity with linearly kinematic hardening} is stated where, according to the papers [P1, P2, P3], the data will be assumed to be independent of time. Moreover, the time derivative of the plastic strain in \eqref{eq:plastic_work_flow} will be neglected so that the model problem represents one time step of \emph{quasi static time discrete elastoplasticity} with homogeneous initial condition $\mymatrix{p}_0 = \mymatrix{0}$, see~\cite{ref:Wiedemann_2013}. For the time discretizations, for instance, an implicit Euler scheme can be used.

% -----------------------------------------------------------------------------------------
%   The Model Problem of Elastoplasticity with Linearly Kinematic Hardening
% -----------------------------------------------------------------------------------------

\section{The Model Problem of Elastoplasticity with Linearly Kinematic Hardening}\label{sec:model_problem}

Let the reference configuration of the elastoplastic body, for $d\in\lbrace 2,3\rbrace$, be represented by the bounded domain $\Omega\in\RR^d$ with Lipschitz-boundary $\Gamma := \partial\Omega$ and outer unit normal $\myvec{n}$. Moreover, let the body be clamped at a \emph{Dirichlet boundary part} $\Gamma_D\subseteq\Gamma$, which corresponds to a given displacement $\overline{\myvec{u}}\equiv\myvec{0}$ on $\Gamma_D$. Then, for given volume force $\myvec{f} : \Omega\longrightarrow\RR^d$ and surface traction $\myvec{g} : \Gamma_N\longrightarrow\RR^d$ on the body and its \emph{Neumann boundary part} $\Gamma_N := \Gamma \setminus \overline{\Gamma_D}$, respectively, the boundary value problem of \emph{quasi static time discrete elastoplasticity with linearly kinematic hardening} is given by: Find a displacement field $\myvec{u} : \Omega\longrightarrow\RR^d$ and a plastic strain $\mymatrix{p} : \Omega\longrightarrow\SS_{d,0}$ that satisfy the \emph{equation of equilibrium}
\begin{align}\label{eq:PDE_model_problem}
 - \Div \big( \mymatrix{\sigma}(\myvec{u},\mymatrix{p}) \big) 
  &= \myvec{f} \qquad \text{in } \Omega,
\end{align}
the boundary conditions
\begin{align}\label{eq:BC_model_problem}
 \myvec{u} 
  = \myvec{0} \quad \text{on } \Gamma_D, \qquad
 \mymatrix{\sigma}(\myvec{u},\mymatrix{p}) \, \myvec{n} 
  = \myvec{g} \quad \text{on } \Gamma_N,
\end{align}
as well as the \emph{plastic flow law}
\begin{align}\label{eq:PFL_model_problem}
 \mymatrix{\sigma}(\myvec{u},\mymatrix{p}) - \HH \, \mymatrix{p}
  \in \partial j(\mymatrix{p}) \qquad \text{in } \Omega,
\end{align}
where the stress $\mymatrix{\sigma}$ is given by the \emph{constitutive relation} $\mymatrix{\sigma}(\myvec{u}) = \CC \, \big( \mymatrix{\varepsilon}(\myvec{u}) - \mymatrix{p} \big)$ with the elasticity tensor $\CC$, the strain tensor $\mymatrix{\varepsilon}$ is given by the \emph{strain-displacement relation} $\mymatrix{\varepsilon}(\myvec{u}) = \frac{1}{2} \, \big( \nabla\myvec{u} + (\nabla\myvec{u})^{\top} \big)$, $\HH$ denotes the hardening tensor, and $\SS_{d,0}$ is the space of symmetric $d\times d$ matrices over $\RR$ with vanishing trace, i.e.
\begin{align*}
 \SS_{d,0} 
  := \bigg\lbrace \mymatrix{\tau} = (\tau_{ij})\in\RR^{d\times d} \; ; \; \mymatrix{\tau} = \mymatrix{\tau}^{\top} \text{ and } \tr \mymatrix{\tau} = \sum_{i=1}^d \tau_{ii} = 0 \bigg\rbrace.
\end{align*}
Finally, the \emph{dissipation function} $j(\cdot)$ is given by $j(\mymatrix{\tau}) := \sigma_y \, \abs{\mymatrix{\tau}}_F$ for any $\mymatrix{\tau}\in\SS_{d,0}$.

% -----------------------------------------------------------------------------------------

\subsubsection*{Derivation of a Weak Formulation}

A classical solution of the boundary value problem \eqref{eq:PDE_model_problem}--\eqref{eq:PFL_model_problem} requires high smoothness assumptions on the primal variables $\myvec{u}$ and $\mymatrix{p}$, as well as on the data $\myvec{f}$, $\myvec{g}$, $\CC$ and might also be unrealistic from a physical point of view as mentioned earlier. Nevertheless, let a sufficiently smooth pair $(\myvec{u},\mymatrix{p})$ solve the boundary value problem \eqref{eq:PDE_model_problem}--\eqref{eq:PFL_model_problem}. Recall that for sufficiently smooth $\myvec{v}$ the identity
\begin{align*}
 \Div \big( \mymatrix{\sigma}(\myvec{u},\mymatrix{p}) \big) \cdot \myvec{v}
  = \Div \big( \mymatrix{\sigma}(\myvec{u},\mymatrix{p}) \, \myvec{v} \big) - \mymatrix{\sigma}(\myvec{u},\mymatrix{p}) : \nabla \myvec{v}
\end{align*}
holds true. Hence, by multiplying \eqref{eq:PDE_model_problem} with a smooth test function $\myvec{v}$ that vanishes on $\Gamma_D$, integrating over $\Omega$ and using the \emph{divergence theorem} of Gauss, one obtains the equation
\begin{align}\label{eq:varEQ_gradient}
 \int_{\Omega} \myvec{f} \cdot\myvec{v} \dd \myvec{x} 
  = \int_{\Omega} \mymatrix{\sigma}(\myvec{u},\mymatrix{p}) : \nabla\myvec{v} \dd \myvec{x} - \int_{\Gamma_N} \myvec{g} \cdot \myvec{v} \dd\myvec{s}.
\end{align}
As in \eqref{eq:varEQ_gradient} the gradient $\nabla\myvec{v}$ can be replaced by $\mymatrix{\varepsilon}(\myvec{v})$ (due to the symmetry of $\mymatrix{\sigma}$) any classical solution $(\myvec{u},\mymatrix{p})$ of \eqref{eq:PDE_model_problem}--\eqref{eq:PFL_model_problem} therefore satisfies the variational equation 
\begin{align}\label{eq:varEQ_01}
 \int_{\Omega} \mymatrix{\sigma}(\myvec{u},\mymatrix{p}) : \mymatrix{\varepsilon}(\myvec{v}) \dd \myvec{x}
  = \int_{\Omega} \myvec{f} \cdot\myvec{v} \dd \myvec{x} + \int_{\Gamma_N} \myvec{g} \cdot \myvec{v} \dd\myvec{s}
\end{align}
for all smooth test functions $\myvec{v}$ that vanish on the Dirichlet boundary part $\Gamma_D$. Moreover, the plastic flow law, \eqref{eq:PFL_model_problem}, yields the variational inequality
\begin{align}\label{eq:varIQ_01}
 \int_{\Omega} j(\mymatrix{\tau}) \dd\myvec{x}
  \geq \int_{\Omega} j(\mymatrix{p}) \dd\myvec{x} + \int_{\Omega} \big( \mymatrix{\sigma}(\myvec{u},\mymatrix{p}) - \HH \, \mymatrix{p} \big) : (\mymatrix{\tau} - \mymatrix{p}) \dd\myvec{x} 
\end{align}
for all integrable $\mymatrix{\tau} : \Omega\longrightarrow\SS_{d,0}$ by the definition of the subdifferential and using \emph{Riesz' representation theorem}.
\vspace{0.1cm}

Note that for the equations \eqref{eq:varEQ_01} and \eqref{eq:varIQ_01} to make sense it is sufficient to require $\myvec{u},\myvec{v}\in H^1(\Omega,\RR^d)$ and $\mymatrix{p}, \mymatrix{q}\in L^2(\Omega,\SS_{d,0})$, as well as $\myvec{f}\in L^2(\Omega,\RR^d)$, $\myvec{g}\in L^2(\Gamma_N,\RR^d)$ and $\sigma_y, \CC_{ijkl}, \HH_{ijkl}\in L^{\infty}(\Omega)$ for the given data. Here, $L^2(\Omega,\RR^d)$ and $L^2(\Omega,\SS_{d,0})$ are the space of vector-valued functions $\myvec{v} = (v_1,\ldots,v_d)^{\top}$ with components $v_i$ in $L^2(\Omega)$ and matrix-valued functions $\mymatrix{q}=(q_{ij})$  with components $q_{ij}$ in $L^2(\Omega)$, respectively, which are endowed with the inner products
\begin{align*}
 (\myvec{v},\myvec{w})_{0,\Omega}
  := \int_{\Omega} \myvec{v} \cdot \myvec{w} \dd \myvec{x}
  = \int_{\Omega} \sum_{i=1}^d v_i(\myvec{x}) \, w_i(\myvec{x}) \dd \myvec{x}, \qquad
 (\mymatrix{p},\mymatrix{q})_{0,\Omega}
  := \int_{\Omega} \mymatrix{p} : \mymatrix{q} \dd \myvec{x}
  = \int_{\Omega} \sum_{i=1}^d \sum_{j=1}^d p_{ij}(\myvec{x}) \, q_{ij}(\myvec{x}) \dd \myvec{x},
\end{align*}
respectively. Furthermore, let $\norm{\cdot}_{0,\Omega}$ denote in both cases the corresponding norm. Analogously, $H^1(\Omega,\RR^d)$ is the space of vector-valued functions $\myvec{v}$ with components $v_i\in H^1(\Omega)$, where $H^1(\Omega)$ is the usual Sobolev space of all functions in $L^2(\Omega)$ having weak first-order partial derivatives in $L^2(\Omega)$. If $\norm{\cdot}_{H^1(\Omega)}$ represents the usual Sobolev norm,
\begin{align*}
 \norm{\myvec{v}}_{H^1(\Omega,\RR^d)}
  := \bigg( \sum_{i=1}^d \norm{v_i}_{H^1(\Omega)}^2 \bigg)^{1/2}
\end{align*}  
is a corresponding norm in $H^1(\Omega,\RR^d)$. Furthermore, for real $s\geq 0$ with $s = m + \lambda$ for some $m\in\NN_0$ and $0<\lambda<1$ let $H^s(\Omega)$ be the space of all functions $v\in L^2(\Omega)$ for which $\norm{v}_s < \infty$, where $\norm{\cdot}_s$ denotes the \emph{Sobolev-Slobodeckij-norm},  see e.g.~\cite{ref:Hackbusch_2012}. Finally, let $H^s(\Omega,\RR^d)$ be the space of vector-valued functions $\myvec{v}$ with components $v_i\in H^s(\Omega)$. In order that $\myvec{u}$ and $\myvec{v}$ satisfy the homogeneous Dirichlet boundary condition they have to be taken from the Hilbert function space
\begin{align*}
 \VV := \big\lbrace \myvec{v}\in H^1(\Omega,\RR^d) \; ; \; \gamma\myvec{v} = \myvec{0} \text{ on } \Gamma_D \big\rbrace,
\end{align*}
where $\gamma : H^1(\Omega,\RR^d)\longrightarrow H^{1/2}(\Gamma_D,\RR^d)$ is the unique trace operator, see~\cite{ref:Hackbusch_2017, ref:Han_2013}. If in the following the trace $\gamma\myvec{v}$ of some function $\myvec{v}\in H^1(\Omega,\RR^d)$ is well defined on some boundary part it is simply written $\myvec{v}$ for $\gamma\myvec{v}$. 
\vspace{0.1cm}

By interpreting the integrals of the right-hand side of \eqref{eq:varEQ_01} as the duality pairing between $\VV$ and its dual space $\VV^{\star}$ and between the trace-space of $\VV$ restricted to $\Gamma_N$ and its dual space $H^{-1/2}(\Gamma_N,\RR^d)$, respectively, a classical solution $(\myvec{u},\mymatrix{p})$ to the boundary value problem \eqref{eq:PDE_model_problem}--\eqref{eq:PFL_model_problem} satisfies the variational equation
\begin{align}\label{eq:varEQ_02}
 \big( \mymatrix{\sigma}(\myvec{u},\mymatrix{p}), \mymatrix{\varepsilon}(\myvec{v}) \big)_{0,\Omega}
  = \langle \myvec{f}, \myvec{v}\rangle + \langle\myvec{g}, \myvec{v} \rangle_{\Gamma_N} \qquad
 \forall \, \myvec{v}\in\VV
\end{align}
as well as the variational inequality
\begin{align}\label{eq:varIQ_02}
 \big( \mymatrix{\sigma}(\myvec{u},\mymatrix{p}) - \HH \, \mymatrix{p}, \mymatrix{p} - \mymatrix{q} \big)_{0,\Omega} + \int_{\Omega} j(\mymatrix{q}) \dd\myvec{x} - \int_{\Omega} j(\mymatrix{p}) \dd \myvec{x}
  \geq 0 \qquad
 \forall \, \mymatrix{q}\in Q,
\end{align}
where $Q := L^2(\Omega,\SS_{d,0})$. By defining the bilinear form $a : \VV\times \VV\longrightarrow\RR$, the so-called \emph{plasticity functional} $\psi : Q\longrightarrow\RR$ and the linear form $\ell : \VV\longrightarrow\RR$ as
\begin{align*}
 a\big( (\myvec{v},\mymatrix{q}), (\myvec{w},\mymatrix{\tau}) \big)
  &:= \big( \mymatrix{\sigma}(\myvec{v},\mymatrix{q}), \mymatrix{\varepsilon}(\myvec{w}) - \mymatrix{\tau}) \big)_{0,\Omega} + ( \HH \, \mymatrix{q}, \mymatrix{\tau} )_{0,\Omega}, \\
 \psi(\mymatrix{q})
  &:= ( \sigma_y, \abs{\mymatrix{q}}_F )_{0,\Omega}, \\
 \ell(\myvec{v})
  &:= \langle \myvec{f}, \myvec{v}\rangle + \langle\myvec{g}, \myvec{v} \rangle_{\Gamma_N},
\end{align*}
respectively, one finds that a pair $(\myvec{u}, \mymatrix{p})$ satisfies the variational equation \eqref{eq:varEQ_02} \emph{and} the variational inequality \eqref{eq:varIQ_02} if and only if $(\myvec{u},\mymatrix{p})$ satisfies the \emph{variational inequality of the second kind}
\begin{align}\label{eq:varIQ_of_second_kind}
 a\big( (\myvec{u},\mymatrix{p}), (\myvec{v} - \myvec{u},\mymatrix{q} - \mymatrix{p}) \big) + \psi(\mymatrix{q}) - \psi(\mymatrix{p})
  \geq \ell(\myvec{v} - \myvec{u}) \qquad
 \forall \, (\myvec{v},\mymatrix{q})\in\VV\times Q.
\end{align}

Any classical solution $(\myvec{u},\mymatrix{p})$ to the boundary value problem \eqref{eq:PDE_model_problem}--\eqref{eq:PFL_model_problem} satisfies the variational inequality \eqref{eq:varIQ_of_second_kind}, which therefore represents a \emph{weak formulation} of \eqref{eq:PDE_model_problem}--\eqref{eq:PFL_model_problem}. Conversely, any \emph{weak solution} $(\myvec{u},\mymatrix{p})\in \VV\times Q$ being sufficiently smooth so that the arguments leading to the variational formulation can be reversed solves the boundary value problem \eqref{eq:PDE_model_problem}--\eqref{eq:PFL_model_problem}. In this sense, the classical formulation \eqref{eq:PDE_model_problem}--\eqref{eq:PFL_model_problem} of the model problem and the variational inequality \eqref{eq:varIQ_of_second_kind} are equivalent.

% -----------------------------------------------------------------------------------------

\subsubsection*{Existence of a Weak Formulation}

In order to guarantee the existence of a weak solution $(\myvec{u},\mymatrix{p})\in\VV\times Q$ of the variational inequality \eqref{eq:varIQ_of_second_kind} one has to require some properties of the elasticity tensor $\CC$ and the hardening tensor $\HH$. First, they have to be symmetric, i.e.~$\CC_{ijkl} = \CC_{jilk} = \CC_{klij}$ and $\HH_{ijkl} = \HH_{jilk} = \HH_{klij}$ for $i,j,k,l\in\underline{d}$ and, secondly, uniformly elliptic, i.e.~there exist constants $c_e,c_h>0$ such that
\begin{align*}
 (\CC \, \mymatrix{\tau}) : \mymatrix{\tau} 
  \geq c_e \, \abs{\mymatrix{\tau}}_F^2, \qquad
 (\HH \, \mymatrix{\tau}) : \mymatrix{\tau} 
  \geq c_h \, \abs{\mymatrix{\tau}}_F^2 \qquad
 \forall \, \mymatrix{\tau}\in Q.
\end{align*}
According to the papers [P1, P2, P3], the \emph{yield stress} $\sigma_y$ in uniaxial tension is assumed to be a positive constant $\sigma_y > 0$. Note that $\VV\times Q$ forms a Hilbert-space, which can be equipped with the norm $\norm{\cdot}$, defined by $\norm{(\myvec{v},\mymatrix{q})} := \big( \norm{\myvec{v}}_{1,\Omega}^2 + \norm{\mymatrix{q}}_{0,\Omega}^2 \big)^{1/2}$ for any $(\myvec{v},\mymatrix{q})\in\VV\times Q$ with $\norm{\myvec{v}}_{1,\Omega} := \big( \norm{\myvec{v}}_{0,\Omega}^2 + \norm{\mymatrix{\varepsilon}(\myvec{v})}_{0,\Omega}^2 \big)^{1/2}$. Then, $a(\cdot,\cdot)$ represents a symmetric, continuous and $(\VV\times Q)$-elliptic bilinear form, see the Appendix for the details. Furthermore, the plasticity functional $\psi(\cdot)$ is convex, Lipschitz-continuous with constant $\norm{\sigma_y}_{0,\Omega}$ and sub-differentiable, cf.~[P2]. Thereby, the convexity and Lipschitz-continuity immediately follow from the convexity of the Frobenius norm $\abs{\cdot}_F$ and by the Cauchy-Schwarz inequality, respectively. Therewith, the \emph{energy functional} $\E : \VV\times Q\longrightarrow\RR$, defined as
\begin{align}\label{eq:defi_energyF}
 \E(\myvec{v},\mymatrix{q})
  := \frac{1}{2} \, a\big( (\myvec{v},\mymatrix{q}), (\myvec{v},\mymatrix{q}) \big) + \psi(\mymatrix{q}) - \ell(\myvec{v}),
\end{align}
is coercive, convex and subdifferentiable, see the Appendix, and therefore, the minimization problem
\begin{align}\label{eq:minProb}
 \E(\myvec{u},\mymatrix{p})
  \leq \E(\myvec{v},\mymatrix{q}) \qquad
 \forall \, (\myvec{v},\mymatrix{q})\in\VV\times Q
\end{align}
has a unique solution $(\myvec{u},\mymatrix{p})\in\VV\times Q$ by \cite[Ch.~II, Prop.~1.2]{ref:Ekeland_1999}. Since a pair $(\myvec{u},\mymatrix{p})\in\VV\times Q$ is a minimizer of \eqref{eq:minProb} if and only if, $(\myvec{u},\mymatrix{p})$ solves the variational inequality of the second kind \eqref{eq:varIQ_of_second_kind}, see e.g.~\cite{ref:Han_1991, ref:Han_2013}, the existence of a weak solution is guaranteed.
\vspace{0.1cm}

Moreover, the $(\VV\times Q)$-ellipticity of $a(\cdot,\cdot)$ immediately yields the uniqueness of a weak solution $(\myvec{u},\mymatrix{p})\in\VV\times Q$ of the variational inequality \eqref{eq:varIQ_of_second_kind}. For, if $(\myvec{u}_1,\mymatrix{p}_1),(\myvec{u}_2,\mymatrix{p}_2)\in\VV\times Q$ are two solutions of \eqref{eq:varIQ_of_second_kind}, it follows that
\begin{align*}
 a\big( (\myvec{u}_1,\mymatrix{p}_1), (\myvec{u}_2-\myvec{u}_1, \mymatrix{p}_2-\mymatrix{p}_1) \big) + \psi(\mymatrix{p}_2) - \psi(\mymatrix{p}_1) 
 &\geq \ell(\myvec{u}_2-\myvec{u}_1), \\
 a\big( (\myvec{u}_2,\mymatrix{p}_2), (\myvec{u}_1-\myvec{u}_2, \mymatrix{p}_1-\mymatrix{p}_2) \big) + \psi(\mymatrix{p}_1) - \psi(\mymatrix{p}_2) 
 &\geq \ell(\myvec{u}_1-\myvec{u}_2).
\end{align*}
Hence, adding these two inequalities gives
\begin{align*}
 0 
  &\leq a\big( (\myvec{u}_2,\mymatrix{p}_2), (\myvec{u}_1-\myvec{u}_2, \mymatrix{p}_1-\mymatrix{p}_2) \big) + a\big( (\myvec{u}_1,\mymatrix{p}_1), (\myvec{u}_2-\myvec{u}_1, \mymatrix{p}_2-\mymatrix{p}_1) \big) \\
 &= - a\big( (\myvec{u}_2-\myvec{u}_1, \mymatrix{p}_2-\mymatrix{p}_1), (\myvec{u}_2-\myvec{u}_1, \mymatrix{p}_2-\mymatrix{p}_1) \big).
\end{align*}
Exploiting the $(\VV\times Q)$-ellipticity of $a(\cdot,\cdot)$ therefore yields
\begin{align*}
 0\leq 
\alpha \, \norm{ (\myvec{u}_2-\myvec{u}_1, \mymatrix{p}_2-\mymatrix{p}_1) }
 \leq a\big( (\myvec{u}_2-\myvec{u}_1, \mymatrix{p}_2-\mymatrix{p}_1), (\myvec{u}_2-\myvec{u}_1, \mymatrix{p}_2-\mymatrix{p}_1) \big) 
 \leq 0,
\end{align*}
from which follows that $\norm{ (\myvec{u}_2-\myvec{u}_1, \mymatrix{p}_2-\mymatrix{p}_1) } = 0$. Thus, the identity $(\myvec{u}_1,\mymatrix{p}_1) = (\myvec{u}_2,\mymatrix{p}_2)$ holds true.

% -----------------------------------------------------------------------------------------
%   CHAPTER :  Contributions of this Thesis
% -----------------------------------------------------------------------------------------

\chapter{Contribution of this Thesis}

While in the papers [P1, P2, P3], representing the first part of this thesis, the numerical analysis of different $hp$-finite element discretizations of a model problem of elastoplasticity with linearly kinematic hardening is studied, a new $hp$-adaptive algorithm for solving variational equations, which \emph{does not} rely on the use of a classical a posteriori error estimator, is introduced in [P4].

% -----------------------------------------------------------------------------------------
%   Part I : $hp$-Finite Element Method for Elastoplasticity
% -----------------------------------------------------------------------------------------

\section{Part I :\, $hp$-Finite Element Method for Elastoplasticity}

As already mentioned in Cahpter~\ref{chapter:introduction}, elastoplasticity with hardening appears in various problems of mechanical engineering, for instance, when modeling the deformation of concrete or metals, see e.g.~\cite{ref:Chen_1988}. Thereby, the specific case of \emph{linearly kinematic hardening} plays an important role. The \emph{holonomic constitutive law} represents a well-established case of elastoplasticity with linearly kinematic hardening, which allows for the incremental computation of the deformation of an elastoplastic body, see e.g.~\cite{ref:Han_1991, ref:Han_1995, ref:Han_2013}. A well-established weak formulation of one (pseudo-)time step of that problem is given by the variational inequality of the second kind, already introduced in Section~\ref{sec:model_problem}, see e.g.~\cite{ref:Carstensen_1999, ref:Han_1991, ref:Han_2013}. This variational inequality, however, contains the non-differentiable \emph{plasticity functional} $\psi(\cdot)$, which causes many difficulties not only in the numerical analysis of a discretization but also in the numeric. One possibility to resolve the non-smoothness of $\psi(\cdot)$ is to apply a suitable regularization, as proposed, for instance, in \cite{ref:Reddy_1987}. The influence of such regularizations can indeed again lead to problems in the numerical analysis and the computation of a discrete solution. In [P1, P3], we present a discretization of the variational inequality \eqref{eq:varIQ_of_second_kind}, in which the Frobenius norm involved in the definition of the plasticity functional $\psi(\cdot)$ is approximated by an appropriate interpolation, as proposed in \cite{ref:Gwinner_2013_2} (within the framework of Tresca friction). This leads to a discrete plasticity functional $\psi_{hp}(\cdot)$, which can exactly be evaluated by a suitable, easy to implement quadrature rule.
\vspace{0.1cm}

Another way to circumvent the difficulties resulting from the non-differentiability of $\psi(\cdot)$ is to reformulate the variational inequality \eqref{eq:varIQ_of_second_kind} as a \emph{mixed formulation}, in which the non-smoothness of $\psi(\cdot)$ is resolved by introducing a Lagrange multiplier, see e.g.~\cite{ref:Han_1991, ref:Han_1995, ref:Schroeder_2011_2} (or the works \cite{ref:Schroeder_2011_4, ref:Schroeder_2011_3} on mixed methods in the similar framework of frictional contact problems). The discrete mixed formulation can then be solved, for instance, by an Uzawa method as proposed in \cite{ref:Schroeder_2011_2} or by the semi-smooth Newton method, presented in [P1]. For a conforming lowest-order finite element method (in particular, conforming in the discrete Lagrange multiplier), see e.g.~\cite{ref:Schroeder_2011_2}, and for a higher-order method, which is non-conforming in the discrete Lagrange multiplier, see~\cite{ref:Wiedemann_2013}. The higher-order mixed finite element methods, presented in [P1, P2, P3], are conforming in the displacement field $\myvec{u}$ and the plastic strain variable $\mymatrix{p}$ but non-conforming in the discrete Lagrange multiplier (except for the lowest-order case). Thereby, the non-conformity in the Lagrange multiplier is necessary to obtain an implementable discretization scheme. Though the use of the discrete Lagrange multiplier as an additional variable naturally leads to a substantial increase of the number of degrees of freedom (as it contains the same number of unknowns as the plastic strain variable), it should be mentioned that the same polynomial degree distribution as well as the same basis functions can be used for $\mymatrix{p}$ and the Lagrange multiplier, which limits the additional effort relating to the implementation. Furthermore, the presented $hp$-discretizations rely on the \emph{same} mesh for all three variables, respectively (i.e.~for the displacement field, the plastic strain and the discrete Lagrange multiplier).
\vspace{0.1cm}

The use of \emph{biorthogonal} basis functions for the discretization of the plastic strain and the Lagrange multiplier allows to show the equivalence between the discrete variational inequality and the discrete mixed formulation, presented in [P1, P3], see [P3] for the proof. Under a slight limitation on the physical elements' shape these two discretizations turn out to be equivalent to a third one, which is again based on the mixed formulation, see [P2, P3]. In this case, the a priori results in [P2], in particular, the \emph{convergence} and the \emph{guaranteed convergence rates} in the mesh size $h$ and polynomial degree $p$, can be applied to all three $hp$-finite element discretizations. It should be mentioned that the non-conformity of the discrete Lagrange multiplier causes a reduction of the guaranteed convergence rates, which is, however, typical for higher-order mixed methods for variational inequalities, see e.g.~\cite{ref:Banz_2022, ref:Banz_2019, ref:Ovcharova_2017}. Furthermore, the use of the biorthogonal basis functions enables to \emph{decouple} the constraints associated with the discrete Lagrange multiplier, which offers the possibility to reformulate the mixed formulation in terms of a system of decoupled nonlinear equations. This, in fact, simplifies the application of solution scheme, in particular, of the semi-smooth Newton solver, proposed in [P1]. Thereby, its applicability as well as its robustness to $h$, $p$ and the projection parameters are shown in the numerical examples in [P1]. For the application of a semi-smooth Newton solver in the context of elastoplastic (contact) problems, see also~\cite{ref:Christensen_2002, ref:Hintermueller_2016}.
\vspace{0.1cm}

As already mentioned in Section~\ref{sec:FEM}, the weak solution of problems in elastoplasticity does not enjoy high regularity properties and typically contains singularities along the free boundary separating the regions of purely elastic deformation from those of plastic deformation. In order to achieve possibly high convergence rates, one therefore has to apply $h$- or $hp$-adaptive finite element methods. Thereby, \emph{a posteriori error control} is an essential tool to steer adaptive refinements, which usually relies on the derivation of a reliable and efficient a posteriori error estimator, see \cite{ref:Ainsworth_1996, ref:verfuerth_2013}. In this context, \emph{reliability} means that the discretization error $\norm{ u_\XX - u_\WW}_\XX$ is bounded from above by the error estimator up to a multiplicative constant and some higher-order terms, whereas \emph{efficiency} is on hand if the reverse holds true, i.e.~if the error estimator is bounded from above by the discretization error up to a multiplicative constant and some higher-order terms. Error control approaches for lowest-order finite element methods for problems of elastoplasticity with hardening can be found, e.g.~in \cite{ref:Alberty_1999, ref:Carstensen_1999, ref:Carstensen_2003, ref:Carstensen_2006_2, ref:Schroeder_2011_2}, and for an optimally converging adaptive finite element method in this context, see \cite{ref:Carstensen_2016}. In [P3], a residual-based a posteriori error estimator for the model problem is derived from upper and lower error estimates relying on an auxiliary problem, which takes the form of a variational equation. Furthermore, its reliability and some (local) efficiency estimates are shown. Thereby, the efficiency estimates are suboptimal for higher-order methods, which is, however, expectable within the considered framework. The error estimator is employable to any discretization that is conforming with respect to the displacement field and the plastic strain and therefore can be applied to all three $hp$-finite element discretizations presented in [P1, P2, P3]. The numerical experiments in [P3] definitely show the potential of $h$- and $hp$-adaptivity for problems in elastoplasticity. We thereby observe, that the finite element spaces corresponding to the adaptive methods are adapted to the singularities of the solution, in particular, to those of the free boundary. \vspace{0.3cm} \\

In order to give an overview of the main results in the papers [P1, P2, P3], the most important theorems and their connections are presented in the following paragraphs. For further details, in particular, the proofs, see [P1, P2, P3]. \\

% -----------------------------------------------------------------------------------------

\subsubsection*{A Mixed Variational Formulation}

The common weak formulation of the model problem \eqref{eq:PDE_model_problem}--\eqref{eq:PFL_model_problem} of \emph{quasi static time discrete elastoplasticity with linearly kinematic hardening} as the variational inequality of the second kind
\begin{align}\label{eq:weakF_VIQ}
 a\big( (\myvec{u},\mymatrix{p}), (\myvec{v} - \myvec{u},\mymatrix{q} - \mymatrix{p}) \big) + \psi(\mymatrix{q}) - \psi(\mymatrix{p})
  \geq \ell(\myvec{v} - \myvec{u}) \qquad
 \forall \, (\myvec{v},\mymatrix{q})\in\VV\times Q
\end{align}
has already been derived in Section~\ref{sec:model_problem}. In addition, we pointed out that there exists a unique weak solution $(\myvec{u},\mymatrix{p})\in\VV\times Q$ of \eqref{eq:weakF_VIQ}. By introducing the nonempty, closed and convex set
\begin{align}\label{eq:representation_Lambda01}
  \Lambda 
 := \big\lbrace \mymatrix{\mu}\in Q \; ; \; \vert\mymatrix{\mu}\vert_F \leq \sigma_y \text{ a.e.~in }\Omega \big\rbrace,
\end{align}
which can equivalently be represented as
\begin{align}\label{eq:representation_Lambda02}
 \Lambda
 = \big\lbrace \mymatrix{\mu}\in Q \; ; \; (\mymatrix{\mu}, \mymatrix{q})_{0,\Omega} \leq \psi(\mymatrix{q}) \text{ for all } \mymatrix{q}\in Q \big\rbrace,
\end{align}
see [P2, Sec.~3], it follows that $\psi(\mymatrix{q}) = \sup_{\mymatrix{\mu}\in\Lambda} (\mymatrix{\mu}, \mymatrix{q})_{0,\Omega}$ for any $\mymatrix{q}\in Q$. Thus, a \emph{mixed variational formulation} of the boundary value problem \eqref{eq:PDE_model_problem}--\eqref{eq:PFL_model_problem} is given by: Find a triple $(\myvec{u},\mymatrix{p},\mymatrix{\lambda})\in \VV\times Q\times \Lambda$ such that
\begin{subequations}\label{eq:mixed_variationalF}
\begin{alignat}{2}
a\big( (\myvec{u},\mymatrix{p}),(\myvec{v},\mymatrix{q}) \big) + (\mymatrix{\lambda}, \mymatrix{q})_{0,\Omega} &= \ell(\myvec{v})
& \qquad &\forall \, (\myvec{v},\mymatrix{q})\in \VV \times Q, \label{eq:mixed_variationalF_01}\\
(\mymatrix{\mu}-\mymatrix{\lambda}, \mymatrix{p})_{0,\Omega} &\leq 0
& \qquad &\forall \, \mymatrix{\mu}\in\Lambda. \label{eq:mixed_variationalF_02}
\end{alignat}
\end{subequations}
Thereby, the unique existence of a solution $(\myvec{u},\mymatrix{p},\mymatrix{\lambda})\in \VV\times Q\times \Lambda$ of the mixed variational problem \eqref{eq:mixed_variationalF} is guaranteed by the following result.

% -----------------------------------------------------------------------------------------
\begin{theorem}[{[P2, Thm.~1]}]
 If $(\myvec{u},\mymatrix{p})\in \VV\times Q$ solves \eqref{eq:weakF_VIQ}, then $(\myvec{u},\mymatrix{p},\mymatrix{\lambda})$ with
\begin{align}\label{eq:repres_continuousLambda}
 \mymatrix{\lambda} 
 = \operatorname{dev}\big( \mymatrix{\sigma}(\myvec{u},\mymatrix{p})-\HH \, \mymatrix{p} \big)
\end{align}
is a solution to \eqref{eq:mixed_variationalF} and, conversely, if $(\myvec{u},\mymatrix{p},\mymatrix{\lambda})\in \VV\times Q \times \Lambda$ solves \eqref{eq:mixed_variationalF}, then $(\myvec{u},\mymatrix{p})$ is a solution of \eqref{eq:weakF_VIQ} and the identity \eqref{eq:repres_continuousLambda} holds true.
\end{theorem}
% -----------------------------------------------------------------------------------------

Furthermore, it is shown in [P2, Lem.~3] that the solution  $(\myvec{u},\mymatrix{p},\mymatrix{\lambda})\in \VV\times Q\times \Lambda$ of \eqref{eq:mixed_variationalF} depends Lipschitz-continuously on the data $\myvec{f}$, $\myvec{g}$ and $\sigma_y$. More precisely,
\begin{align*}%\label{eq:Lipschitz_depend}
 \norm{ (\myvec{u}_2 - \myvec{u}_1, \mymatrix{p}_2 - \mymatrix{p}_1) } + \norm{ \mymatrix{\lambda}_2 -\mymatrix{\lambda}_1 }_{0,\Omega}  
 \lesssim  \norm{ \sigma_{y,2} - \sigma_{y,1} }_{0,\Omega} + \norm{ \myvec{f}_2 - \myvec{f}_1 }_{\VV^{\star}} + \norm{ \myvec{g}_2 - \myvec{g}_1 }_{H^{-1/2}(\Gamma_N,\RR^d)}, 
\end{align*}
where $(\myvec{u}_i,\mymatrix{p}_i,\mymatrix{\lambda}_i)$, for $i=1,2$, denotes the discrete solution according to the data $(\myvec{f}_i,\myvec{g}_i,\sigma_{y,i})$, respectively.

% -----------------------------------------------------------------------------------------

\subsubsection*{$hp$-Finite Element Discretizations}

In order to present an $hp$-finite element discretization of the weak formulations \eqref{eq:weakF_VIQ} and \eqref{eq:mixed_variationalF}, respectively, let $\T_{h}$ be a locally quasi-uniform finite element mesh of $\Omega$ consisting of convex and shape regular quadrilaterals or hexahedrons. Moreover, let $\widehat{T}:=[-1,1]^d$ be the reference element and $\myvec{F}_T:\widehat{T}\longrightarrow T$ denote the bi/tri-linear bijective mapping for $T\in\T_h$. We set $h := (h_T)_{T\in\T_h}$ and $p := (p_T)_{T\in\T_h}$ where $h_T$ and $p_T$ denote the local element size and polynomial degree, respectively. We assume that the local polynomial degrees of neighboring elements are comparable in the sense of \cite{ref:Melenk_2005}. For the discretization of the displacement field $\myvec{u}$ and of the plastic strain $\mymatrix{p}$ in all three papers [P1, P2, P3] we use the $hp$-finite element spaces
\begin{align*}
 V_{hp} &:= \Big\lbrace \myvec{v}_{hp}\in \VV \; ; \; \myvec{v}_{hp \, | \, T}\circ \myvec{F}_T\in \big(\PP_{p_T}(\widehat{T})\big)^d \text{ for all } T\in\T_h\Big\rbrace, \\
 Q_{hp} &:= \Big\lbrace \mymatrix{q}_{hp}\in Q \; ; \; \mymatrix{q}_{hp \, | \, T}\circ \myvec{F}_T\in \big(\PP_{p_T-1}(\widehat{T})\big)^{d\times d} \text{ for all } T\in\T_h\Big\rbrace.
\end{align*}
Furthermore, let $\hat{\myvec{x}}_{k,T}\in \widehat{T}$ for $k \in \ul{n_T}$ be the tensor product Gauss quadrature points on $\widehat{T}$ with corresponding positive weights $\hat{\omega}_{k,T}\in\RR$, where $n_T := p_T^d$ for $T\in\T_h$. Thereby, we introduce the quadrature rule
\begin{align}\label{eq:quadrature_rule}
 \Q_{hp}(\cdot) := \sum_{T\in\T_h} \Q_{hp, T}(\cdot),
\end{align}
where the local quantities $\Q_{hp, T}(\cdot)$ are given by
\begin{align*}
 \Q_{hp, T}(f) :=
 \begin{cases}
  |T| \, f\big(\myvec{F}_T(\myvec{0})\big), & \text{if } p_T = 1, \\
  \sum_{k=1}^{n_T} \hat{\omega}_{k,T} \, |\det \nabla \myvec{F}_T(\hat{\myvec{x}}_{k,T})| \, f\big(\myvec{F}_T(\hat{\myvec{x}}_{k,T})\big), & \text{if } p_T \geq 2,
 \end{cases}
 \qquad T\in\T_h.
\end{align*}
In order to obtain a discretization of \eqref{eq:weakF_VIQ}, in which the non-differentiable plasticity functional $\psi(\cdot)$ is approximated, we interpolate the Frobenius norm involved in the definition of $\psi(\cdot)$ by a nodal interpolation operator, see [P1, Sec.~2; P2, Sec.~5]. In fact, the resulting \emph{discrete plasticity functional} $\psi_{hp} : Q_{hp}\longrightarrow\RR$ can exactly be evaluated by the quadrature rule \eqref{eq:quadrature_rule}, i.e.
\begin{align*}
 \psi_{hp}(\mymatrix{q}_{hp})
  = \Q_{hp}\big( \psi(\mymatrix{q}_{hp}) \big),
\end{align*}
cf.~[P2, Sec.~5]. Thus a \emph{discrete variational inequality of the second kind}, approximating the variational inequality \eqref{eq:weakF_VIQ}, is given by: Find a pair $(\myvec{u}_{hp}, \mymatrix{p}_{hp})\in V_{hp}\times Q_{hp}$ such that
\begin{align}\label{eq:discrete_variq_second_kind}
 a\big( (\myvec{u}_{hp},\mymatrix{p}_{hp}), (\myvec{v}_{hp}-\myvec{u}_{hp}, \mymatrix{q}_{hp}-\mymatrix{p}_{hp}) \big) + \psi_{hp}(\mymatrix{q}_{hp}) - \psi_{hp}(\mymatrix{p}_{hp}) 
 \geq \ell(\myvec{v}_{hp} - \myvec{u}_{hp})
 \quad \forall \, (\myvec{v}_{hp},\mymatrix{q}_{hp})\in V_{hp}\times Q_{hp}.
\end{align}
Thereby, the next result guarantees the existence of a unique discrete solution of \eqref{eq:discrete_variq_second_kind}.

% -----------------------------------------------------------------------------------------
\begin{theorem}[{[P3, Thm.~4]}]\label{thm:uniqueEx_VIQ}
The discrete variational inequality \eqref{eq:discrete_variq_second_kind} has a unique solution $(\myvec{u}_{hp},\mymatrix{p}_{hp})\in V_{hp}\times Q_{hp}$.
\end{theorem}
% -----------------------------------------------------------------------------------------

For the discretization of the mixed variationa problem \eqref{eq:mixed_variationalF} one needs to choose a nonempty, convex and closed set of \emph{admissible discrete Lagrange multipliers}. A discretization of the representation \eqref{eq:representation_Lambda02} is given by
\begin{align*}
 \Lambda_{hp}^{(w)} 
  := \big\lbrace \mymatrix{\mu}_{hp}\in Q_{hp} \; ; \; (\mymatrix{\mu}_{hp},\mymatrix{q}_{hp})_{0,\Omega} \leq \psi_{hp}(\mymatrix{q}_{hp}) \text{ for all } \mymatrix{q}_{hp}\in Q_{hp}  \big\rbrace,
\end{align*}
leading to the \emph{discrete mixed formulation}: Find a triple $(\myvec{u}_{hp},\mymatrix{p}_{hp},\mymatrix{\lambda}_{hp})\in V_{hp}\times Q_{hp}\times\Lambda_{hp}^{(w)}$ such that
\begin{subequations}\label{eq:discrete_mixedF_01}
\begin{alignat}{2}
 a\big( (\myvec{u}_{hp},\mymatrix{p}_{hp}),(\myvec{v}_{hp},\mymatrix{q}_{hp}) \big) + (\mymatrix{\lambda}_{hp}, \mymatrix{q}_{hp})_{0,\Omega} &= \ell(\myvec{v}_{hp}) & \qquad & \forall \, (\myvec{v}_{hp},\mymatrix{q}_{hp}) \in V_{hp}\times Q_{hp},\label{eq:discrete_mixed_variationalF_01}\\
 (\mymatrix{\mu}_{hp}-\mymatrix{\lambda}_{hp}, \mymatrix{p}_{hp})_{0,\Omega}
 &\leq 0 & \quad & \forall \, \mymatrix{\mu}_{hp} \in \Lambda_{hp}^{(i)}\label{eq:discrete_mixed_variationalF_02}
\end{alignat}
\end{subequations}
with $i=w$. It turns out that the discrete formulations \eqref{eq:discrete_variq_second_kind} and \eqref{eq:discrete_mixedF_01} are equivalent in the sense of the next result, Theorem~\ref{thm:discreteEquiv_01}, where $\mathcal{P}_{hp} : Q\longrightarrow Q_{hp}$ denotes the $L^2$-projection operator.

% ----------------------------------------------------------------------------------------- 
\begin{theorem}[{[P3, Thm.~2]}]\label{thm:discreteEquiv_01}
If the pair $(\myvec{u}_{hp},\mymatrix{p}_{hp})\in V_{hp}\times Q_{hp}$ solves \eqref{eq:discrete_variq_second_kind}, then the triple $(\myvec{u}_{hp},\mymatrix{p}_{hp},\mymatrix{\lambda}_{hp})$ with 
\begin{align}\label{eq:repres_lagrangeM}
 \mymatrix{\lambda}_{hp} 
  = \mathcal{P}_{hp} \left( \operatorname{dev}(\mymatrix{\sigma}(\myvec{u}_{hp},\mymatrix{p}_{hp}) - \HH \, \mymatrix{p}_{hp}) \right)
\end{align}
is a solution of \eqref{eq:discrete_mixedF_01}. Conversely, if $(\myvec{u}_{hp},\mymatrix{p}_{hp},\mymatrix{\lambda}_{hp})\in V_{hp}\times Q_{hp}\times \Lambda_{hp}^{(w)}$ solves \eqref{eq:discrete_mixedF_01} then $(\myvec{u}_{hp},\mymatrix{p}_{hp})$ solves \eqref{eq:discrete_variq_second_kind} and the identity \eqref{eq:repres_lagrangeM} holds true.
\end{theorem}
% -----------------------------------------------------------------------------------------

\pagebreak
In particular, Theorem~\ref{thm:discreteEquiv_01} guarantees the existence of a discrete solution $(\myvec{u}_{hp},\mymatrix{p}_{hp},\mymatrix{\lambda}_{hp})\in V_{hp}\times Q_{hp}\times\Lambda_{hp}^{(w)}$ of the mixed problem \eqref{eq:discrete_mixedF_01}, the first two components of which are uniquely determined by means of Theorem~\ref{thm:uniqueEx_VIQ}. The remaining uniqueness of the discrete Lagrange multiplier $\mymatrix{\lambda}_{hp}$ is then guaranteed by the \emph{discrete inf-sup condition}
\begin{align}\label{BBS_eq:babusca_brezzi_cond}
 \sup_{\substack{(\myvec{v}_{hp},\mymatrix{q}_{hp})\in V_{hp}\times Q_{hp} \\ \norm{ (\myvec{v}_{hp},\mymatrix{q}_{hp}) } \neq 0}} \frac{(\mymatrix{\mu}_{hp},\mymatrix{q}_{hp})_{0,\Omega}}{\norm{ (\myvec{v}_{hp},\mymatrix{q}_{hp}) }}
  = \norm{ \mymatrix{\mu}_{hp} }_{0,\Omega} \qquad 
  \forall \, \mymatrix{\mu}_{hp}\in Q_{hp},
\end{align}
see [P2, Lem.~4], which particularly holds true for all $\mymatrix{\mu}_{hp}\in\Lambda_{hp}^{(w)}$ as $\Lambda_{hp}^{(w)}\subseteq Q_{hp}$. In the papers [P2, P3], beside \eqref{eq:mixed_variationalF} a second discretization of the mixed variational formulation \eqref{eq:mixed_variationalF} is considered, which only differs in the choice of admissible discrete Lagrange multipliers. For this purpose, we discretize the representation \eqref{eq:representation_Lambda01} by
\begin{align*}
 \Lambda_{hp}^{(s)} 
  := \big\lbrace \mymatrix{\mu}_{hp}\in Q_{hp} \; ; \; \abs{ \mymatrix{\mu}_{hp} \big( \myvec{F}_T (\hat{\myvec{x}}_{k,T})\big) }_F \leq \sigma_y \text{ for all } 1\leq k \leq n_T \text{ and } T \in \mathcal{T}_h \big\rbrace.
\end{align*}
Hence, the corresponding discrete mixed formulation is to find a triple $(\myvec{u}_{hp},\mymatrix{p}_{hp},\mymatrix{\lambda}_{hp})\in V_{hp}\times Q_{hp}\times\Lambda_{hp}^{(s)}$ that satisfies \eqref{eq:discrete_mixedF_01} with $\Lambda_{hp}^{(s)}$. The unique existence of a discrete solution is proved in [P2, Thm.~5], which also contains some \emph{stability estimates}. Under the requirement
\begin{align}\label{eq:requirement_detF}
 \det \nabla \myvec{F}_T \in \PP_1(\wh{T}) \qquad 
 \forall \, T\in\T_h \text{ with } p_T\geq 2
\end{align}
the following relation between the two presented sets of admissible discrete Lagrange multipliers holds true.

% -----------------------------------------------------------------------------------------
\begin{theorem}[{[P2, Thm.~7]}]\label{thm:set_discreteLM}
 On condition that \eqref{eq:requirement_detF} holds true, it follows that $\Lambda_{hp}^{(s)} = \Lambda_{hp}^{(w)}$. 
\end{theorem}
% -----------------------------------------------------------------------------------------
%
While \eqref{eq:requirement_detF} is no restriction for lower-order methods, i.e.~in the case that $p_T = 1$ for all $T\in\T_h$, or in the two-dimensional case, it slightly limits the mesh elements' shape in the case that $d=3$. Note that under the requirement \eqref{eq:requirement_detF}, the two discrete mixed formulations \eqref{eq:discrete_mixedF_01} coincide for $i=s,w$ and, thus, all three discretizations turn out to be equivalent in that case.

% -----------------------------------------------------------------------------------------

\subsubsection*{Representation as a System of Decoupled Nonlinear Equations}

In order to rewrite the discrete mixed problem \eqref{eq:discrete_mixedF_01} with $\Lambda_{hp}^{(w)}$ in terms of a system of nonlinear equations we first decouple the constraints in $\Lambda_{hp}^{(w)}$ and \eqref{eq:discrete_mixed_variationalF_02} with $\Lambda_{hp}^{(w)}$. For this purpose, let $\{\wh{\phi}_{k,T}\}_{k\in\ul{n_T}}$ be the Lagrange basis functions on $\wh{T}$ defined via the Gauss points $\hat{\myvec{x}}_{l,T}$, i.e.
\begin{align*}
 \wh{\phi}_{k,T}\in\PP_{p_T-1}(\wh{T}), \quad
 \wh{\phi}_{k,T}(\hat{\myvec{x}}_{l,T}) = \delta_{kl} \qquad
 \forall \, k,l\in\ul{n_T} \quad \forall \, T\in\T_h,
\end{align*}
where $\delta_{kl}$ is the usual Kronecker delta symbol. Moreover, let $\phi_1,\ldots,\phi_N$ be piece-wise defined as
\begin{align*}
 \phi_{\zeta(k,T')\, | \, T} 
 := \begin{cases}
                \wh{\phi}_{k,T'} \circ \myvec{F}_{T'}^{-1},& \text{if } T=T',\\
                0,& \text{if } T\neq T', \\
    \end{cases} 
 \qquad  T,T'\in\T_h, \quad k\in\ul{n_T},
\end{align*}
where $\zeta:\big\lbrace (k,T) \; ; \;  T\in\mathcal{T}_h, \, k\in\ul{n_T} \big\rbrace\rightarrow \lbrace 1,\ldots,N \rbrace$ with $N :=\sum_{T\in\mathcal{T}_h} n_T$ is a one-to-one numbering. Finally, let $\varphi_1,\ldots,\varphi_N$ be the \emph{biorthogonal basis functions} to $\phi_1,\ldots,\phi_N$ that are uniquely determined by the conditions $\varphi_{\zeta(k,T)\, | \, T} \circ\myvec{F}_T\in\PP_{p_T-1}(\wh{T})$ for $T\in\T_h$, $k\in\ul{n_T}$ and 
\begin{align*}
 (\phi_i,\varphi_j)_{0,\Omega} 
  = \delta_{ij} \, (\phi_i, 1)_{0,\Omega} \qquad 
 \forall\, i,j\in\ul{N}.
\end{align*}
We note that under assumption \eqref{eq:requirement_detF} we have $\varphi_{\zeta(k,T)} = \phi_{\zeta(k,T)}$ for $k\in\ul{n_T}$ and $T\in\T_h$ and from [P3, Sec.~3.3] we obtain
\pagebreak
\begin{align*}
 Q_{hp} 
  = \left\lbrace \sum_{i=1}^N \mymatrix{q}_i \, \phi_i \; ; \; \mymatrix{q}_i\in\SS_{d,0} \right\rbrace
  = \left\lbrace \sum_{i=1}^N \mymatrix{\mu}_i \, \varphi_i \; ; \; \mymatrix{\mu}_i\in\SS_{d,0} \right\rbrace.
\end{align*}
Therewith, and by defining the quantities $D_i := (\phi_i, 1)_{0,\Omega}$ and $\sigma_i := D_i^{-1} (\sigma_y, \phi_i)_{0,\Omega}$ for all $i\in\ul{N}$ the following result is obtained, which allows to decouple the constraints in $\Lambda_{hp}^{(w)}$ and \eqref{eq:discrete_mixed_variationalF_02} with $\Lambda_{hp}^{(w)}$.

% -----------------------------------------------------------------------------------------
\begin{theorem}[{[P1, Thm.~2]}]\label{thm:decoupling}
 It holds
 \begin{align*}
  \Lambda_{hp}^{(w)} 
   = \left\lbrace \sum_{i=1}^N \mymatrix{\mu}_i \, \varphi_i \; ; \;  \mymatrix{\mu}_i \in \SS_{d,0} \text{ and } \abs{ \mymatrix{\mu}_i }_F \leq \sigma_i \right\rbrace.
 \end{align*}
 Furthermore, by representing $\mymatrix{\lambda}_{hp}\in \Lambda_{hp}^{(w)}$ and $\mymatrix{p}_{hp}\in Q_{hp}$ as $\mymatrix{\lambda}_{hp} = \sum_{i\in\ul{N}} \mymatrix{\lambda}_i \, \varphi_i$ and $\mymatrix{p}_{hp} = \sum_{i\in\ul{N}} \mymatrix{p}_i \, \phi_i$, respectively, $\mymatrix{\lambda}_{hp}$ satisfies the inequality \eqref{eq:discrete_mixed_variationalF_02} with $\Lambda_{hp}^{(w)}$ if and only if 
 %$\mymatrix{\lambda}_i : \mymatrix{p}_i = \sigma_i \, \abs{ \mymatrix{p}_i }_F$ for all $i\in\ul{N}$.
\begin{align*}%\label{eq:char_IQconstraint}
 \mymatrix{\lambda}_i : \mymatrix{p}_i 
  = \sigma_i \, \abs{ \mymatrix{p}_i }_F \qquad
 \forall \, i\in\ul{N}.
\end{align*}
\end{theorem}
% -----------------------------------------------------------------------------------------

If $(\myvec{u}_{hp},\mymatrix{p}_{hp},\mymatrix{\lambda}_{hp})\in V_{hp}\times Q_{hp}\times \Lambda_{hp}^{(w)}$ is the discrete solution of \eqref{eq:discrete_mixedF_01} with $\Lambda_{hp}^{(w)}$ then Theorem~\ref{thm:decoupling} together with the Cauchy-Schwarz inequality yield the two implications
\begin{align*}
 \abs{ \mymatrix{\lambda}_i }_F < \sigma_i 
  \;\Longrightarrow \; \mymatrix{p}_i = \mymatrix{0} , \qquad
 \abs{ \mymatrix{\lambda}_i }_F = \sigma_i 
  \; \Longrightarrow \; \exists\, c \geq 0 \text{ with } \mymatrix{p}_i = c \, \mymatrix{\lambda}_i, 
\end{align*}
see [P1, Sec.~2], which suggest to introduce the nonlinear functions $\mymatrix{\chi}_i : \SS_{d,0}\longrightarrow\SS_{d,0}$, for some $\rho > 0$, as
\begin{align*}
 \mymatrix{\chi}_i(\mymatrix{p}_i, \mymatrix{\lambda}_i) 
  := \max \big\lbrace \sigma_i, \abs{ \mymatrix{\lambda}_i + \rho \, \mymatrix{p}_i }_F \big\rbrace \, \mymatrix{\lambda}_i - \sigma_i \, (\mymatrix{\lambda}_i + \rho \, \mymatrix{p}_i) \qquad
 \forall \, i\in\ul{N}. 
\end{align*} 
The next result allows to rewrite \eqref{eq:discrete_mixedF_01} with $\Lambda_{hp}^{(w)}$ as a nonlinear system of equations.

% -----------------------------------------------------------------------------------------
\begin{theorem}[{[P1, Thm.~3]}]
 The discrete Lagrange multiplier $\mymatrix{\lambda}_{hp}\in\Lambda_{hp}^{(w)}$ satisfies the condition \eqref{eq:discrete_mixed_variationalF_02} with $\Lambda_{hp}^{(w)}$ if and only if
 \begin{align*}
  \mymatrix{\chi}_i(\mymatrix{p}_i, \mymatrix{\lambda}_i) 
   = \mymatrix{0} \qquad
  \forall \, i\in\ul{N}.
 \end{align*}
\end{theorem}
% -----------------------------------------------------------------------------------------

In order to introduce the nonlinear system of equations, an appropriate basis of $\SS_{d,0}$ has to be chosen first. If $d = 2$ we take $\mymatrix{\Phi}_1 := \frac{1}{\sqrt{2}} \, \big( \begin{smallmatrix} 1 & 0 \\ 0 & -1 \end{smallmatrix} \big)$ and $\mymatrix{\Phi}_2 := \frac{1}{\sqrt{2}} \, \big( \begin{smallmatrix} 0 & 1 \\ 1 & 0 \end{smallmatrix} \big)$. For $d=3$ let
\begin{align*}
 \mymatrix{\Phi}_1 
  := \tfrac{1}{\sqrt{2}} & \, \bigg( \begin{smallmatrix} 1 & 0 & 0 \\ 0 & -1 & 0 \\ 0 & 0 & 0 \end{smallmatrix} \bigg), \quad
 \mymatrix{\Phi}_2 
  := \tfrac{1}{\sqrt{6}} \, \bigg( \begin{smallmatrix} 1 & 0 & 0 \\ 0 & 1 & 0 \\ 0 & 0 & -2 \end{smallmatrix} \bigg), \quad
 \mymatrix{\Phi}_3 
  := \tfrac{1}{\sqrt{2}} \, \bigg( \begin{smallmatrix} 0 & 1 & 0 \\ 1 & 0 & 0 \\ 0 & 0 & 0 \end{smallmatrix} \bigg), \\
 & \quad \mymatrix{\Phi}_4 
  := \tfrac{1}{\sqrt{2}} \, \bigg( \begin{smallmatrix} 0 & 0 & 1 \\ 0 & 0 & 0 \\ 1 & 0 & 0 \end{smallmatrix} \bigg), \quad
 \mymatrix{\Phi}_5 
  := \tfrac{1}{\sqrt{2}} \, \bigg( \begin{smallmatrix} 0 & 0 & 0 \\ 0 & 0 & 1 \\ 0 & 1 & 0 \end{smallmatrix} \bigg). 
\end{align*}
Note that these basis functions are orthonormal with respect to the Frobenius inner product. Therewith, we can represent $\mymatrix{p}_{hp}$ and $\mymatrix{\lambda}_{hp}$ in terms of
\begin{align*}
 \mymatrix{p}_{hp}
  = \sum_{i=1}^N \sum_{k=1}^L p_{L(i-1)+k} \, \mymatrix{\Phi}_k \, \phi_i, \qquad
 \mymatrix{\lambda}_{hp}
  = \sum_{i=1}^N \sum_{k=1}^L \lambda_{L(i-1)+k} \, \mymatrix{\Phi}_k \, \varphi_i
\end{align*}
with $L := \frac{1}{2} \, (d-1) \, (d+2)$. Furthermore, by choosing functions $\vartheta_1,\ldots,\vartheta_M$ such that $\lbrace \myvec{e}_k \, \vartheta_i \; ; \; k\in\ul{d} \text{ and } i\in\ul{M} \rbrace$ forms a basis of $V_{hp}$, where $\myvec{e}_k$ denotes the $k$-th Euclidean unit vector in $\RR^d$, we have
\begin{align*}
 \myvec{u}_{hp} = \sum_{i=1}^M \sum_{k=1}^d u_{d(i-1)+k} \, \myvec{e}_k \, \vartheta_i.
\end{align*}
Thus, the discrete solution $(\myvec{u}_{hp},\mymatrix{p}_{hp},\mymatrix{\lambda}_{hp})\in V_{hp}\times Q_{hp}\times \Lambda_{hp}^{(w)}$ of \eqref{eq:discrete_mixedF_01} with $\Lambda_{hp}^{(w)}$ can completely be represented by the coefficient vectors $\myvec{v}_{\myvec{u}} := (u_1,\ldots,u_{dM})^{\top}\in\RR^{dM}$, $\myvec{v}_{\mymatrix{p}} := (p_1,\ldots,p_{LN})^{\top}\in\RR^{LN}$ and $\myvec{v}_{\mymatrix{\lambda}} := (\lambda_1,\ldots,\lambda_{LN})^{\top}\in\RR^{LN}$. Finally, writing $\mymatrix{\chi}_i(\mymatrix{p}_i, \mymatrix{\lambda}_i) = \sum_{k=1}^L \chi_{i,k} \, \mymatrix{\Phi}_k$ and defining the vectors $\chi_i := (\chi_{i,1},\ldots,\chi_{i,L})^{\top}\in\RR^{L}$ turns the discrete mixed problem \eqref{eq:discrete_mixedF_01} with $\Lambda_{hp}^{(w)}$ into a system of decoupled nonlinear equations, given by
\begin{align}\label{eq:semismooth_newtonF}
\myvec{F} (\myvec{v}_{\myvec{u}}, \myvec{v}_{\mymatrix{p}}, \myvec{v}_{\mymatrix{\lambda}} ) 
 := \begin{pmatrix}
\mymatrix{K} \, \myvec{v}_{\myvec{u}} - \mymatrix{B} \, \myvec{v}_{\mymatrix{p}} - \myvec{l} \\
-\mymatrix{B}^{\top} \, \myvec{v}_{\myvec{u}} + \mymatrix{C} \, \myvec{v}_{\mymatrix{p}} + \mymatrix{D} \, \myvec{v}_{\mymatrix{\lambda}} \\
{\chi}_1\left( \sum_{k=1}^L p_k \, \mymatrix{\Phi}_k, \sum_{k=1}^L \lambda_k \, \mymatrix{\Phi}_k \right) \\
\vdots \\
{\chi}_N\left( \sum_{k=1}^L p_{L(N-1)+k} \, \mymatrix{\Phi}_k, \sum_{k=1}^L \lambda_{L(N-1)+k} \, \mymatrix{\Phi}_k \right)
\end{pmatrix} = \myvec{0},
\end{align}
where the symmetric, positive definite matrices $\mymatrix{K}\in\RR^{dM\times dM}$ and $\mymatrix{C}\in\RR^{LN\times LN}$, the positive definite diagonal matrix $\mymatrix{D}\in\RR^{LN\times LN}$ and the coupling matrix $\mymatrix{B}\in\RR^{dM\times LN}$ are given component-wise by
\begin{alignat*}{2}
 K_{d(i-1)+k,d(j-1)+l} 
  &= a\big( (\mathfrak{e}_l \, \vartheta_j,\mathbf{0}) , (\mathfrak{e}_k \, \vartheta_i,\mathbf{0}) \big)
   \qquad &&\forall \, i,j\in\ul{M} \quad \forall \, k,l\in\ul{d}, \\
 C_{L(i-1)+k,L(j-1)+l} 
  &= a\big( (\mathfrak{0},\boldsymbol{\Phi}_l \, \phi_i) , (\mathfrak{0},\boldsymbol{\Phi}_k \, \phi_j) \big)
   \qquad &&\forall \, i,j\in\ul{N} \quad \forall \, k,l\in\ul{L}, \\
 D_{L(i-1)+k,L(j-1)+l} 
  &= \delta_{lk} \, \delta_{ij} \, D_i
   \qquad &&\forall \, i,j\in\ul{N} \quad \forall \, k,l\in\ul{L}, \\
 B_{d(i-1)+k,L(j-1)+l} 
  &= a\big( (\myvec{o},\mymatrix{\Phi}_l \, \phi_j) , (\myvec{e}_k \, \vartheta_i, \mathbf{0}) \big)
   \qquad &&\forall \, i\in\ul{N} \quad \forall \, j\in\ul{M} \quad \forall \, k\in\ul{d} \quad \forall \, l \in \ul{L},
\end{alignat*}
and the vector $\myvec{l}\in\RR^{dM}$ has the components $l_{d(i-1)+k} = \ell(\myvec{e}_k \, \vartheta_i)$ for $i\in\ul{M}$ and $k\in\ul{d}$. Note that the first equation lines in \eqref{eq:semismooth_newtonF} are linear and those according to $\chi_1,\ldots,\chi_N$ are the nonlinear but semi-smooth ones. Hence, to solve \eqref{eq:semismooth_newtonF} by an iterative solver one may apply the \emph{semi-smooth Newton solver}, proposed in [P1, Sec.~3], which is given by
\begin{subequations} %\label{BBS_eq:semismooth_Newton}
\begin{align*}
 \mymatrix{H}_k \Big( \Delta\myvec{v}_{\myvec{u}}^{(k)}, \Delta\myvec{v}_{\mymatrix{p}}^{(k)},  \Delta\myvec{v}_{\mymatrix{\lambda}}^{(k)} \Big)^{\top} 
  &= - \myvec{F} \Big( \myvec{v}_{\myvec{u}}^{(k)}, \myvec{v}_{\mymatrix{p}}^{(k)}, \myvec{v}_{\mymatrix{\lambda}}^{(k)} \Big) \\
\Big( \myvec{v}_{\myvec{u}}^{(k+1)}, \myvec{v}_{\mymatrix{p}}^{(k+1)}, \myvec{v}_{\mymatrix{\lambda}}^{(k+1)} \Big)^{\top} 
 &= \Big( \myvec{v}_{\myvec{u}}^{(k)}, \myvec{v}_{\mymatrix{p}}^{(k)}, \myvec{v}_{\mymatrix{\lambda}}^{(k)} \Big)^{\top} + t_k \Big( \Delta\myvec{v}_{\myvec{u}}^{(k)}, \Delta\myvec{v}_{\mymatrix{p}}^{(k)}, \Delta\myvec{v}_{\mymatrix{\lambda}}^{(k)} \Big)^{\top}
\end{align*}
\end{subequations}
with $\mymatrix{H}_k\in\partial\myvec{F}\big( \myvec{v}_{\myvec{u}}^{(k)}, \myvec{v}_{\mymatrix{p}}^{(k)}, \myvec{v}_{\mymatrix{\lambda}}^{(k)} \big)$, where $\partial\myvec{F}(\cdot)$ denotes the Clarke subdifferential of $\myvec{F}(\cdot)$. Thereby, the step length parameter $t_k \in (0,1]$ has to be chosen by an adequate step length selection procedure. The numerical examples in [P1, Sec.~5] show the applicability of the semi-smooth Newton solver as well as its robustness to $h$, $p$ and the projection parameters. Furthermore, we observe superlinear convergence properties, which will be investigated in a forthcoming work. 

% -----------------------------------------------------------------------------------------

\subsubsection*{A Priori Error Analysis}

The convergence analysis in [P2, Sec.~5] is derived for the discrete mixed problem \eqref{eq:discrete_mixedF_01} with $\Lambda_{hp}^{(s)}$ where, in addition, the requirement \eqref{eq:requirement_detF} is assumed to hold true, and is based on the following \emph{a priori error estimate}.

% -----------------------------------------------------------------------------------------
\begin{theorem}[{[P2, Thm.~6]}]\label{thm:apriori_errorEst}
 There exist two positive constants $c_1$ and $c_2$ such that for all $(\myvec{v}_{hp},\mymatrix{q}_{hp},\mymatrix{\mu}_{hp})\in V_{hp}\times Q_{hp}\times \Lambda_{hp}$ and any $\mymatrix{\mu}\in\Lambda$ the following a priori estimate holds true
 \begin{align*}
   \norm{ (\myvec{u}-\myvec{u}_{hp}, \mymatrix{p}-\mymatrix{p}_{hp}) }^2 + \norm{ \mymatrix{\lambda}-\mymatrix{\lambda}_{hp} }_{0,\Omega}^2
   &\leq  c_1 \, \Big( \norm{ (\myvec{u}-\myvec{v}_{hp}, \mymatrix{p}-\mymatrix{q}_{hp}) }^2 +  \norm{ \mymatrix{\lambda}-\mymatrix{\mu}_{hp} }_{0,\Omega}^2 \Big)  \\
   &\qquad + c_2 \, (\mymatrix{p}, \mymatrix{\lambda}_{hp}-\mymatrix{\mu} + \mymatrix{\lambda}-\mymatrix{\mu}_{hp})_{0,\Omega}.
 \end{align*}
\end{theorem}
% -----------------------------------------------------------------------------------------

\pagebreak
By exploiting the equivalent representations for the set of admissible discrete Lagrange multipliers, cf.~Theorem~\ref{thm:set_discreteLM}, convergence is obtained under minimal regularity conditions.

% -----------------------------------------------------------------------------------------
\begin{theorem}[{[P2, Thm.~8]}]
 The following norm convergence holds true
 \begin{align*}
  \lim_{h/p \to 0} \Big( \norm{ \myvec{u}-\myvec{u}_{hp} }_{1,\Omega}^2 + \norm{ \mymatrix{p}-\mymatrix{p}_{hp} }_{0,\Omega}^2 + \norm{ \mymatrix{\lambda}-\mymatrix{\lambda}_{hp} }_{0,\Omega}^2 \Big) 
   = 0.
 \end{align*}
\end{theorem}
% -----------------------------------------------------------------------------------------

For the derivation of convergence rates one has to assume a certain regularity of the weak solution $(\myvec{u},\mymatrix{p},\mymatrix{\lambda})\in \VV\times Q\times \Lambda$. For this purpose, let $(\myvec{u},\mymatrix{p},\mymatrix{\lambda})\in H^s(\Omega,\RR^d)\times H^t(\Omega,\RR^{d\times d})\times H^l(\Omega,\RR^{d\times d})$ for some $s\geq 1$ and $t,l\geq 0$. As in the lowest-order case the discretization is particularly conforming in the discrete Lagrange multiplier, i.e.~$\Lambda_{hp}^{(s)}\subseteq \Lambda$, we obtain optimal convergence rates. Thereby, the assumption \eqref{eq:requirement_detF} can be neglected so that there is no restriction to certain elements' shapes in that case. In the following, the symbol $\lesssim$ is used to hide a constant $c>0$ in the expression $A \leq c \, B$ which is independent of $h$ and $p$.

% -----------------------------------------------------------------------------------------
\begin{theorem}[{[P2, Thm.~9]}]
If $p_T = 1$ for all $T \in \mathcal{T}_h$ we obtain the optimal order of convergence
\begin{align*}
   \norm{ (\myvec{u}-\myvec{u}_{hp}, \mymatrix{p}-\mymatrix{p}_{hp}) }^2  + \norm{ \mymatrix{\lambda}-\mymatrix{\lambda}_{hp} }_{0,\Omega}^2 
   & \lesssim h^{2 \, \min(1,s-1,t,l)} \, \big( |\myvec{u}|_{s,\Omega}^2 + |\mymatrix{p}|_{t,\Omega}^2 + |\mymatrix{\lambda}|_{l,\Omega}^2 \big),
\end{align*}
where $\abs{\cdot}_k$ denotes the Sobolev-seminorm on $H^k(\Omega,X)$ for $X\in\lbrace \RR^d, \RR^{d\times d}\rbrace$ and $k\geq 0$.
\end{theorem}
% -----------------------------------------------------------------------------------------

When applying a higher-order method the non-conformity error $(\mymatrix{p}, \mymatrix{\lambda}_{hp} - \mymatrix{\mu})_{0,\Omega}$ arising in the a priori error estimate of Theorem~\ref{thm:apriori_errorEst} has to be estimated for any $\mymatrix{\mu}\in\Lambda$. As in that case we have to use a nodal interpolation operator we may not achieve optimal convergence rates.

% -----------------------------------------------------------------------------------------
\begin{theorem}[{[P2, Thm.~11]}]
 For $t,l > d/2$ it holds
\begin{align*}
 \norm{ (\myvec{u}-\myvec{u}_{hp}, \mymatrix{p}-\mymatrix{p}_{hp}) }^2  + \norm{ \mymatrix{\lambda}-\mymatrix{\lambda}_{hp} }_{0,\Omega}^2 
 \lesssim \frac{h^{\min(p,2s-2,t,l)}}{p^{\min(2s-2,t,l)}}.
\end{align*}
\end{theorem}

Thus, the convergence rates are potentially suboptimal by a factor of Two compared to the best possible rates for the finite element discretization.

% -----------------------------------------------------------------------------------------

\subsubsection*{A Posteriori Error Estimates}

The following estimates are employable to any triple $(\myvec{u}_N,\mymatrix{p}_N,\mymatrix{\lambda}_N)\in \VV\times Q \times Q$, in particular, to the discrete solution of the variational inequality \eqref{eq:discrete_variq_second_kind} as well as to the discrete solutions of the mixed problem \eqref{eq:discrete_mixedF_01} with $\Lambda_{hp}^{(w)}$ or $\Lambda_{hp}^{(s)}$. The derivation of upper and lower bounds is thereby based on the following \emph{auxiliary problem}: Find a pair $(\myvec{u}^*,\mymatrix{p}^*)\in \VV\times Q$ such that
\begin{align*}%\label{eq:additional_problem}
 a\big( (\myvec{u}^*,\mymatrix{p}^*),(\myvec{v},\mymatrix{q}) \big) 
 = \ell(\myvec{v}) - (\mymatrix{\lambda}_N, \mymatrix{q})_{0,\Omega} \qquad 
 \forall \, (\myvec{v},\mymatrix{q})\in \VV\times Q.
\end{align*}
Moreover, for $\mymatrix{\mu}\in\Lambda$, let the \emph{global plasticity error contribution} be denoted by
\begin{align*}
 \sfE(\mymatrix{\mu}) 
 := \Vert \mymatrix{\mu}-\mymatrix{\lambda}_N\Vert_{0,\Omega}^2 + \psi(\mymatrix{p}_N) - (\mymatrix{\mu}, \mymatrix{p}_N)_{0,\Omega}.
\end{align*}
%

% -----------------------------------------------------------------------------------------
\begin{theorem}[{[P3, Thm.~5]}]\label{thm:upperB}
 For every $\mymatrix{\mu}\in\Lambda$ the following upper bound holds true
 \begin{align*} %\label{eq:reliability_estimate}
  \norm{ (\myvec{u}-\myvec{u}_N, \mymatrix{p}-\mymatrix{p}_N) }^2 + \norm{ \mymatrix{\lambda}-\mymatrix{\lambda}_N }_{0,\Omega}^2 
   \lesssim \norm{ (\myvec{u}^*-\myvec{u}_N, \mymatrix{p}^*-\mymatrix{p}_N) }^2 + \sfE(\mymatrix{\mu}).
 \end{align*}
\end{theorem}
% -----------------------------------------------------------------------------------------

In order to obtain lower error estimates, the error contribution $\sfE(\cdot)$ can explicitly be minimized over $\Lambda$. By [P2, Lem.~15], the unique minimizer $\mymatrix{\mu}^{*}\in\Lambda$ is given by
\begin{align*}%\label{eq:defi_minimizer_mu}
 \mymatrix{\mu}^* 
 := \min\left\lbrace 1, \frac{\sigma_y} {\abs{ \widehat{\mymatrix{\mu}} }_F^{-1}} \right\rbrace \widehat{\mymatrix{\mu}}, \qquad 
 \widehat{\mymatrix{\mu}}:= \mymatrix{\lambda}_{N} + \frac{1}{2} \, \mymatrix{p}_{N}.
\end{align*}
Therewith, the following lower bounds can be shown.

% -----------------------------------------------------------------------------------------
\begin{theorem}[{[P3, Thm.~8]}]\label{thm:lowerB}
 For the minimizer $\mymatrix{\mu}^{*}\in\Lambda$ the following lower bound holds true
\begin{align*} %\label{eq:efficiency_est}
 \norm{ (\myvec{u}^*-\myvec{u}_N, \mymatrix{p}^*-\mymatrix{p}_N) }^2 + \sfE(\mymatrix{\mu}^*)
  \lesssim \norm{ (\myvec{u}-\myvec{u}_N, \mymatrix{p}-\mymatrix{p}_N) }^2 + \norm{ \mymatrix{\lambda}-\mymatrix{\lambda}_N }_{0,\Omega}^2 +  \norm{ \mymatrix{p}-\mymatrix{p}_N }_{0,\Omega}.
\end{align*}
\end{theorem}
% -----------------------------------------------------------------------------------------

To derive a \emph{residual-based a posteriori error estimator}, let $V_N := V_{hp}$, $Q_N := Q_{hp}$, take $(\myvec{u}_N,\mymatrix{p}_N)\in V_N\times Q_N$ and choose $\mymatrix{\lambda}_N$ such that either
\begin{align*} %\label{eq:find_lambda}
  ( \mymatrix{\lambda}_N, \mymatrix{q}_N )_{0,\Omega}
 = \ell(\myvec{v}_N) - a\big( ( \myvec{u}_N, \mymatrix{p}_N ) , (\myvec{v}_N,\mymatrix{q}_N) \big) \qquad  
 \forall \, (\myvec{v}_N,\mymatrix{q}_N) \in V_N\times Q_N
\end{align*}
holds true or take $\mymatrix{\lambda}_N := \operatorname{dev}\big( \mymatrix{\sigma}(\myvec{u}_N,\mymatrix{p}_N) - \HH \, \mymatrix{p}_N \big)$. If one solves a discretization of the mixed variational formulation \eqref{eq:mixed_variationalF}, clearly, the thereby obtained discrete Lagrange multiplier $\mymatrix{\lambda}_N\in Q_N$ can be used as well. For $T\in\T_h$, define the local error quantities 
\begin{align*} %\label{eq:eta_T}
 \eta_T^2 
 := \frac{h_T^2}{p_T^2} \, \norm{ \myvec{f}_N + \operatorname{div}\mymatrix{\sigma}( \myvec{u}_N,\mymatrix{p}_N ) }_{0,T}^2 
 + \sum_{e \in \E_T^{I}} \frac{h_e}{2p_e} \, \norm{ \llbracket \mymatrix{\sigma}( \myvec{u}_N,\mymatrix{p}_N ) \, \myvec{n}_e \rrbracket }_{0,e}^2 
 + \sum_{e \in \E_T^{N}}  \frac{h_e}{p_e} \, \norm{ \mymatrix{\sigma} ( \myvec{u}_N,\mymatrix{p}_N ) \, \myvec{n}_e - \myvec{g}_N }_{0,e}^2,
\end{align*}
where $\E_T^{I}$ and $\E_T^{N}$ are the edges of $T$ that lie in the interior of $\Omega$ and on the Neumann-boundary $\Gamma_N$, respectively, $h_e$ and $p_e$ are the local edge size and polynomial degree and $\llbracket\cdot\rrbracket$ represents the usual jump function, cf.~[P3, Sec.~4]. Moreover, for $\mymatrix{\mu}\in\Lambda$, let
\begin{align*}
 \eta^2_T(\mymatrix{\mu}) 
  := \eta_T^2 + \norm{ \operatorname{dev} \big( \mymatrix{\sigma} (\myvec{u}_N,\mymatrix{p}_N ) - \HH \, \mymatrix{p}_N \big) - \mymatrix{\lambda}_N }_{T,\Omega}^2 
 + \norm{ \mymatrix{\mu} - \mymatrix{\lambda}_N }_{0,T}^2 + (\sigma_y, \vert \mymatrix{p}_N\vert_F )_{0,T} - (\mymatrix{\mu} , \mymatrix{p}_N)_{0,T}
\end{align*}
and set $\eta^2(\mymatrix{\mu}) := \sum_{T\in\T_h} \eta^2_T(\mymatrix{\mu})$. Then, for $d=2$ the \emph{a posteriori error estimator} $\eta^2(\mymatrix{\mu})$ satisfy the following reliability estimate, where $\osc^2$ denotes the typical data oscillation terms, see [P3, Sec.~4] for the definition.

% -----------------------------------------------------------------------------------------
\begin{theorem}[{[P3, Thm.~11]}] %\label{thm:residualErrorEst}
For any $\mymatrix{\mu}\in\Lambda$ the following reliability estimate holds true
\begin{align} \label{eq:residual_est_reliable}
 \Vert ( \myvec{u}-\myvec{u}_N, \mymatrix{p}-\mymatrix{p}_N ) \Vert^2 + \Vert \mymatrix{\lambda}-\mymatrix{\lambda}_N \Vert_{0,\Omega}^2  
 &\lesssim  \eta^2(\mymatrix{\mu}) + \osc^2.
\end{align}
\end{theorem}
% -----------------------------------------------------------------------------------------

To obtain an efficiency estimate of $\eta^2(\mymatrix{\mu})$ let the elasticity tensor $\CC$ be constant and assume that the element mapping $\myvec{F}_T(\cdot)$ is affine for any $T\in\T_h$. Unfortunately, we cannot avoid the suboptimality for the estimate, which is, however, common in the framework of elastoplastic problems. \pagebreak

% -----------------------------------------------------------------------------------------
\begin{theorem}[{[P3, Thm.~14]}]\label{thm:globalLowerEst}
For the minimizer $\mymatrix{\mu}^{*}\in\Lambda$ the following (suboptimal) efficiency estimate holds true
\begin{align*}
 \eta^2(\mymatrix{\mu}^*)
 \lesssim_p \Vert (\myvec{u}-\myvec{u}_N, \mymatrix{p}-\mymatrix{p}_N)\Vert^2 + \Vert \mymatrix{\lambda}-\mymatrix{\lambda}_N\Vert_{0,\Omega}^2 + \Vert \mymatrix{p}-\mymatrix{p}_N\Vert_{0,\Omega} + \osc^2.
\end{align*}
\end{theorem}

By using the same arguments leading to Theorem~\ref{thm:upperB} and Theorem~\ref{thm:lowerB} on an element level yields a local version of Theorem~\ref{thm:globalLowerEst}, see [P3, Sec.~4].

% -----------------------------------------------------------------------------------------
%   Part II : $hp$-Adaptivity based on Local Error Reductions
% -----------------------------------------------------------------------------------------

\section{Part II :\, $hp$-Adaptivity Based on Local Error Reductions}

As it was already pointed out in Section~\ref{sec:FEM}, it may be necessary to apply a combination of local $h$- and $p$-refinements in order to recover the optimal algebraic convergence rates or to obtain exponential convergence (even for weak solutions of a low regularity). The key advantage of $hp$-adaptive strategies is their ability to approximate singularities of the weak solution of a boundary value problem. Typically, a posteriori error estimators are used to steer the adaptive refinements as it is done, for instance, in the numerical experiments of the papers [P2, P3]. While a posteriori error estimators can excellently be used to steer $h$-adaptivity, the derivation of effective computable error bounds in the $hp$-adaptive framework is generally known as challenging due to considerable technical difficulties, see e.g.~\cite{ref:Dolejsi_2016, ref:Melenk_2001}. In the case of a pure $h$-adaptive scheme physical elements only have to be flagged for refinement whereas in the context of $hp$-adaptivity one has to choose carefully between various possible $hp$-refinements on the flagged elements. This can be done, for instance, by employing suitable smoothness testing strategies, see \cite{ref:Eibner_2007, ref:Fankhauser_2014, ref:houston_2005, ref:Wihler_2011}. In contrast to this apporaches, the $hp$-adaptive procedure proposed in the paper [P4] does \emph{neither} rely on classical a posteriori error estimation \emph{nor} on smoothness indicators. Instead, it is based on a \emph{prediction strategy} for the reduction of the (global) \emph{energy error} which results from local \emph{$p$-enrichments} or \emph{$hp$-refinements}. The strategy is therefore closely related to the energy minimization technique, presented in \cite{ref:houston_2016}. The energy error represents the discretization error measured in the norm, which is induced by the involved symmetric bilinear form. 
\vspace{0.1cm}

In order to compare different $hp$-refinements on the individual elements of the current mesh, for each element $Q$, the current (global) $hp$-finite element solution $u_{hp}$ is decomposed into two parts: a locally supported part $u_{hp}^{\text{loc}}$ (with support $Q$ or a patch around $Q$) and a globally supported part $\wt{u}_{hp} := u_{hp} - u_{hp}^{\text{loc}}$. The crucial idea of the locally predicting strategy is to compare the current solution $u_{hp}$ with a locally modified solution, in which the local part $u_{hp}^{\text{loc}}$ is replaced by a linear combination of (locally supported) so-called \emph{enrichment functions}, the definition of which relies either on an increase of the local polynomial degree on $Q$ ($p$-enrichment) or on an $hp$-refinement of $Q$ (i.e.~decomposing $Q$ in a few subelements and distributing appropriate polynomial degrees on these subelements). Thereby, the span of the enrichment functions together with the unchanged part $\wt{u}_{hp}$ defines a (low-dimensional) so-called \emph{enrichment} of \emph{replacement space} $\YY$, in which the locally modified solution $u_\YY$ is sought. As the globally supported part $\wt{u}_{hp}$ is explicitly included in the definition of $\YY$ the solution $u_\YY$ can be expressed in terms of a linear combination of the enrichment functions and $\wt{u}_{hp}$, and, thus, represents a \emph{global} solution as well. Therefore, it is reasonable to compare the low-dimensional solution $u_{\YY}$ with $u_{hp}$. It turns out that the computation of the locally \emph{predicted error reduction}, which is defined as the difference between the discretization errors of $u_{hp}$ and the locally obtained modified solution $u_\YY$, involves low-dimensional linear problems that are computationally inexpensive and highly parallelizable. As the locally predicted (but globally effective) error reduction can be therefore computed at a negligible cost, different $p$-enrichments and $hp$-refinements can be compared on the element $Q$ to determine an optimal one, leading to the highest possible predicted global error reduction. The resulting $hp$-adaptive algorithm in [P4] passes through all elements of a current mesh (which can be done in parallel), and compares the predicted error reductions for different $p$-enrichments and $hp$-refinements in order to find an optimal one on each element. Then, by using an appropriate marking strategy, all those elements, from which the most substantial (global) error reduction can be expected, are enriched. Schematically, the proposed algorithm follows the structure:

% -----------------------------------------------------------------------------------------
\begin{center}	
\includegraphics[scale=0.4]{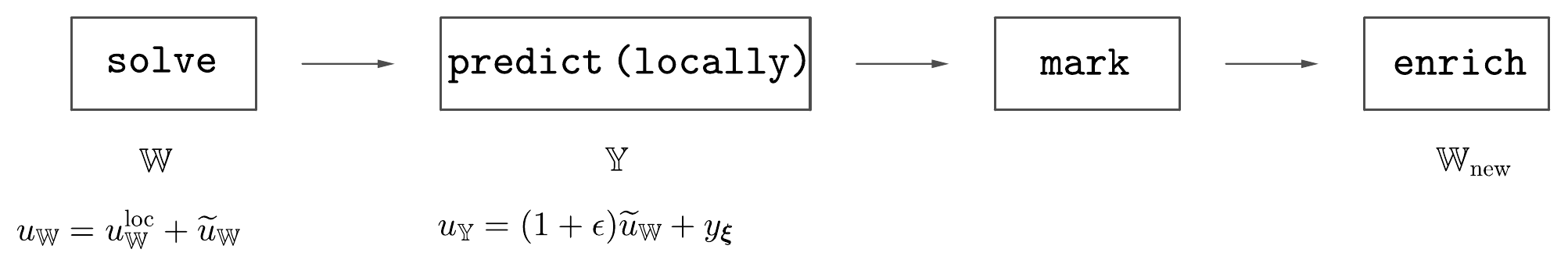}
\end{center}
% -----------------------------------------------------------------------------------------

The following paragraphs present the core ideas of the paper [P4].

% -----------------------------------------------------------------------------------------

\subsubsection*{Locally Predicted Error Reduction in an Abstract Framework}

On a (real) Hilbert space $\XX$ consider the weak formulation: Find a $u_\XX\in\XX$ such that the variational equation
\begin{align*}
 a(u_\XX, v) 
  = b(v) \qquad
 \forall \, v\in\XX
\end{align*}
holds true, where $a : \XX\times \XX\longrightarrow\RR$ is a bounded, symmetric and $\XX$-elliptic bilinear form inducing the norm $v \longmapsto \norm{v}_\XX := a(v,v)^{1/2}$ and $b : \XX\longrightarrow\RR$ is a bounded linear form. Furthermore, consider its Riesz-Galerkin discretization related to a finite-dimensional subspace $\WW\subseteq\XX$, i.e.: Find a $u_\WW\in\WW$ such that
\begin{align}\label{eq:lowDim_problem}
 a(u_\WW, w) 
  = b(w) \qquad
 \forall \, w\in\WW.
\end{align}
By the finite dimension of $\WW$ it is spanned by finitely many basis functions $\phi_1,\ldots,\phi_N$ and for some index set $\mathcal{I}^{\text{loc}}\subset\ul{N}$ let $\WW^{\text{loc}} := \operatorname{span} \lbrace \phi_i \; ; \; i\in\mathcal{I}^{\text{loc}} \rbrace \subset \WW$. By introducing the linear projection operator
\begin{align*}
 \Pi^{\text{loc}} : \WW \longrightarrow \WW^{\text{loc}}, \qquad
 v = \sum_{i\in\ul{N}} v_i \, \phi_i \longmapsto \Pi^{\text{loc}} v := \sum_{i\in\mathcal{I}^{\text{loc}}} v_i \, \phi_i,
\end{align*}
one can decompose the solution $u_\WW$ of \eqref{eq:lowDim_problem} in terms of $u_\WW = \uloc + \ut$, where
\begin{align*}
 \uloc
  := \Pi^{\text{loc}} u_\WW \in\WW^{\text{loc}},\qquad
 \ut
  := u_\WW - \Pi^{\text{loc}} u_{\WW}.
\end{align*}
The idea to improve the Galerkin-approximation $u_\WW$ is to enrich or replace the (local) space $\WW^{\text{loc}}$. For this purpose, choose a (small) set $\mymatrix{\xi} := \lbrace \xi_1,\ldots,\xi_L\rbrace$ of linearly independent elements in $\XX$ -- so-called \emph{enrichment functions} -- such that 
\begin{align*}
 \YY 
  := \operatorname{span} \lbrace \ut, \xi_1,\ldots,\xi_L \rbrace
  \subseteq\XX
\end{align*}
has dimension $L+1$, which implicitly implies that $\ut\neq 0$. If $\WW^{\text{loc}} \subset \YY$ we call $\YY$ a \emph{local enrichment space} and otherwise a \emph{local replacement space}. Consider now the low-dimensional problem: Find a $u_\YY\in\YY$ such that
\begin{align*}
 a(u_{\YY}, v) 
  = b(v) \qquad 
 \forall \, v\in \YY.
\end{align*}
Since $\dim\YY = L+1$ the solution $u_\YY$ can be expressed in terms of a linear combination $u_\YY = (1+\epsilon) \, \ut + y_{\mymatrix{\xi}}$, where $\epsilon\in\RR$ and $y_{\mymatrix{\xi}}\in\operatorname{span} \mymatrix{\xi}$. By introducing the discretization errors $e_\WW := u_\XX - u_\WW$ and $e_\YY := u_\XX - u_\YY$ we define the \emph{predicted error reduction} $\Delta e_{\WW,\YY}$ to be
\begin{align*}
 \Delta e_{\WW,\YY}^2 
  := \norm{e_{\WW}}_\XX^2 - \norm{e_{\YY}}_\XX^2.
\end{align*}

% -----------------------------------------------------------------------------------------
\begin{theorem}[{[P4, Prop.~1]}]
 For the predicted error reduction the following identities hold true
\begin{align*}% \label{eq:errid}
 \Delta e_{\WW,\YY}^2
  = \norm{ u_{\YY} - u_{\WW} }_\XX^2 - 2 \, \res[\YY]{\uloc}
  = \res{y_{\mymatrix{\xi}}} - \res[\YY]{\uloc}
\end{align*}
with the residuals $\rho_\WW$ and $\rho_\YY$ defined as $\res{\cdot} := b(\cdot) - a(u_\WW,\cdot)$ and $\res[\YY]{\cdot} := b(\cdot) - a(u_{\YY},\cdot)$, respectively. 
\end{theorem}
% -----------------------------------------------------------------------------------------

In the case that $\YY$ represents a local enrichment space it follows that $\res[\YY]{\uloc} = 0$ and, thus,
\begin{align*}
 \Delta e_{\WW,\YY}^2
  = \norm{ u_{\YY}-u_{\WW} }_\XX^2 = \res{y_{\myvec{\xi}}},
\end{align*}
which leads to an \emph{energy reduction} property, see [P4, Rem.~1]. For the computation of the predicted error reduction by means of linear algebra let the matrix $\mymatrix{A} = (A_{ij})\in\RR^{L\times L}$ and the vectors $\myvec{b}=(b_i),\myvec{c}=(c_i)\in\RR^L$ be given component-wise by
\begin{align*}
 A_{ij} := a(\xi_j, \xi_i), \qquad
    b_i := b(\xi_i), \qquad
    c_i := a(\ut, \xi_i) \qquad
 \forall \, i,j\in\ul{L}
\end{align*}
and define the quantities
\begin{align*}
 a_{00} := \norm{u_\WW}_\XX^2 - \Vert \uloc \Vert_\XX^2 - 2 \, \delta, \qquad
 \delta := b(\uloc) - \Vert \uloc \Vert_\XX^2.
\end{align*}

% -----------------------------------------------------------------------------------------
\begin{theorem}[{[P4, Prop.~2]}]\label{thm:compErrorRed}
 Consider the (symmetric) linear system
\begin{align}\label{eq:linsys}
\begin{pmatrix}
a_{00} & \myvec{c}^\top\\
\myvec{c} & \mymatrix{A}
\end{pmatrix}
\begin{pmatrix}
\epsilon\\ \myvec{y}
\end{pmatrix}
=
\begin{pmatrix}
\delta\\
\myvec{b}-\myvec{c}
\end{pmatrix},
\end{align}
where $\epsilon\in\mathbb{R}$ and~$\myvec{y}\in\mathbb{R}^L$. Then, the predicted error reduction can be computed by the formula
\begin{align*} %\label{eq:formula}
 \Delta e_{\WW,\YY}^2
 =\myvec{y}^\top(\myvec{b}-\myvec{c})- \Vert \uloc \Vert^2+\epsilon \, \delta.
 \end{align*}
\end{theorem}
% -----------------------------------------------------------------------------------------

According to the ideas, presented in Section~\ref{sec:implementation}, we discuss the assembling of the quantities $\mymatrix{A}$, $\myvec{b}$ and $\myvec{c}$ in the case that $\XX = \XX(\Omega)$ represents a Hilbert function space over a bounded domain $\Omega\subset\RR^d$. For this purpose, let $\D = \D_h$ be a decomposition of $\Omega$ into closed subsets $K\subseteq\Omega$ such that
\begin{align*}
  \overline{\Omega} 
   = \bigcup_{K\in\D_h} K,
\end{align*}
where $\operatorname{int}(K)\cap \operatorname{int}(K') = \emptyset$ for any $K,K'\in\D_h$ with $K\neq K'$ and set $\XX_K := \lbrace v_{\, | \, K} \; ; \; v\in\XX \rbrace$ for $K\in\D_h$. On each physical element $K\in\D_h$ let $\lbrace \zeta_1^K,\ldots,\zeta_{M_K}^K \rbrace \subseteq \XX_K$ be a set of functions such that there exist representation matrices $\mymatrix{C}_K = (c_{ij}^K)\in\RR^{N\times M_K}$ and $\mymatrix{D}_K = (d_{ij}^K)\in\RR^{L\times M_K}$ satisfying
\begin{align*}
  \phi_{i \, | \, K} 
   = \sum_{j\in\ul{M_K}} c_{ij}^K \, \zeta_j^K, \quad i\in\ul{N}
  \qquad \text{and} \qquad
  \xi_{i \, | \, K}  
   =  \sum_{j\in\ul{M_K}} d_{ij}^K \, \zeta_j^K, \quad i\in\ul{L}.
\end{align*}
If the bilinear form $a(\cdot,\cdot)$ and the linear form $b(\cdot)$ are decomposable in the sense of \eqref{eq:decomp_ab} with local bilinear forms $a_K : \XX_K\times \XX_K\longrightarrow\RR$ and linear forms $b_K : \XX_K\longrightarrow\RR$ on $\XX_K$ for $K\in\D_h$ let us introduce the local matrices $\mymatrix{A}_K = (A_{ij}^K)\in\RR^{M_K\times M_K}$ and the local vectors $\myvec{b}_K = (b_i^K)\in\RR^{M_K}$ component-wise by
\begin{align*}
 A_{ij}^K 
  := a_K \big( \zeta_j^K, \zeta_i^K \big), \qquad
 b_i^K 
  := b_K \big( \zeta_i^K \big) \qquad
 \forall \, i,j\in\ul{M}_K.
\end{align*}
By expressing $\ut\in\WW$ in terms of a linear combination $\ut = \sum_{i\in\ul{N}} u_i \, \phi_i$ we receive the coefficient vector $\wt{\myvec{u}} := (u_1,\ldots,u_N)^{\top}\in\RR^N$, where $u_i = 0$ for $i\in\mathcal{I}^{\text{loc}}$. Therewith, the global quantities arising in the linear system \eqref{eq:linsys} can be assembled element-wise.

% -----------------------------------------------------------------------------------------
\begin{theorem}[{[P4, Prop.~3]}]\label{thm:assemblingFormulars}
 The following identities hold true
\begin{align}\label{eq:assembling_1}
 \mymatrix{A} 
  = \sum_{K\in\D_h} \mymatrix{D}_K \, \mymatrix{A}_K \, \mymatrix{D}_K^{\top}, \qquad
 \myvec{b} 
  = \sum_{K\in\D_h} \mymatrix{D}_K \, \myvec{b}_K, \qquad
 \myvec{c} 
  = \sum_{K\in\D_h} \mymatrix{D}_K \, \mymatrix{A}_K \, \mymatrix{C}_K^{\top} \, \wt{\myvec{u}}.
\end{align}
\end{theorem}

% -----------------------------------------------------------------------------------------

\subsubsection*{Application to the $\boldsymbol{hp}$-Finite Element Context}

In the following, let $\Omega$ be a bounded domain with Lipschitz-boundary $\Gamma := \partial\Omega$, which contains a boundary part $\Gamma_D\subseteq\Gamma$ of positive surface measure, and set
\begin{align*}
 \XX 
  := \lbrace v\in H^1(\Omega) \; ; \; v = 0 \text{ on } \Gamma_D\rbrace.
\end{align*}
Furthermore, let $\Q$ be a decomposition of $\Omega$ into transformed hexahedrons in the sense of [P4, Sec.~3.2] and let $\mathfrak{F}_Q : \wh{Q}\longrightarrow Q$ be the bijective element mapping from the reference element $\wh{Q} := [-1,1]^d$ onto the physical element $Q\in\Q$. Therewith, we introduce the $hp$-finite element space
\begin{align*}
 \WW
  := \big\lbrace v\in\XX \; ; \; v_{\, | \, Q}\circ \mathfrak{F}_Q\in\PP_{p_Q}(\wh{Q}) \text{ for all } Q\in\Q \big\rbrace.
\end{align*}
In order to determine a basis of $\PP_r(\wh{Q})$ for $r\in\NN$, let the functions $\psi_j : [-1,1]\longrightarrow\RR$, for $j\in\NN_0$, be given by
\begin{align*}%\label{eq:defi_1D_psi}
 \psi_0(t) := \frac{1}{2} \, (1-t), \qquad
 \psi_1(t) := \frac{1}{2} \, (1+t), \qquad
 \psi_j(t) := \int_{-1}^s L_{j-1}(s) \dd s, \quad j\geq 2,
\end{align*}
where $L_j:[-1,1]\longrightarrow\RR$ denotes the $j$-th Legendre polynomial, normalized such that $L_j(-1)=(-1)^j$ for $j\geq 1$, see [P4, Sec.~3.1]. By means of the one-dimensional functions $\psi_j$, for any multi-index $\myind{j} = (j_1,\ldots,j_d)\in\NN_0^d$, we define the functions $\wh{\psi}_{\hspace{0.1em}\myind{j}} : \wh{Q}\longrightarrow \RR$ by
\begin{align*}
 \wh{\psi}_{\hspace{0.1em}\myind{j}}(\myvec{x}) 
  := \prod_{k\in\ul{d}} \psi_{j_k}(x_k), \qquad
 \myvec{x} = (x_1,\ldots,x_d)^{\top}\in\wh{Q},
\end{align*}
for which we obtain $\PP_r(\wh{Q}) = \operatorname{span}\lbrace \wh{\psi}_{\hspace{0.1em}\myind{j}} \; ; \; \myind{j}\in\NN_0^d \text{ with } 0 \leq j_k \leq r \rbrace$.
\vspace{0.1cm}

\begin{figure}[ht]
\centering
\includegraphics[trim = 2cm 2cm 3cm 2.75cm, clip, width = 12.25cm]{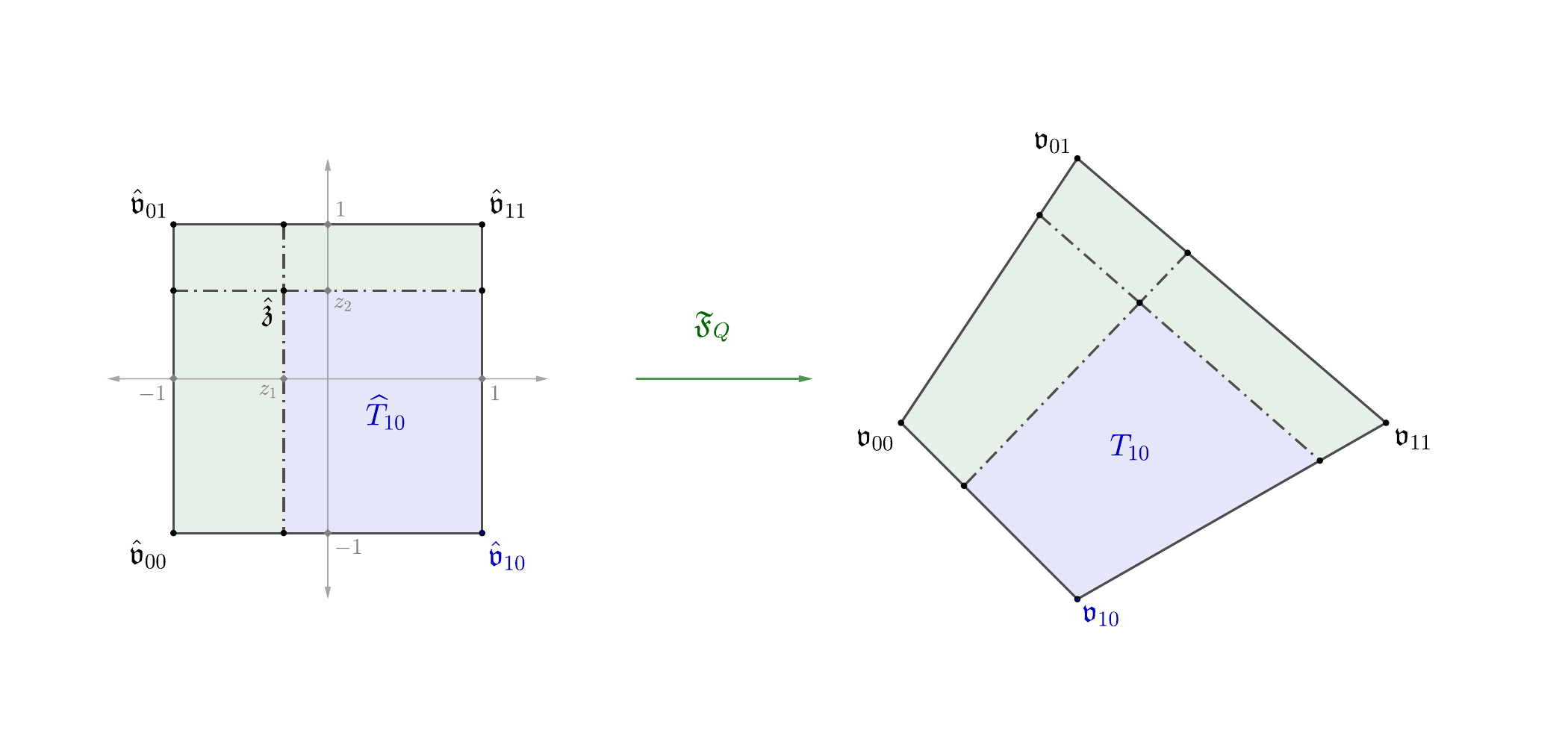}
\caption{\it Refinement of $\wh{Q}$ with respect to $\hat{\myvec{z}}\in(-1,1)^2$ and corresponding refinement of $Q = \mathfrak{F}_Q(\wh{Q})$.}
\label{fig:refinement}
\end{figure}

In [P4, Sec.~3.4], we explicitly construct enrichment functions $\xi_1,\ldots,\xi_L$ that in each case can be associated with exactly one element of the mesh $\Q$. For this purpose, let us focus on some $Q\in\Q$ with local polynomial degree $p_Q$ and let $\R(\wh{Q})$ be a refinement of the reference element $\wh{Q}$ with respect to some dividing point $\hat{\myvec{z}}\in(-1,1)^d$ as specified in [P4, Sec.~3.4] and illustrated in Figure~\ref{fig:refinement} for the two-dimensional case. 

% -----------------------------------------------------------------------------------------
%
%
% -----------------------------------------------------------------------------------------

\begin{figure}[ht]
	\centering
	\includegraphics[scale=0.185]{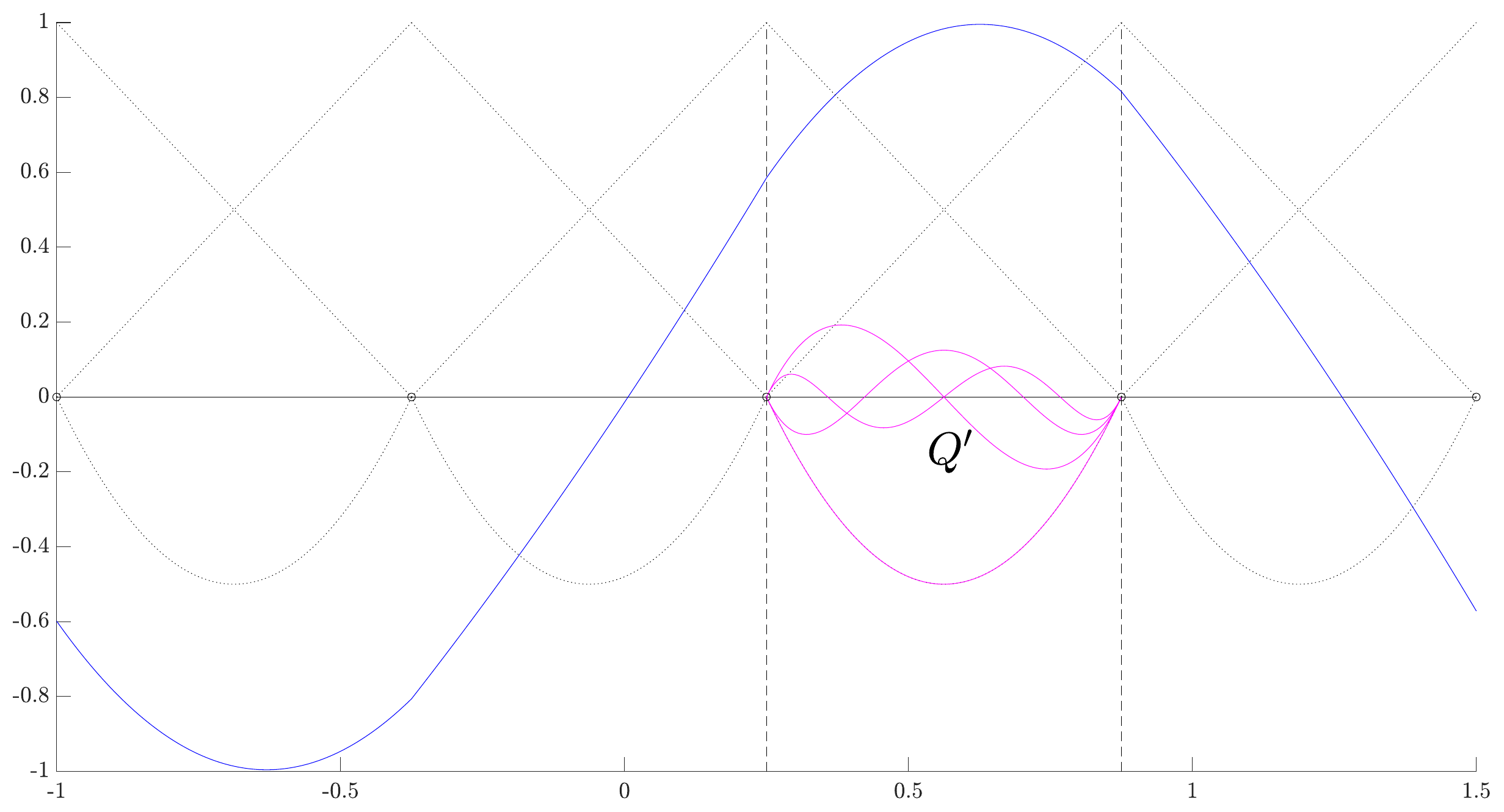} \qquad\quad
	\includegraphics[scale=0.185]{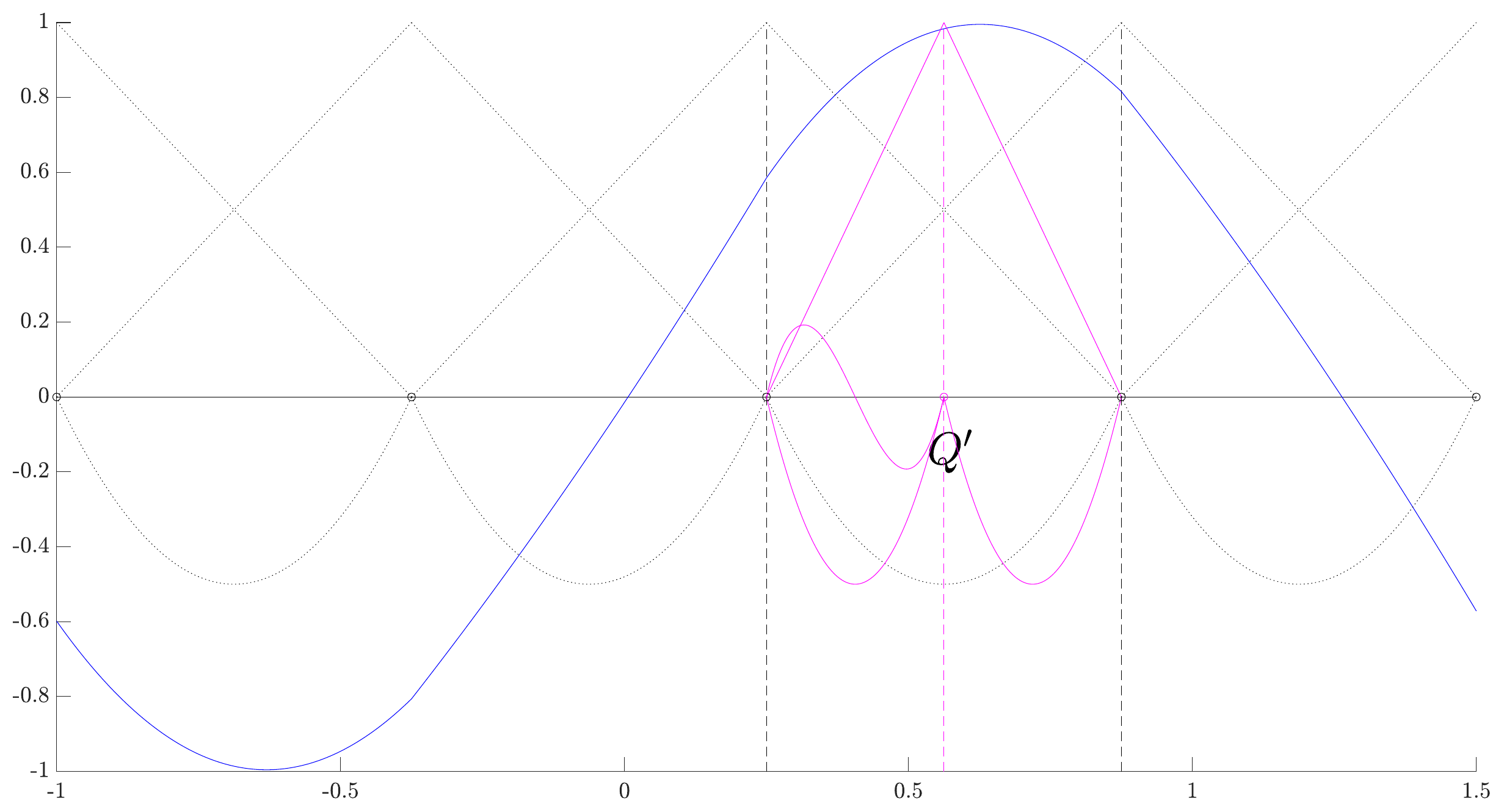}
	\caption{\it A $p$-enrichment (left) vs.~an $hp$-refinement (right).}
	\label{fig:p_vs_h_enrichment}
\end{figure}

We introduce two different types of enrichment functions on $Q$:
\begin{enumerate}[(i)]
 \item For the definition of \emph{$p$-enrichment functions on $Q$}, we consider polynomials on the reference element $\widehat{Q}$, with polynomial degrees larger than $p_{Q}$, and transform them to the physical element $Q$.
 
 \item For the construction of \emph{$hp$-enrichment functions on $Q$}, we consider polynomials on the $2^d$ sub-hexahedra $\wh{T}_{\myind{i}}\in\R(\wh{Q})$ of a refinement $\R(\wh{Q})$ of $\wh{Q}$ with respect to some $\hat{\myvec{z}}\in(-1,1)^d$, which are then transformed to $Q$.
\end{enumerate}
These two scenarios of a $p$-enrichment and an $hp$-refinement are illustrated for an one-dimensional element in Figure~\ref{fig:p_vs_h_enrichment}. Thereby, the discrete solution $u_{\WW}$ is higlighted in blue, the in each case four enrichment functions are depicted in magenta and the basis $\lbrace \phi_1,\ldots,\phi_N\rbrace$ of $\WW$ is indicated in dotted lines, respectively. The functions resulting from transforming polynomials on $\wh{Q}$ or on sub-hexahedra of $\wh{Q}$ to $Q$ will be termed \emph{transformed polynomials}.

% -----------------------------------------------------------------------------------------
%
%
% -----------------------------------------------------------------------------------------

For any multi-index $\myind{j}\in\NN_0^d$ let us introduce the functions $\xi_{\hspace{0.1em}\myind{j}} : \Omega\longrightarrow\RR$ by
\begin{align*} %\label{eq:defi_xi_p_enrichment}
 \xi_{\hspace{0.1em}\myind{j}}(\myvec{x}) 
  := \begin{cases} 
       \wh{\psi}_{\hspace{0.1em}\myind{j}} \circ \myvec{F}_{Q}^{-1}(\myvec{x}), & \text{if } \myvec{x}\in Q, \\
       0, & \text{if } \myvec{x}\in\Omega\setminus Q.
 \end{cases}
\end{align*}
Then, for a \emph{$p$-enrichment on $Q$} finitely many of these transformed polynomials are chosen. Thereby, only such transformed polynomials are considered that vanish along the boundary of $Q$, i.e.~the \emph{$p$-enrichment functions} are given by a set $\mathfrak{E}_p = \lbrace \xi_{\hspace{0.1em}\myind{j}} \; ;\; \myind{j}\in\myind{J}_d\rbrace$, determined by some multi-index set
\begin{align*}%\label{eq:indexSet_pEnr}
 \myind{J}_d \subseteq \big\lbrace \myind{j}\in\NN_0^d \; ; \; j_k \geq 2 \text{ for } k\in\ul{d} \big\rbrace, \qquad
 \abs{\myind{J}_d} = L < \infty.
\end{align*}
Note that any $\xi\in\mathfrak{E}_p$ is continuous on $\Omega$ and has support $Q$ by construction, see [P4, Prop.~4].
\vspace{0.1cm}

In order to define $hp$-enrichment functions on $Q$, let $\R(Q) := \lbrace T_{\myind{i}} \; ; \; \myind{i}\in\lbrace 0,1\rbrace^d \rbrace$ with $T_{\myind{i}} := \myvec{F}_Q(\wh{T}_{\myind{i}})$ be the refinement of $Q$ corresponding to the refinement $\R(\wh{Q})$ of the reference element $\wh{Q}$, cf.~Figure \ref{fig:refinement}. Then, on each sub-hexahedra $T_{\myind{i}}\in\R(Q)$ we introduce, for $\myind{j}\in\NN_0^d$, the functions $\zeta_{\hspace{0.1em}\myind{j}}^{\hspace{0.1em}\myind{i}} : T_{\myind{i}}\longrightarrow\RR$ by
\begin{align*}
 \zeta_{\hspace{0.1em}\myind{j}}^{\hspace{0.1em}\myind{i}}(\myvec{x})
  := \wh{\psi}_{\hspace{0.1em}\myind{j}} \circ \myvec{F}_{\myind{i}}^{-1}(\myvec{x}), \qquad
 \myvec{x}\in T_{\myind{i}},
\end{align*}
where $\myvec{F}_{\myind{i}} : \wh{Q}\longrightarrow T_{\myind{i}}$ denotes the bijective mapping from the reference element onto the sub-hexahedra $T_{\myind{i}}$. The $hp$-enrichment functions are now constructed by means of the transformed polynomials $\zeta_{\myind{j}}^{\myind{i}}$ in such a way that each of them can be associated with a $r$-dimensional node of the mesh $\R(Q)$ which \emph{does not} lie on the boundary $\partial Q$ ($r\in\lbrace 0,\ldots, d\rbrace$). We call such nodes \emph{internal nodes} and point out that each $r$-dimensional internal node can be characterized by choosing a so-called \emph{orientation tuple}
\begin{align*}
 \myind{a} \in D_r := \big\lbrace (a_1,\ldots,a_r)\in\ul{d}^r \; ; \; a_1 < \cdots < a_r \big\rbrace
\end{align*}
and a \emph{location tuple} $\ell = (\ell_1,\ldots,\ell_r)\in\lbrace 0,1\rbrace^r$. Therefore, each internal node may be identified with a tuple $\myind{n} = (\myind{a},\ell)$, see [P4, Sec.3.4.3] for the details. To specify an $hp$-enrichment function $\xi_{\hspace{0.1em}\myind{n},\myind{p}}$ on a $r$-dimensional internal node $\myind{n}=(\myind{a},\ell)$ one has to choose a certain polynomial distribution $\myind{p} = (p_1,\ldots,p_r)$ with $p_k\geq 2$ for $k\in\ul{r}$, containing the polynomial degrees for each of the $r$ directions $k\in A(\myind{a}) := \lbrace a_k \; ; \; k\in\ul{r} \rbrace$ of the $r$-dimensional object $\myind{n}$. Then, by introducing the set
\begin{align*}
 \myind{I}(\myind{n})
  := \big\lbrace \myind{i}=(i_1,\ldots,i_d)\in\lbrace 0,1\rbrace \; ; \; i_k = \ell_k \text{ for each } k\in A(\myind{a}) \big\rbrace
\end{align*}
and, for each $\myind{i}\in\myind{I}(\myind{n})$, the $d$-tuple $\myind{j}(\myind{i},\myind{p}) = (j_1,\ldots,j_d)\in\NN_0^d$ component-wise by
\begin{align*}
 j_k := \begin{cases}
  p_k, & \text{if }k\in A(\myind{a}), \\
  1-i_k, & \text{if } k\notin A(\myind{a}),
 \end{cases}
\end{align*}
the $hp$-enrichment function $\xi_{\myind{n},\myind{p}} : \Omega\longrightarrow\RR$ associated with the internal node $\myind{n}$ and the corresponding polynomial distribution $\myind{p}$ is defined as
\begin{align*}
 \xi_{\myind{n},\myind{p}}(\myvec{x}) 
  := \begin{cases}
  \zeta_{\hspace{0.1em}\myind{j}(\myind{i},\myind{p})}^{\hspace{0.1em}\myind{i}}(\myvec{x}), & \text{if } \myvec{x}\in T_{\myind{i}} ,\text{ for } \myind{i}\in\myind{I}(\myind{n}), \\
  0, & \text{if } \myvec{x} \in \Omega \setminus T(\myind{n}),
 \end{cases}
\end{align*}
where $T(\myind{n}) := \bigcup_{\myind{i}\in\myind{I}(\myind{n})} T_{\myind{i}}$. By choosing, for each internal node $\myind{n}$, a finite set
\begin{align*}
 \myind{P}(\myind{n}) \subseteq \big\lbrace \myind{p}\in\NN_0^r \; ; \; p_k\geq 2 \text{ for } k\in\ul{r} \big\rbrace
\end{align*}
of polynomial distributions one can determine a set $\mathfrak{E}_{hp,\myind{n}} := \lbrace \xi_{\hspace{0.1em}\myind{n},\myind{p}} \; ; \; \myind{p}\in\myind{P}(\myind{n})\rbrace$ of \emph{$hp$-enrichment functions} associated with the node $\myind{n}$. An \emph{$hp$-enrichment on $Q$} is then given by the set
\begin{align*}
 \mathfrak{E}_{hp}
  := \bigcup_{\myind{n}\in\mathcal{N}} \mathfrak{E}_{hp,\myind{n}},
\end{align*}
where $\mathcal{N}$ denotes the set of all internal nodes $\myind{n}$ of the refinement $\R(Q)$.

% -----------------------------------------------------------------------------------------

\begin{theorem}[{[P4, Prop.~5]}]
 Any $\xi_{\myind{n},\myind{p}}\in\mathfrak{E}_{hp}$ is continuous in $\Omega$, and it holds $\operatorname{supp}(\xi_{\myind{n},\myind{p}}) = T(\myind{n})$.
\end{theorem}

% -----------------------------------------------------------------------------------------

Recall that the global quantities $\mymatrix{A}$, $\myvec{b}$ and $\myvec{c}$ for the computation of the predicted error reduction, cf.~Theorem~\ref{thm:compErrorRed}, can be assembled element-wise by the formulas of Theorem~\ref{thm:assemblingFormulars}. As both, the $p$-enrichment as well as the $hp$-enrichment functions on $Q$ have support only on $Q$, by representing the $p$-enrichment functions in terms of the functions $\zeta_{\hspace{0.1em}\myind{j}}^{\hspace{0.1em}\myind{i}}$, the summation over all $Q\in\Q$ in the assembling procedure reduces to a sum over the sub-hexahedrons $T_{\myind{i}}$ of the refinement $\R(Q)$. Obviously, this reduces the computational effort for calculating the predicted error reduction for an enrichment on $Q$. As the involved representation matrices only have to be specified for the sub-hexahedrons $T_{\myind{i}}\in\R(Q)$ for $\myind{i}\in\lbrace 0,1\rbrace^d$ we have
\begin{align*}
 \mymatrix{A} 
  = \sum_{T_{\myind{i}}\in\R(Q)} \mymatrix{D}_{\myind{i}} \, \mymatrix{A}_{\myind{i}} \, \mymatrix{D}_{\myind{i}}^{\top}, \qquad
 \myvec{b} 
  = \sum_{T_{\myind{i}}\in\R(Q)} \mymatrix{D}_{\myind{i}} \, \myvec{b}_{\myind{i}}, \qquad
 \myvec{c} 
  = \sum_{T_{\myind{i}}\in\R(Q)} \mymatrix{D}_{\myind{i}} \, \mymatrix{A}_{\myind{i}} \, \mymatrix{C}_{\myind{i}}^{\top} \, \wt{\myvec{u}}.
\end{align*}

By using the \emph{constraint coefficient technique}, which is presented in [P4, Sec.~3.3], one can compute the representation matrices $\mymatrix{C}_{\myind{i}}$ and $\mymatrix{D}_{\myind{i}}$ in a highly efficient way, which reduces the computational effort even further. In order to compute the matrices $\mymatrix{C}_\myind{i}$ for $\myind{i}\in\lbrace 0,1\rbrace^d$, first the restrictions of the degrees of freedom $\phi_1,\ldots,\phi_N$ to $Q$ are represented in terms of transformed polynomials $\psi_{\hspace{0.1em}\myind{j}}^Q = \wh{\psi}_{\hspace{0.1em}\myind{j}} \circ \myvec{F}_Q^{-1}$ leading to the representation matrix $\mymatrix{C}_Q$. Secondly, for each $\myind{i}\in\lbrace 0,1\rbrace^d$, matrices $\mymatrix{B}_{\myind{i}}$ containing the constraint coefficients to represent the restrictions of the functions $\psi_{\hspace{0.1em}\myind{j}}^Q$ to the sub-hexahedras $T_{\myind{i}}$ in terms of the functions $\zeta_{\hspace{0.1em}\myind{j}}^{\hspace{0.1em}\myind{i}}$ have to be defined, see [P4, Sec.~3.5] for the details. Therewith, one obtains the following result.

% -----------------------------------------------------------------------------------------

\begin{theorem}[{[P4, Prop.~6]}]
 For any $\myind{i}\in\lbrace 0,1\rbrace^d$ the formula $\mymatrix{C}_{\myind{i}} =  \mymatrix{C}_{Q} \, \mymatrix{B}_{\myind{i}}$ holds true.
\end{theorem}

% -----------------------------------------------------------------------------------------

The computation of the representation matrices $\mymatrix{D}_{\myind{i}}$ depends on the kind of enrichment, which is considered on $Q$. In the case of an $p$-enrichment on $Q$ the components of the matrices $\mymatrix{D}_{\myind{i}}$ are given by constraint coefficients, cf.~[P4, Prop.~7]. For an $hp$-enrichment on $Q$ the components are either zero or one, cf.~[P4, Prop.~8].
\vspace{0.1cm}

We emphasize that the representation matrices have to be computed \emph{once only} if different $p$- and $hp$-enrichments are compared on an element $Q$, which makes the comparison of different enrichment strategies cheap. Furthermore, the computation of the predicted error reductions for the different $p$- and $hp$-enrichments on $Q$ can be done in parallel. The resulting $hp$-adaptive algorithm is presented in detail in [P4, Sec.~4] and follows the following steps:
\begin{enumerate}[(i)]
\item \emph{Solving step}. -- The Galerkin approximation $u_\WW$ on a current mesh is computed.

\item \emph{Prediction step}. -- On each element $Q$ of the current mesh the (locally supported) space $\WW_Q^{\text{loc}}$ is constructed to decompose the current solution in terms of $u_\WW = \ut + \uloc$.
\begin{enumerate}
 \item[--] The predicted error reductions $\Delta e_{p,i}$ for finitely many different $p$-enrichments on $Q$ are computed and their maximum $\Delta e_{p,\max}^Q$ is determined. 
 \item[--] The predicted error reductions $\Delta e_{hp,i}$ for finitely many different $hp$-enrichments on $Q$ are computed and their maximum $\Delta e_{hp,\max}^Q$ is determined.
\end{enumerate}
Finally, let $\Delta e_{\max}^Q := \max\big\lbrace \Delta e_{p,\max}^Q, \Delta e_{hp,\max}^Q \big\rbrace$.

\item \emph{Marking step}. -- Mark a subset of the current mesh elements to be flagged for an enrichment.

\item \emph{Enrichment step}. -- On each of the marked elements apply that $p$- or $hp$ enrichment leading to the maximal error reduction $\Delta e_{\max}^Q$.
\end{enumerate}

After sufficiently many iterations of these steps, starting with some initial mesh of $\Omega$, the final solution $u_\WW$ is outputted, cf.~[P4, Alg.~1].

%\newpage
%%\thispagestyle{plain}
%\quad
%\newpage

% -----------------------------------------------------------------------------------------
%   CHAPTER :  Appendix
% -----------------------------------------------------------------------------------------

\chapter{Appendix}

Here, the basic properties of the bilinear form $a(\cdot,\cdot)$, appearing in the weak formulation \eqref{eq:varIQ_of_second_kind} of the boundary value problem \eqref{eq:PDE_model_problem}--\eqref{eq:PFL_model_problem}, and the corresponding energy functional $\E(\cdot)$, cf.~\eqref{eq:defi_energyF}, are shown. Though these properties are generally known and similar calculations can be found, for instance, in \cite{ref:Han_2013}, we present the proofs for our specific model problem and in our notation.

% -----------------------------------------------------------------------------------------

\subsubsection*{Preliminaries}

Note that by the definition of the norm $\norm{\cdot}$ on $\VV\times Q$, cf.~Section~\ref{sec:model_problem}, the subadditivity of the square root together with the Cauchy-Schwarz inequality yield
\begin{align}\label{eq:norm_est}
 \norm{(\myvec{v},\mymatrix{q})}
  \leq \norm{\myvec{v}}_{1,\Omega} + \norm{\mymatrix{q}}_{0,\Omega}
  \leq \sqrt{2} \, \norm{(\myvec{v},\mymatrix{q})} \qquad
 \forall \, (\myvec{v},\mymatrix{q})\in\VV\times Q.
\end{align}
Moreover, recall that for $a,b\in\RR$ it holds $2ab \leq a^2 + b^2$, from which immediately follows that
\begin{align*}
 2ab
  \leq \delta \, a^2 + \frac{1}{\delta} \, b^2 \qquad
 \forall \, a,b\in\RR \quad \forall \, \delta>0
\end{align*}
and therefore
\begin{align}\label{eq:hilfsIQ_01}
 2 \, \mymatrix{\varepsilon}(\myvec{v}) : \mymatrix{q}
  \leq \delta \, \abs{\mymatrix{\varepsilon}(\myvec{v})}_F^2 + \frac{1}{\delta} \, \abs{\mymatrix{q}}_F^2 \qquad
 \forall \, (\myvec{v},\mymatrix{q})\in\VV\times Q.
\end{align}
If $\Omega$ represents a bounded domain with Lipschitz-boundary $\Gamma := \partial\Omega$, then by \emph{Korn's first inequality} there exists a constant $c_{\Omega}>0$ (only depending on the domain $\Omega$) such that
\begin{align}\label{eq:first_kornIQ}
 \norm{\myvec{v}}_{H^1(\Omega,\RR^d)}^2
  \leq c_{\Omega} \, \norm{ \mymatrix{\varepsilon}(\myvec{v}) }_{0,\Omega}^2 \qquad
 \forall \, \myvec{v}\in\VV.
\end{align}
Due to \eqref{eq:first_kornIQ} the norms $\norm{\cdot}_{H^1(\Omega,\RR^d)}$ and $\norm{\cdot}_{1,\Omega} = \big( \norm{\myvec{v}}_{0,\Omega}^2 + \norm{\mymatrix{\varepsilon}(\myvec{v}) }_{0,\Omega}^2 \big)^{1/2}$, c.f.~Section~\ref{sec:model_problem}, are equivalent on $\VV$. Furthermore, for a connected boundary part $\Gamma'\subseteq\Gamma$ by the \emph{trace theorem} there exists a constant $c_{tr}>0$ with
\begin{align}\label{eq:trace_thm}
 \norm{\myvec{v}}_{\frac{1}{2}, \Gamma'} 
  \leq c_{tr} \, \norm{\myvec{v}}_{1,\Omega} \qquad
 \forall \, \myvec{v}\in\VV.
\end{align}
Finally, if $\abs{\KK} := \sum_{i,j,k,l\in\underline{d}} \KK_{ijkl}$ for any forth-order tensor $\KK$ then there holds
\begin{align}\label{eq:hilfsIQ_02}
 \abs{\KK \, \mymatrix{\tau}}_F 
  \leq \abs{\KK} \, \abs{\mymatrix{\tau}}_F \qquad
 \forall \, \mymatrix{\tau}\in\RR^{d\times d}.
\end{align}

% -----------------------------------------------------------------------------------------

\subsubsection*{Properties of the Bilinear Form $\boldsymbol{a(\cdot,\cdot)}$}

The \emph{symmetry} of $a(\cdot,\cdot)$ directly follows from the required symmetry of the elasticity tensor $\CC$ and the hardening tensor $\HH$ as by definition it holds
\begin{align*}
 a\big( (\myvec{u},\mymatrix{p}), (\myvec{v},\mymatrix{q}) \big) 
  = \int_{\Omega} \CC \, (\mymatrix{\varepsilon}(\myvec{u})-\mymatrix{p}) : (\mymatrix{\varepsilon}(\myvec{v})-\mymatrix{q}) \dd \myvec{x} + \int_{\Omega} \HH \,\mymatrix{p} : \mymatrix{q} \dd \myvec{x}
\end{align*}
for any $(\myvec{u},\mymatrix{p}), (\myvec{v},\mymatrix{q})\in\VV\times Q$ and $\KK \, \mymatrix{p} : \mymatrix{q} = \KK \, \mymatrix{q} : \mymatrix{p}$ for a symmetric forth-order tensor $\KK$.
\vspace{0.1cm}

Since by assumption it holds $\CC_{ijkl},\HH_{ijkl}\in L^{\infty}(\Omega)$ for all $i,j,k,l\in\underline{d}$ the constants
\begin{align*}
 c_1 
  := \esssup_{\myvec{x}\in\Omega} \abs{\CC(\myvec{x})}, \qquad
 c_2 
  := \esssup_{\myvec{x}\in\Omega} \abs{\HH(\myvec{x})}
\end{align*}
are both finite. Thus, the inequality \eqref{eq:hilfsIQ_02} yields
\begin{align*}
\norm{\mymatrix{\sigma}(\myvec{u},\mymatrix{p})}_{0,\Omega}
 = \bigg( \int_{\Omega} \big| \CC \, (\mymatrix{\varepsilon}(\myvec{u})-\mymatrix{p}) \big|_F^2 \dd \myvec{x} \bigg)^{1/2} 
 &\leq \bigg( \int_{\Omega} \abs{\CC}^2 \, \abs{\mymatrix{\varepsilon}(\myvec{u})-\mymatrix{p}}_F^2 \dd \myvec{x} \bigg)^{1/2} \\
 &\leq \esssup_{\myvec{x}\in\Omega} \abs{\CC(\myvec{x})} \, \norm{\mymatrix{\varepsilon}(\myvec{u})-\mymatrix{p}}_{0,\Omega}  \\
 &= c_1 \, \norm{\mymatrix{\varepsilon}(\myvec{u})-\mymatrix{p}}_{0,\Omega} 
\end{align*}
as well as
\begin{align*}
 \norm{\HH \, \mymatrix{p}}_{0,\Omega}
  = \bigg( \int_{\Omega} \abs{\HH \, \mymatrix{p}}_F^2 \dd \myvec{x} \bigg)^{1/2} 
  \leq \esssup_{\myvec{x}\in\Omega} \abs{\HH(\myvec{x})} \, \norm{\mymatrix{p}}_{0,\Omega}
  = c_2 \, \norm{\mymatrix{p}}_{0,\Omega}.
\end{align*}
Therefore, the \emph{continuity} of the bilinear form $a(\cdot,\cdot)$ can be deduced from
\begin{align*}
a\big( (\myvec{u},\mymatrix{p}), (\myvec{v},\mymatrix{q}) \big)
 &= \big( \mymatrix{\sigma}(\myvec{u},\mymatrix{p}), \mymatrix{\varepsilon}(\myvec{v})-\mymatrix{q} \big)_{0,\Omega} + (\mathbb{H} \, \mymatrix{p}, \mymatrix{q})_{0,\Omega} \\
 &\leq \norm{ \mymatrix{\sigma}(\myvec{u},\mymatrix{p}) }_{0,\Omega} \, \norm{ \mymatrix{\varepsilon}(\myvec{v}) - \mymatrix{q} }_{0,\Omega} + \norm{ \HH \, \mymatrix{p}}_{0,\Omega} \, \norm{ \mymatrix{q} }_{0,\Omega} \\
 &\leq c_1 \, \norm{ \mymatrix{\varepsilon}(\myvec{u})-\mymatrix{p} }_{0,\Omega} \, \norm{ \mymatrix{\varepsilon}(\myvec{v}) - \mymatrix{q} }_{0,\Omega} + c_2 \, \norm{ \mymatrix{p}}_{0,\Omega} \, \norm{ \mymatrix{q} }_{0,\Omega} \\
 &\leq c_1 \, \big( \norm{\myvec{u}}_{1,\Omega} + \norm{\mymatrix{p}}_{0,\Omega} \big) \,  \big( \norm{\myvec{v}}_{1,\Omega} + \norm{\mymatrix{q}}_{0,\Omega} \big) + c_2 \, \norm{\mymatrix{p}}_{0,\Omega} \, \norm{\mymatrix{q}}_{0,\Omega} \\
 &\leq (c_1 + c_2) \, \big( \norm{\myvec{u}}_{1,\Omega} + \norm{\mymatrix{p}}_{0,\Omega} \big) \, \big( \norm{\myvec{v}}_{1,\Omega} + \norm{\mymatrix{q}}_{0,\Omega} \big) \\
 &\leq 2 \, (c_1 + c_2) \, \norm{ (\myvec{u},\mymatrix{p}) } \, \norm{ (\myvec{v},\mymatrix{q})},
\end{align*}
where for the last inequality \eqref{eq:norm_est} was used. Let $c_a := 2 \, (c_1 + c_2)$.
\vspace{0.1cm}

To show the \emph{$(\VV\times Q)$-ellipticity} of $a(\cdot,\cdot)$ first note that
\begin{align*}
 \norm{ \mymatrix{\varepsilon}(\myvec{v})-\mymatrix{q} }_{0,\Omega}^2
 = \int_{\Omega} \abs{ \mymatrix{\varepsilon}(\myvec{v}) - \mymatrix{q} }_F^2 \dd \myvec{x}
 &= \int_{\Omega} \abs{ \mymatrix{\varepsilon}(\myvec{v}) }_F^2 - 2 \, \mymatrix{\varepsilon}(\myvec{v}) : \mymatrix{q} + \abs{ \mymatrix{q} }_F^2 \dd \myvec{x} \\
 &\geq \int_{\Omega} (1-\delta) \, \abs{ \mymatrix{\varepsilon}(\myvec{v}) }_F^2 + \Big(1-\frac{1}{\delta} \Big) \abs{ \mymatrix{q} }_F^2 \dd \myvec{x} \\
 &= (1-\delta) \, \norm{ \mymatrix{\varepsilon}(\myvec{v}) }_{0,\Omega}^2 + \Big(1-\frac{1}{\delta} \Big) \norm{ \mymatrix{q} }_{0,\Omega}^2
\end{align*}
by exploiting the inequality \eqref{eq:hilfsIQ_01}. Therewith, the required uniform ellipticity of $\CC$ and $\HH$ yields
\begin{align*}
 a\big( (\myvec{v},\mymatrix{q}), (\myvec{v},\mymatrix{q}) \big) 
 &= \big( \mymatrix{\sigma}(\myvec{v},\mymatrix{q}), \mymatrix{\varepsilon}(\myvec{v})-\mymatrix{q} \big)_{0,\Omega} + (\mathbb{H} \, \mymatrix{q}, \mymatrix{q})_{0,\Omega} \\
 &= \int_{\Omega} \CC \, (\mymatrix{\varepsilon}(\myvec{v}) - \mymatrix{q}) : (\mymatrix{\varepsilon}(\myvec{v}) - \mymatrix{q}) + \mathbb{H} \, \mymatrix{q} : \mymatrix{q} \dd \myvec{x} \\
 &\geq \int_{\Omega} c_e \, \abs{ \mymatrix{\varepsilon}(\myvec{v}) - \mymatrix{q} }_F^2 + c_h \, \abs{ \mymatrix{q} }_F^2 \dd \myvec{x} \\
 &\geq c_e \, (1-\delta) \, \norm{ \mymatrix{\varepsilon}(\myvec{v}) }_{0,\Omega}^2 + \Big( c_e \Big(1-\frac{1}{\delta} \Big) + c_h \Big) \norm{ \mymatrix{q} }_{0,\Omega}^2.
\end{align*}
Hence, for a sufficient small $\delta$ it follows that
\begin{align}\label{eq:proof_ellipticity01}
 a\big( (\myvec{v},\mymatrix{q}), (\myvec{v},\mymatrix{q}) \big)
  \geq c \, \big( \norm{ \mymatrix{\varepsilon}(\myvec{v}) }_{0,\Omega}^2 + \norm{ \mymatrix{q} }_{0,\Omega}^2 \big) 
\end{align}
with a positive constant $c :=  \, \min \big\lbrace c_e \, (1-\delta), c_e \, (1-\frac{1}{\delta}) + c_h \big\rbrace$. As by the \emph{triangle inequality}
\begin{align*}
 \norm{ \mymatrix{\varepsilon}(\myvec{v}) }_{0,\Omega}^2 
 = \int_{\Omega} \abs{ \mymatrix{\varepsilon}(\myvec{v}) }_F^2 \dd \myvec{x}
 \leq \int_{\Omega} \Big(\frac{1}{2} \, \abs{ \nabla\myvec{v} }_F + \frac{1}{2} \, \abs{ (\nabla\myvec{v})^{\top} }_F \Big)^2 \dd \myvec{x}
 = \int_{\Omega} \abs{ \nabla\myvec{v} }_F^2 \dd \myvec{x}
 = \norm{ \nabla\myvec{v} }_{0,\Omega}^2,
\end{align*}
\emph{Korn's first inequality}, \eqref{eq:first_kornIQ}, gives
\pagebreak
\begin{align*}
 \norm{ \mymatrix{\varepsilon}(\myvec{v}) }_{0,\Omega}^2 
 \geq \frac{1}{c_\Omega} \, \norm{ \myvec{v} }_{H^1(\Omega,\RR^d)}^2
 = \frac{1}{c_\Omega} \, \big( \norm{ \myvec{v} }_{0,\Omega}^2 + \norm{ \nabla\myvec{v} }_{0,\Omega}^2 \big)
 \geq \frac{1}{c_\Omega} \, \big( \norm{ \myvec{v} }_{0,\Omega}^2 + \norm{ \mymatrix{\varepsilon}(\myvec{v}) }_{0,\Omega}^2 \big)
 = \frac{1}{c_\Omega} \, \norm{ \myvec{v} }_{1,\Omega}^2,
\end{align*}
from which together with \eqref{eq:proof_ellipticity01} finally follows that
\begin{align}\label{eq:ellipticity_bilinearF}
 a\big( (\myvec{v},\mymatrix{q}), (\myvec{v},\mymatrix{q}) \big) 
 \geq c \, \big( \norm{ \mymatrix{\varepsilon}(\myvec{v}) }_{0,\Omega}^2 + \norm{ \mymatrix{q} }_{0,\Omega}^2 \big)
 \geq \frac{c}{c_\Omega} \, \norm{ \myvec{v} }_{1,\Omega}^2 + c \, \norm{ \mymatrix{q} }_{0,\Omega}^2
 \geq \alpha \, \norm{ (\myvec{v},\mymatrix{q}) }^2
\end{align}
with $\alpha := \min\lbrace \frac{c}{c_\Omega}, c \rbrace$. Hence, the bilinear form $a(\cdot,\cdot)$ is $(\VV\times Q)$-elliptic.

% -----------------------------------------------------------------------------------------

\subsubsection*{Properties of the Energy Functional $\boldsymbol{\E(\cdot)}$}

For the \emph{coercivity} of $\E(\cdot)$ one has to show that $\E(\myvec{v},\mymatrix{q})\longrightarrow\infty$ for $(\myvec{v},\mymatrix{q})\in\VV\times Q$ with $\norm{(\myvec{v},\mymatrix{q})}\longrightarrow\infty$. First of all, the $(\VV\times Q)$-ellipticity of $a(\cdot,\cdot)$ yields
\begin{align*}
 a\big( (\myvec{v},\mymatrix{q}), (\myvec{v},\mymatrix{q}) \big)
  \geq \alpha \, \norm{ (\myvec{v},\mymatrix{q}) }
  = \alpha \, \big( \norm{\myvec{v}}_{1,\Omega}^2 + \norm{ \mymatrix{q} }_{0,\Omega}^2 \big) \qquad
 \forall \, (\myvec{v},\mymatrix{q})\in\VV\times Q.
\end{align*}
Furthermore, by the definition of the plasticity functional it holds
\begin{align*}
 \psi(\mymatrix{q})
  = \int_{\Omega} \sigma_y \, \abs{\mymatrix{q}}_F \dd \myvec{x} 
  \geq 0 \qquad
 \forall \, \mymatrix{q}\in Q
\end{align*}
and finally, the \emph{trace theorem}, cf.~\eqref{eq:trace_thm}, implies
\begin{align*}
 \ell(\myvec{v})
  = \langle \myvec{f}, \myvec{v} \rangle + \langle \myvec{g}, \myvec{v} \rangle_{\Gamma_N}
  \leq \norm{\myvec{f}}_{\VV^{\star}} \norm{\myvec{v}}_{1,\Omega} + \norm{\myvec{g}}_{-\frac{1}{2},\Gamma_N} \norm{\myvec{v}}_{\frac{1}{2},\Gamma_N}
  \leq \Big( \norm{\myvec{f}}_{\VV^{\star}} + c_{tr} \, \norm{\myvec{g}}_{-\frac{1}{2},\Gamma_N} \Big) \, \norm{\myvec{v}}_{1,\Omega},
\end{align*}
where $\norm{\cdot}_{\VV^{\star}}$ and $\norm{\cdot}_{-\frac{1}{2},\Gamma_N}$ denote the dual norm on $\VV^{\star}$ and $H^{-1/2}(\Omega,\RR^d)$, respectively. By defining the nonnegative constant $c := \norm{\myvec{f}}_{\VV^{\star}} + c_{tr} \, \norm{\myvec{g}}_{-\frac{1}{2},\Gamma_N}$ it therefore follows that 
\begin{align}\label{eq:coercivity_E}
 \E(\myvec{v}, \mymatrix{q})
  \geq \alpha \, \big( \norm{\myvec{v}}_{1,\Omega}^2 + \norm{ \mymatrix{q} }_{0,\Omega}^2 \big) - c \, \norm{\myvec{v}}_{1,\Omega} \qquad
 \forall \, (\myvec{v}, \mymatrix{q})\in\VV\times Q.
\end{align}
Now, as the quadratic term $\norm{\myvec{v}}_{1,\Omega}^2$ grows faster than the linear term $\norm{\myvec{v}}_{1,\Omega}$ for $\norm{\myvec{v}}_{1,\Omega}\longrightarrow\infty$ the inequality \eqref{eq:coercivity_E} shows $\E(\myvec{v}, \mymatrix{q})\longrightarrow\infty$ for $\norm{\myvec{v}}_{1,\Omega}\longrightarrow\infty$ or $\norm{ \mymatrix{q} }_{0,\Omega}\longrightarrow\infty$, which completes the argument.
\vspace{0.1cm}

In order to show the \emph{convexity} of $\E(\cdot)$ first note that for arbitrary $t\in[0,1]$ it holds
\begin{align}\label{eq:hilfsUeberlegung}
 t^2 
  = t \, \big( 1 - (1-t) \big)
  = t - t \, (1-t), \qquad
 (1-t)^2
  = (1-t) - t \, (1-t).
\end{align}
For arbitrary $(\myvec{u}, \mymatrix{p}), (\myvec{v}, \mymatrix{q})\in\VV\times Q$ let $(\myvec{w}, \mymatrix{\tau}) := t \, (\myvec{u}, \mymatrix{p}) + (1-t) \, (\myvec{v}, \mymatrix{q})$. Then, the bilinearity of $a(\cdot,\cdot)$ gives
\begin{align*}
 &a\big( (\myvec{w},\mymatrix{\tau}),(\myvec{w},\mymatrix{\tau}) \big) \\
  & \qquad\quad  = t^2 \, a\big( (\myvec{u},\mymatrix{p}),(\myvec{u},\mymatrix{p}) \big) + t \, (1-t) \, \Big( a\big( (\myvec{u},\mymatrix{p}), (\myvec{v},\mymatrix{q}) \big) + a\big( (\myvec{v},\mymatrix{q}), (\myvec{u},\mymatrix{p}) \big) \Big) + (1-t)^2 \, a\big( (\myvec{v},\mymatrix{q}), (\myvec{v},\mymatrix{q}) \big),
\end{align*}
which equals
\begin{align*}
 &t \, a\big( (\myvec{u},\mymatrix{p}),(\myvec{u},\mymatrix{p}) \big) 
 + t \, (1-t) \, \Big( a\big( (\myvec{u},\mymatrix{p}), (\myvec{v},\mymatrix{q}) \big) - a\big( (\myvec{u},\mymatrix{p}), (\myvec{u},\mymatrix{p}) \big) - a\big( (\myvec{v},\mymatrix{q}), (\myvec{v},\mymatrix{q}) \big) + a\big( (\myvec{v},\mymatrix{q}), (\myvec{u},\mymatrix{p}) \big) \Big) \\
 & \qquad\quad + (1-t) \, a\big( (\myvec{v},\mymatrix{q}), (\myvec{v},\mymatrix{q}) \big)
\end{align*}
due to \eqref{eq:hilfsUeberlegung} and can be simplified by means of the bilinearity of $a(\cdot,\cdot)$ to
\begin{align*}
 t \, a\big( (\myvec{u},\mymatrix{p}),(\myvec{u},\mymatrix{p}) \big) 
 + (1-t) \, a\big( (\myvec{v},\mymatrix{q}), (\myvec{v},\mymatrix{q}) \big)
 - t \, (1-t) \, a\big( (\myvec{u}-\myvec{v},\mymatrix{p}-\mymatrix{q}), (\myvec{u}-\myvec{v},\mymatrix{p}-\mymatrix{q}) \big).
\end{align*}
Hence, together with the $(\VV\times Q)$-ellipticity of $a(\cdot,\cdot)$ it follows that
\begin{align*}
 a\big( (\myvec{w},\mymatrix{\tau}),(\myvec{w},\mymatrix{\tau}) \big)
  &\leq t \, a\big( (\myvec{u},\mymatrix{p}),(\myvec{u},\mymatrix{p}) \big) 
   + (1-t) \, a\big( (\myvec{v},\mymatrix{q}), (\myvec{v},\mymatrix{q}) \big)
   - \alpha \, t \, (1-t) \, \norm{(\myvec{u} - \myvec{v}), \mymatrix{p} - \mymatrix{q}} \\
  &\leq t \, a\big( (\myvec{u},\mymatrix{p}),(\myvec{u},\mymatrix{p}) \big) 
   + (1-t) \, a\big( (\myvec{v},\mymatrix{q}), (\myvec{v},\mymatrix{q}) \big).
\end{align*}
Moreover, by the convexity of $\psi(\cdot)$ and the linearity of $\ell(\cdot)$ it holds
\pagebreak
\begin{align*}
 \psi \big( t \, \mymatrix{p} + (1-t) \, \mymatrix{q} \big)
  \leq t \, \psi(\mymatrix{p}) + (1-t) \, \psi(\mymatrix{q}), \qquad
 \ell \big( t \, \myvec{u} + (1-t) \, \myvec{v} \big)
  = t \, \ell(\myvec{u}) + (1-t) \, \ell(\myvec{v}).
\end{align*}
Overall, this yields
\begin{align*}
 \E(\myvec{w}, \mymatrix{\tau})
  &= \frac{1}{2} \, a\big( (\myvec{w}, \mymatrix{\tau}), (\myvec{w}, \mymatrix{\tau}) \big) + \psi \big( t \, \mymatrix{p} + (1-t) \, \mymatrix{q} \big) - \ell \big( t \, \myvec{u} + (1-t) \, \myvec{v} \big) \\
  &\leq \frac{t}{2} \, a\big( (\myvec{u},\mymatrix{p}), (\myvec{u},\mymatrix{p}) \big) 
   + \frac{1-t}{2} \, a\big( (\myvec{v},\mymatrix{q}), (\myvec{v},\mymatrix{q}) \big) 
   + t \, \psi(\mymatrix{p}) + (1-t) \, \psi(\mymatrix{q})
   - t \, \ell(\myvec{u}) + (1-t) \, \ell(\myvec{v}) \\
  &= t \, \E(\myvec{u},\mymatrix{p}) + (1-t) \, \E(\myvec{v},\mymatrix{q}),
\end{align*}
which leads the convexity of $\E(\cdot)$.
\vspace{0.1cm}

The functional $\E(\cdot)$ is \emph{subdifferentiable} on $\VV\times Q$ if and only if for any $(\myvec{u},\mymatrix{p})\in\VV\times Q$ there exists a $\partial\E(\myvec{u},\mymatrix{p})\in(\VV\times Q)^{\star}$ such that
\begin{align*}
 \E(\myvec{v},\mymatrix{q})
  \geq \E(\myvec{u},\mymatrix{p}) + \big\langle \partial\E(\myvec{u},\mymatrix{p}), (\myvec{v} - \myvec{u}, \mymatrix{q} - \mymatrix{p}) \big\rangle \qquad
 \forall \, (\myvec{v},\mymatrix{q})\in \VV\times Q.
\end{align*}
First of all, the bilinearity, symmetry and $(\VV\times Q)$-ellipticity of $a(\cdot,\cdot)$ gives
\begin{align*}
 0 
  \leq a\big( (\myvec{u}-\myvec{v}, \mymatrix{p}-\mymatrix{q}), (\myvec{u}-\myvec{v}, \mymatrix{p}-\mymatrix{q}) \big) 
  = a\big( (\myvec{u},\mymatrix{p}), (\myvec{u},\mymatrix{p}) \big) - 2 \, a\big( (\myvec{u},\mymatrix{p}), (\myvec{v},\mymatrix{q}) \big) + a\big( (\myvec{v},\mymatrix{q}), (\myvec{v},\mymatrix{q}) \big),
\end{align*}
from which immediately follows that
\begin{align*}
 a\big( (\myvec{u},\mymatrix{p}), (\myvec{v},\mymatrix{q}) \big)
  \leq \frac{1}{2} \, a\big( (\myvec{u},\mymatrix{p}), (\myvec{u},\mymatrix{p}) \big) + \frac{1}{2} \, a\big( (\myvec{v},\mymatrix{q}), (\myvec{v},\mymatrix{q}) \big)
\end{align*}
and subtracting $a\big( (\myvec{u},\mymatrix{p}), (\myvec{u},\mymatrix{p}) \big)$ one both sides finally yields
\begin{align*}
 a\big( (\myvec{u},\mymatrix{p}), (\myvec{v}-\myvec{u},\mymatrix{q}-\mymatrix{p}) \big)
  \leq \frac{1}{2} \, a\big( (\myvec{v},\mymatrix{q}), (\myvec{v},\mymatrix{q}) \big) - \frac{1}{2} \, a\big( (\myvec{u},\mymatrix{p}), (\myvec{u},\mymatrix{p}) \big).
\end{align*}
Thus, together with the subdifferentiability of $\psi(\cdot)$ it follows that
\begin{align*}
 \E(\myvec{v},\mymatrix{q}) - \E(\myvec{u},\mymatrix{p})
  &= \frac{1}{2} \Big( a\big( (\myvec{v},\mymatrix{q}), (\myvec{v},\mymatrix{q}) \big) - a\big( (\myvec{u},\mymatrix{p}), (\myvec{u},\mymatrix{p}) \big) \Big) + \psi(\mymatrix{q}) - \psi(\mymatrix{p}) - \ell(\myvec{v}) + \ell(\myvec{u}) \\
&\geq \big\langle A(\myvec{u},\mymatrix{p}), (\myvec{v}-\myvec{u},\mymatrix{q}-\mymatrix{p}) \big\rangle + \big\langle \partial\psi(\mymatrix{p}), \mymatrix{q}-\mymatrix{p} \big\rangle - \langle \ell, \myvec{v}-\myvec{u} \rangle,
\end{align*}
where $A(\myvec{u},\mymatrix{p})\in(\VV\times Q)^{\star}$ denotes the uniquely determined operator that is characterized by
\begin{align*}
 \big\langle A(\myvec{u},\mymatrix{p}), (\myvec{v},\mymatrix{q}) \big\rangle
  = a\big( (\myvec{u},\mymatrix{p}), (\myvec{v},\mymatrix{q}) \big) \qquad
 \forall \, (\myvec{v},\mymatrix{q})\in\VV\times Q.
\end{align*}
Note that the existence of the operator $A(\myvec{u},\mymatrix{p})$ for $(\myvec{u},\mymatrix{p})\in\VV\times Q$ is guaranteed by the continuity of the bilinear form $a(\cdot,\cdot)$, see e.g.~\cite{ref:Hackbusch_2017}. Hence, the subdifferentiability of $\E(\cdot)$ follows by defining $\partial\E(\cdot)$ as
\begin{align*}
 \partial\E(\myvec{u},\mymatrix{p})
  := \big\langle A(\myvec{u},\mymatrix{p}), (\myvec{v}-\myvec{u},\mymatrix{q}-\mymatrix{p}) \big\rangle + \big\langle \partial\psi(\mymatrix{p}), \mymatrix{q}-\mymatrix{p} \big\rangle - \langle \ell, \myvec{v}-\myvec{u} \rangle \qquad
 \forall \, (\myvec{v},\mymatrix{q})\in \VV\times Q.
\end{align*}
%\vspace{0.1cm}

Finally, to show that $\E(\cdot)$ is \emph{weakly lower semi-continuous} let $(\myvec{v}_n, \mymatrix{q}_n)\rightharpoonup (\myvec{v}, \mymatrix{q})$ for $n\to\infty$, which means that for any $G^{\star}\in(\VV\times Q)^{\star}$ it follows that $\langle G^{\star}, (\myvec{v}_n,\mymatrix{q}_n) \rangle\longrightarrow \langle G^{\star}, (\myvec{v},\mymatrix{q})\rangle$ as $n\to\infty$. Thus, the subdifferetiability of $\E(\cdot)$ implies
\begin{align*}
 \E(\myvec{v}_n, \mymatrix{q}_n) - \E(\myvec{v}, \mymatrix{q})
  \geq \big\langle \partial\E(\myvec{v}, \mymatrix{q}), (\myvec{v}_n-\myvec{v}, \mymatrix{q}_n - \mymatrix{q}) \big\rangle \qquad
 \forall \, n\in\NN.
\end{align*}
Taking the limit inferior on both sides therefore yields
\begin{align*}
 \liminf_{n\to\infty} \E(\myvec{v}_n, \mymatrix{q}_n) - \E(\myvec{v}, \mymatrix{q})
  \geq \liminf_{n\to\infty} \big\langle \partial\E(\myvec{v}, \mymatrix{q}), (\myvec{v}_n-\myvec{v}, \mymatrix{q}_n-\mymatrix{q}) \big\rangle 
  = 0,
\end{align*}
which implies $\liminf_{n\to\infty} \E(\myvec{v}_n, \mymatrix{q}_n) \geq \E(\myvec{v}, \mymatrix{q})$ and, thus, the weakly lower semi-continuity of $\E(\cdot)$.

%\newpage
%%\thispagestyle{plain}
%\quad
%\newpage

\newpage

% -----------------------------------------------------------------------------------------
%    LITERATUR
% -----------------------------------------------------------------------------------------

\addcontentsline{toc}{chapter}{Bibiliography}
\bibliographystyle{plaindin}
\bibliography{bac}

\newpage

{\small
\listoffigures 
}

All Figures are self-generated (the one on page nine with the software \emph{Matlab} and the remaining ones with the Open-Source-Software \emph{GeoGebra}).

\newpage

% -----------------------------------------------------------------------------------------
%    Paper 1  -  hp-Finite Elements with Decoupled Constraints for Elastoplasticity
% -----------------------------------------------------------------------------------------

\chapter{Publication P1}

{\large
\begin{center}
$hp$-Finite Elements with Decoupled Constraints for Elastoplasticity
\end{center}
}
\vspace{1cm}
\textbf{Published in:} \emph{Spectral and High Order Methods for Partial Differential Equations ICOSAHOM 2020+1}. \\

\copyright~2023 Springer %, reprinted with permission.

% -----------------------------------------------------------------------------------------
%    Paper 2  -  hp-FEM for a Mixed Variational Approach in Elastoplasticity
% -----------------------------------------------------------------------------------------

\chapter{Publication P2}

{\large
\begin{center}
Mixed Finite Elements of Higher-Order in Elastoplasticity
\end{center}
}
\vspace{1cm}
\textbf{Submitted to:} \emph{Applied Numerical Mathematics} (under review). \\

Available as arXiv-preprint \emph{arXiv:2401.09080} (2024).

% -----------------------------------------------------------------------------------------
%    Paper 3  -  A Posteriori Error Estimates for hp-FE Discretizations in Elastoplasticity
% -----------------------------------------------------------------------------------------

\chapter{Publication P3}

{\large
\begin{center}
A Posteriori Error Estimates for $hp$-FE Discretizations in Elastoplasticity
\end{center}
}
\vspace{1cm}
\textbf{Submitted to:} \emph{Computers \& Mathematics with Applications}. \\

Available as arXiv-preprint: \emph{arXiv:2401.09105} (2024).

% -----------------------------------------------------------------------------------------
%    Paper 4  -  An hp-adaptive strategy based on locally predicted error reductions
% -----------------------------------------------------------------------------------------

\chapter{Publication P4}

{\large
\begin{center}
An $hp$-adaptive strategy based on locally predicted error reductions
\end{center}
}
\vspace{1cm}
\textbf{Submitted to:} \emph{Computational Methods in Applied Mathematics} (in revision). \\

Available as arXiv-preprint \emph{arXiv:2311.13255} (2023).

\end{document}